\pdfoutput=1 

\documentclass{article}%
\usepackage{amssymb}
\usepackage{amsmath}
\usepackage{graphicx}
\usepackage{chicago}
\usepackage{amsfonts}%
\providecommand{\U}[1]{\protect\rule{.1in}{.1in}}

\newtheorem{theorem}{Theorem}

\newtheorem{corollary}{Corollary}

\newtheorem{example}[theorem]{Example}

\newtheorem{lemma}[theorem]{Lemma}

\newtheorem{proposition}{Proposition}

\begin{document}

\title{The Logic of Partitions: \\Introduction to the Dual of the Logic of Subsets}
\author{David Ellerman\\Department of Philosophy \\University of California/Riverside}
\maketitle

\begin{abstract}
Modern categorical logic as well as the Kripke and topological models of
intuitionistic logic suggest that the interpretation of ordinary
"propositional" logic should in general be the logic of subsets of a given
universe set. Partitions on a set are dual to subsets of a set in the sense of
the category-theoretic duality of epimorphisms and monomorphisms--which is
reflected in the duality between quotient objects and subobjects throughout
algebra. If "propositional" logic is thus seen as the logic of subsets of a
universe set, then the question naturally arises of a dual logic of partitions
on a universe set. This paper is an introduction to that logic of partitions
dual to classical subset logic. The paper goes from basic concepts up through
the correctness and completeness theorems for a tableau system of partition logic.

\end{abstract}
\tableofcontents

\section{Introduction to partition logic}

\subsection{The idea of a dual logic of partitions}

In ordinary propositional logic, the atomic variables and compound formulas
are usually interpreted as representing propositions. But in terms of
mathematical entities, the variables and formulas may be taken as representing
subsets of some fixed universe set $U$ with the propositional interpretation
being isomorphic to the special case of a one element set $U$ with subsets $0$
and $1$. The logic of subsets of a universe set is used to model
truth-functional reasoning with propositions as well as other binary on-off
phenomena such as switching circuits. But these specific applications should
not obscure the fact that from the purely mathematical viewpoint, it is the
logic of subsets of a universe set.

\begin{quote}
The propositional calculus considers "Propositions" $p$, $q$, $r$,... combined
under the operations "and","or", "implies", and "not", often written as
$p\wedge q$, $p\vee q$, $p\Rightarrow q$, and $\lnot p$. Alternatively, if
$P$, $Q$, $R$,... are subsets of some fixed set $U$ with elements $u$, each
proposition $p$ may be replaced by the proposition $u\in P $ for some subset
$P\subset U$; the propositional connectives above then become operations on
subsets; intersection $\wedge$, union $\vee$, implication ($P\Rightarrow Q$ is
$\lnot P\vee Q$), and complement of subsets. \cite[p. 48]{macm:sh}
\end{quote}

\noindent In order not to take a specific application for the whole theory,
"propositional logic" would be better called "subset logic" or "Boolean logic."

Modern logic was reformulated and generalized largely by F. William Lawvere in
what is now called \textit{categorical logic}.\footnote{See Lawvere and
Rosebrugh \citeyear[Appendix A]{law:sfm}\ for a good treatment. For the
generalization to topos theory see Mac Lane and Moerdijk
\citeyear{macm:sh}
and for the category theoretic background, the best references are Mac Lane
\citeyear{mac:cwm}
and Awodey
\citeyear{awod:ct}%
.} Subsets were generalized to subobjects or "parts"\ (equivalence classes of
monomorphisms) so that this generalized logic became the logic of subobjects
or parts in a topos (such as the category of sets).\footnote{Sometimes the
propositional and subset interpretations are "connected" by interpreting $U$
as the set of possible worlds and a subset as the set of possible worlds where
a proposition is true. While this interpretation may be pedagogically useful,
it is conceptually misleading since $U$ is simply an abstract set. None of the
philosophical problems involved in "possible worlds" semantics have anything
to do with the subset interpretation of ordinary logic. Starting with the
subset interpretation, each subset $P\subseteq U$ has an associated
proposition $u\in P$ (see the tableau rules given below) so that the
operations on the subsets (e.g., union, intersection, etc.) have corresponding
truth-functional operations on these propositions.}

There is a category-theoretic duality between subsets of a set and partitions
on a set. "The dual notion (obtained by reversing the arrows) of `part' is the
notion of \textit{partition}."\ \cite[p. 85]{law:sfm} In category theory, this
reverse-the-arrows duality gives the duality between monomorphisms, e.g.,
injective set functions, and epimorphisms, e.g., surjective set functions, and
between subobjects and quotient objects. If modern logic is formulated as the
logic of subsets (or more generally, subobjects or "parts"), then the question
naturally arises of a dual logic of partitions.

Quite aside from category theory duality, Gian-Carlo Rota emphasized the
seminal analogies between the subsets of a set and the partitions on a set.
Just as subsets of a set are partially ordered by inclusion, so partitions on
a set are partially ordered by refinement. Moreover, both partial orderings
are lattices (i.e., have meets and joins) with a top element and a bottom
element. This work on partition logic was inspired by both Rota's program to
develop the subset-partition analogies and by the category-theoretic treatment
of logic together with the reverse-the-arrows duality between subsets and partitions.

This paper is an introduction to the "propositional" part of partition logic.
The first part is an introduction that defines the logical operations on
partitions and develops the algebra of partitions (along with the dual algebra
of equivalence relations) that is analogous to the Boolean algebra of subsets
of a set. The second part of the paper develops a proof theory for partition
logic using semantic tableaus and proves the correctness and completeness
theorems for those tableaus.

\subsection{Duality of elements of a subset and distinctions of a partition}

The set-of-blocks definition of a \textit{partition} on a set $U$ is a set of
non-empty subsets ("blocks") of $U$ where the blocks are mutually exclusive
(the intersection of distinct blocks is empty) and jointly exhaustive (the
union of the blocks is $U$). If subsets are dual to partitions (in the sense
of monomorphisms being dual to epimorphisms), then what is the dual concept
that corresponds to the notion of \textit{elements of a subset}? The dual
notion is the notion of a \textit{distinction of a partition} which is a pair
of elements in distinct blocks of the partition. The duality between elements
of a subset and distinctions of a partition already appears in the very notion
of a function between sets. What binary relations, i.e., subsets $R\subseteq
X\times Y$, specify functions $f:X\rightarrow Y$?

A binary relation $R\subseteq X\times Y$ \textit{transmits elements} if for
each element $x\in X$, there is an ordered pair $\left(  x,y\right)  \in R$
for some $y\in Y$.

A binary relation $R\subseteq X\times Y$ \textit{reflects elements} if for
each element $y\in Y$, there is an ordered pair $\left(  x,y\right)  \in R$
for some $x\in X$.

A binary relation $R\subseteq X\times Y$ \textit{transmits distinctions} if
for any pairs $\left(  x,y\right)  $ and $\left(  x^{\prime},y^{\prime
}\right)  $ in $R$, if $x\not =x^{\prime}$, then $y\not =y^{\prime}$.

A binary relation $R\subseteq X\times Y$ \textit{reflects distinctions} if for
any pairs $\left(  x,y\right)  $ and $\left(  x^{\prime},y^{\prime}\right)  $
in $R$, if $y\not =y^{\prime}$, then $x\not =x^{\prime}$.

The dual role of elements and distinctions can be seen if we translate the
usual characterization of the binary relations that define functions into the
elements-and-distinctions language. A binary relation $R\subseteq X\times Y$
defines a \textit{function} $X\rightarrow Y$ if it is defined everywhere on
$X$ and is single-valued. But "being defined everywhere" is the same as
transmitting elements, and being single-valued is the same as reflecting distinctions:

\begin{center}
a binary relation $R$ is a \textit{function} if it transmits elements and
reflects distinctions.
\end{center}

What about the other two special types of relations, i.e., those which
transmit distinctions or reflect elements? The two important special types of
functions are the injections and surjections, and they are defined by the
other two notions:

\begin{center}
an \textit{injective function} is a function that transmits distinctions, and

a \textit{surjective function} is a function that reflects elements.
\end{center}

In view of the dual role of subsets and partitions (and elements and
distinctions), it is interesting to note that many basic ideas expressed using
subsets such as the notion of a "function" could just as well be expressed in
a dual manner using partitions. The dual to the product $X\times Y$ is the
coproduct $X%
{\textstyle\biguplus}
Y$ which in the category of sets is the disjoint union. If a binary relation
is defined as a subset $R$ of the product $X\times Y$, then a \textit{binary
corelation} would be a partition $\pi$ on the coproduct $X%
{\textstyle\biguplus}
Y$. Instead of defining a function as a certain type of binary relation (i.e.,
which transmits elements and reflects distinctions), a function could just as
well be defined as a certain type of binary corelation. Let $\left[  u\right]
_{\pi}$ denote the block of a partition $\pi$ containing an element $u$ from
the universe set of the partition. Then a binary corelation $\pi$ (a partition
on $X%
{\textstyle\biguplus}
Y$) is \textit{functional} if 1) every element $x\in X$ is transmitted to some
$y $-block, i.e., $\exists y\in Y$, $x\in\left[  y\right]  _{\pi}$, and 2)
distinctions on $Y$ are reflected as distinctions of $\pi$, i.e., if
$y\not =y^{\prime}$ for $y,y^{\prime}\in Y$, then $\left[  y\right]  _{\pi
}\not =\left[  y^{\prime}\right]  _{\pi}$.

Moreover, this definition of a function is quite familiar (with different
terminology) in combinatorics. For a functional corelation $\pi$, there is one
and only one block of the partition for each element $y\in Y$ so the blocks
$\left[  y\right]  _{\pi}$ can be thought of as "boxes." Then the elements of
$X$ can be thought of as "balls" and then a function is just a distribution of
the balls into the boxes. Thus the functional corelation definition of a
function is just a "disguised" version of the balls-in-boxes definition of a
function used in combinatorial theory \cite[p. 31]{stan:encomb1}. A functional
corelation is \textit{injective} if distinctions between balls are transmitted
as distinctions between boxes ("different balls to different boxes"), i.e.,
$x\not =x^{\prime}$ implies $\left[  x\right]  _{\pi}\not =\left[  x^{\prime
}\right]  _{\pi}$, and is \textit{surjective} if each box contains at least
one ball (i.e., each $y$ is reflected as an $x$). Although functions were
historically defined as functional binary relations, from the mathematical
viewpoint, functions could just as well be defined as functional binary corelations.

The duality between the two definitions of functions is clear in category
theory. Given the diagram $f:X\rightarrow Y$ in the category of sets, its
limit is the functional relation corresponding to $f$ and its colimit is the
functional corelation corresponding to $f$. The functional relation
corresponding to a function is its \textit{graph} and the functional
corelation corresponding to a function is its \textit{cograph} \cite[p.
29]{law:sfm}.

\subsection{Partitions and equivalence relations}

An \textit{equivalence relation} on a set $U$ is a subset $E\subseteq U\times
U$ that is reflexive, symmetric, and transitive. Every equivalence relation on
a set $U$ determines a partition on $U$ where the equivalence classes are the
mutually exclusive and jointly exhaustive blocks of the partition. Conversely,
every partition on a set determines an equivalence relation on the set (two
elements are equivalent if they are in the same block of the partition). The
notions of a partition on a set and an equivalence relation on a set are thus
interdefinable ("cryptomorphic"\ as Gian-Carlo Rota would say). Indeed,
equivalence relations and partitions are often considered as the "same." But
for our purposes it is important to keep the notions distinct (as in the above
definitions) so that we may consider the complementary type of binary
relation. A \textit{partition relation }$R\subseteq U\times U$\textit{\ }is
irreflexive (i.e., $\left(  u,u\right)  \not \in R$ for any $u\in U$),
symmetric [i.e., $\left(  u,u^{\prime}\right)  \in R$ implies $\left(
u^{\prime},u\right)  \in R$], and \textit{anti-transitive}\ in the sense that
if $\left(  u,u^{\prime}\right)  \in R$, then for any $a\in U$, either
$\left(  u,a\right)  \in R$ or $\left(  a,u^{\prime}\right)  \in R$ [i.e.,
$U\times U-R=R^{c}$ is transitive]. Thus as binary relations, equivalence
relations and partition relations are complementary. That is, $E\subseteq
U\times U$ is an equivalence relation if and only if (iff) $E^{c}\subseteq
U\times U$ is a partition relation. A partition relation is the set of
distinctions of a partition.

In a similar manner, the closed and open sets of a topological space can each
be defined in terms of the other and are complementary as subsets of the
space. Indeed, this is a useful analogy. There is a natural ("built-in")
closure operation on $U\times U=U^{2}$ which makes it a closure space. A
subset $C\subseteq U^{2}$ is \textit{closed} (1) if $C$ contains the diagonal
$\Delta=\left\{  \left(  u,u\right)  \mid u\in U\right\}  $ (reflexivity), (2)
if $\left(  u,u^{\prime}\right)  \in C$, then $\left(  u^{\prime},u\right)
\in C$ (symmetry), and (3) if $\left(  u,u^{\prime}\right)  $ and $\left(
u^{\prime},u^{\prime\prime}\right)  $ are in $C$, then $\left(  u,u^{\prime
\prime}\right)  $ is in $C$ (transitivity). Thus the closed sets of $U^{2}$
are the reflexive, symmetric, and transitive relations, i.e., the equivalence
relations on $U$. The intersection of any number of closed sets is closed.
Given a subset $S\subseteq U^{2}$, the \textit{closure} $\overline{S}$ is the
reflexive, symmetric, and transitive closure of $S$. The formation of the
closure $\overline{S}$ can be divided into two steps. First $S^{\ast}$ is
formed from $S$ by adding any diagonal pairs $\left(  u,u\right)  $ not
already in $S$ and by symmetrizing $S$, i.e., adding $\left(  u^{\prime
},u\right)  $ if $\left(  u,u^{\prime}\right)  \in S$. To form the transitive
closure of $S^{\ast}$, for any finite sequence $u=u_{1},u_{2},...,u_{n}%
=u^{\prime}$ with $\left(  u_{i},u_{i+1}\right)  \in S^{\ast}$ for
$i=1,...,n-1$, add $\left(  u,u^{\prime}\right)  $ and $\left(  u^{\prime
},u\right)  $ to the closure. The result is the reflexive, symmetric, and
transitive closure $\overline{S}$ of $S$. The complements of the closed sets
in $U\times U$ are defined as the \textit{open} sets which are the partition
relations on $U$. As usual, the \textit{interior} $\operatorname*{int}(S)$ of
any subset $S$ is defined as the complement of the closure of its complement:
$\operatorname*{int}(S)=\left(  \overline{S^{c}}\right)  ^{c}$.

It should, however, be carefully noted that the closure space $U\times U$ is
not a topological space, i.e., the closure operation on $U^{2}$ is not a
\textit{topological} closure operation in the sense that the union of two
closed set is not necessarily closed (or, equivalently, the intersection of
two open sets is not necessarily open). Since the lattice of open sets (or of
closed sets) of a topological space is distributive, this failure of the
closure operation on $U\times U$ to be topological is behind the
non-distributivity of the lattice of partitions (or of equivalence relations)
on a set $U$.

The set-of-blocks definition of a \textit{partition} $\pi$ on a set $U$ is a
set $\left\{  B\right\}  _{B\in\pi}$ of non-empty subsets or
"blocks"\ $B\subseteq U$ that are disjoint and whose union is $U$%
.\footnote{Just as the usual treatment of the Boolean algebra of all subsets
of a universe $U$ assumes that $U$ has one or more elements, so our treatment
of the lattice of all partitions on $U$ will assume that $U$ has \textit{two}
or more elements. This avoids the "degenerate" special cases of there being
only one subset of an empty $U$ and only one partition on a singleton $U$.} A
pair $\left(  u,u^{\prime}\right)  \in U\times U$ is a \textit{distinction} or
\textit{dit} (from DIsTinction) of the partition $\pi$ if there are distinct
blocks $B,B^{\prime}\in\pi$ with $u\in B$ and $u^{\prime}\in B^{\prime}$. The
set of distinctions of a partition $\pi,$ its \textit{dit set} denoted
$\operatorname*{dit}\left(  \pi\right)  \subseteq U\times U$, is the partition
seen as a partition relation:

\begin{center}
$\operatorname*{dit}\left(  \pi\right)  =\bigcup\limits_{B,B^{\prime}\in
\pi,B\not =B^{\prime}}B\times B^{\prime}$
\end{center}

\noindent(where it is understood that the union includes both the cartesian
products $B\times B^{\prime}$ and $B^{\prime}\times B$ for $B\not =B^{\prime}%
$).\footnote{Strictly speaking, one could argue that a "distinction" should be
an unordered pair $\left\{  u,u^{\prime}\right\}  $ but it is analytically
more convenient to deal with ordered pairs. In finite probability theory with
equiprobable elements in the sample space, the relative count of elements in a
subset (or event) defines the \textit{probability} $\operatorname*{Prob}%
\left(  S\right)  $ of the subset $S$. Dualizing, the count of the
distinctions of a partition relative to the total number of ordered pairs with
a finite universe $U$ defines the \textit{logical entropy} $h(\pi)$ of a
partition $\pi$ \cite{ell:cd}. In this logical information theory, it is also
analytically better to deal with ordered pairs. Then the logical entropy
$h\left(  \pi\right)  $ of a partition $\pi$ is simply the probability that a
random draw of a pair (with replacement) is a distinction of the partition
just as $\operatorname*{Prob}\left(  S\right)  $ is the probability that a
random draw is an element of the subset.}

A pair $\left(  u,u^{\prime}\right)  \in U\times U$ is an
\textit{indistinction} or \textit{indit} (from INDIsTinction) of a partition
$\pi$ if $u$ and $u^{\prime}$ belong to the same block of $\pi$. The set of
indistinctions of a partition $\pi$, its \textit{indit set} denoted
$\operatorname*{indit}\left(  \pi\right)  =U\times U-\operatorname*{dit}%
\left(  \pi\right)  $, is the complementary equivalence relation:

\begin{center}
$\operatorname*{indit}\left(  \pi\right)  =\bigcup\limits_{B\in\pi}B\times
B=U\times U-\operatorname*{dit}\left(  \pi\right)  =\operatorname*{dit}\left(
\pi\right)  ^{c}$.
\end{center}

In terms of\ the closure space structure on $U\times U$, let $\mathcal{O}%
\left(  U\times U\right)  $ be the open sets (partition relations) which are
the dit sets $\operatorname*{dit}(\pi)$ of partitions while the complementary
closed sets (equivalence relations) are the indit sets $\operatorname*{indit}%
\left(  \pi\right)  $ of partitions.

Partitions on $U$ are partially ordered by the \textit{refinement} relation:
given two partitions $\pi=\left\{  B\right\}  _{B\in\pi}$ and $\sigma=\left\{
C\right\}  _{C\in\sigma}$,

\begin{center}
$\sigma\preceq\pi$ (read "$\pi$ refines $\sigma$" or "$\sigma$ is refined by
$\pi$") if for any block $B\in\pi$, there is a block $C\in\sigma$ with
$B\subseteq C$.\footnote{Note that the opposite partial order is called the
"refinement" ordering in the customary "upside down" treatment of the lattice
of partitions. Gian-Carlo Rota used to joke that it should be called the
"unrefinement" relation. Indeed, in a recent book on Rota-style combinatorial
theory, that relation is sensibly called "reverse refinement" \cite[p.
30]{kung:rota}. It could also be called the "coarsening" \cite[p. 38]{law:sfm}
relation.}
\end{center}

\noindent The equivalent definition using dit sets (i.e., partition relations)
is just inclusion:

\begin{center}
$\sigma\preceq\pi$ iff $\operatorname*{dit}\left(  \sigma\right)
\subseteq\operatorname*{dit}\left(  \pi\right)  $.
\end{center}

Partitions might be represented by surjections $U\rightarrow\pi$ and every
refinement relation $\sigma\preceq\pi$ is realized by the unique map
$\pi\rightarrow\sigma$ that takes each block $B\in\pi$ to the block
$C\in\sigma$ containing it. The refinement map makes the following triangle commute:

\begin{center}
$%
\begin{array}
[c]{ccc}%
U & \rightarrow & \pi\\
\downarrow & \swarrow & \\
\sigma &  &
\end{array}
$

Refinement as a map
\end{center}

\noindent and thus it gives a morphism in the ("coslice") category of sets
under $U$ \cite[p. 15]{awod:ct}.

The partial ordering of partitions on $U$ has a least element or bottom which
is the indiscrete partition $0=\left\{  U\right\}  $ (nicknamed the "blob")
with the null dit set $\operatorname*{dit}(0)=\emptyset$ (no distinctions).
The blob distinguishes nothing and is refined by all partitions on $U$. The
partial ordering also has a greatest element or top which is the discrete
partition $1=\left\{  \left\{  u\right\}  :u\in U\right\}  $ where all blocks
are singletons and whose dit set is all ordered pairs off the diagonal, i.e.,
$\operatorname*{dit}(1)=U\times U-\Delta$ where $\Delta=\left\{  \left(
u,u\right)  :u\in U\right\}  $. The discrete partition refines all partitions
on $U$.

In any partial order with a least element $0$, an element $\alpha$ is an
\textit{atom} in the partial ordering if there is no element between it and
the bottom $0$, i.e., if $0\leq\pi\leq\alpha$ implies $\pi=0$ or $\pi=\alpha$.
In the inclusion partial order of subsets of $U$, the atoms are the singleton
subsets. In the refinement partial order of partitions, the atomic partitions
are the binary partitions, the partitions with two blocks. Any partition less
refined than a partition $\pi$ must fuse two or more blocks of $\pi$. Hence
the binary partitions are the partitions so that any less refined partition
has to be the blob.

\subsection{Category-theoretic duality of subsets and partitions}

In addition to the basic monomorphism-epimorphism duality between subsets and
partitions, a set of dual relationships between subset and partition concepts
as well as between element and distinction concepts will be described in this
section using basic category-theoretic notions in the category of sets. This
duality in the category of sets extends beyond the basic reverse-the-arrows
duality that holds in all categories, and it underlies the duality between
subset logic and partition logic.

In the category of sets, the singleton $1$ might be thought of as the generic
element. We have seen that functions preserve (or transmit) elements and
reflect (or transmit in the backwards direction) distinctions. The basic
property of the generic element $1$ is that for every element $u\in U$, there
is a function $1\overset{u}{\rightarrow}U$ that transmits "elementness" from
the generic element to $u\in U$. The partition-dual to the generic element $1$
is $2=\left\{  0,1\right\}  $ which might be thought of as the generic
distinction. The basic property of the generic distinction $2$ is that for any
pair $u,u^{\prime}$ of distinct elements of $U$, there is a function
$\alpha:U\rightarrow2$ that reflects or backwards-transmits "distinctness"
from the generic distinction $2$ to the pair $u,u^{\prime}$.

Given two parallel functions $f,g:X\rightarrow Y$, if they are different,
$f\not =g$, then there is an element $x\in X$ such that the two functions
carry $x$ to a distinction $f\left(  x\right)  \not =g\left(  x\right)  $ of
$Y$. By the basic property of the generic element $1$, there is a function
$1\overset{x}{\rightarrow}X$ that transmits the generic element to that
element $x$. Thus the generic element $1$ is a \textit{separator} in the sense
that given two set functions $f,g:X\rightarrow Y$, if $f\not =g$, then
$\exists x:1\rightarrow X$ (an injection) such that $1\overset{x}{\rightarrow
}X\overset{f}{\rightarrow}Y\not =1\overset{x}{\rightarrow}%
X\overset{g}{\rightarrow}Y$. Dually, by the basic property of the generic
distinction, there is a function $\alpha:Y\rightarrow2$ that reflects the
generic distinction to the distinction $f\left(  x\right)  \not =g\left(
x\right)  $ of $Y$. Thus the generic distinction $2$ is a \textit{coseparator}
\cite[pp. 18-19]{law:sfm} in the sense that given two set functions
$f,g:X\rightarrow Y $, if $f\not =g$, then $\exists\alpha:Y\rightarrow2$ (a
surjection) such that $X\overset{f}{\rightarrow}Y\overset{\alpha}{\rightarrow
}2\not =X\overset{g}{\rightarrow}Y\overset{\alpha}{\rightarrow}2$.

Other dual roles of the generic element $1$ and generic distinction $2$ follow
from the dual basic properties. Consider the product of $X$ and $Y$ in the
category of sets. A set $P$ with maps $p_{1}:P\rightarrow X$ and
$p_{2}:P\rightarrow Y$ is the \textit{product}, denoted $X\times Y$, if for
any set $Z$ and pair of maps $f:Z\rightarrow X$ and $g:Z\rightarrow Y$ with
domain $Z$, there is a unique map $\left\langle f,g\right\rangle :Z\rightarrow
P$ such that $p_{1}\left\langle f,g\right\rangle =g$ and $p_{2}\left\langle
f,g\right\rangle =g$. The generic element $1$ has the property that it
suffices as the test set $Z=1$. That is, if the set $P$ with its pair of maps
had the universal mapping property for pairs of maps with domain $1$, then it
has the universal mapping property for any pairs of maps with a common domain
$Z$, i.e., it is the product. This property of the generic element $1$ extends
to all limits in the category of sets.

The dual construction is the coproduct, denoted $X%
{\textstyle\biguplus}
Y$ or $X+Y$, which can be constructed as the disjoint union of $X$ and $Y$
with the two insertion maps. A set $C$ with maps $i_{1}:X\rightarrow C$ and
$i_{2}:Y\rightarrow C$ is the \textit{coproduct} $X%
{\textstyle\biguplus}
Y$ if for any set $Z$ and pair of maps $f:X\rightarrow Z$ and $g:Y\rightarrow
Z$ with codomain $Z$, there is a unique map (which we will denote)
$\left\rangle f,g\right\langle :C\rightarrow Z$ such that: $\left\rangle
f,g\right\langle i_{1}=f$ and $\left\rangle f,g\right\langle i_{2}%
=g$.\footnote{There seems to be no standard notation for the coproduct factor
map so we have just reversed the angle brackets from the product factor map.}
The generic distinction $2$ has the property that it suffices as the test set
$Z=2$. That is, if the set $C$ and its pair of maps had the universal mapping
property for pairs of maps with codomain $2$, then it has the universal
mapping property for any pair of maps with a common codomain $Z$, i.e., it is
the coproduct.\cite[p. 272]{law:cm} This property of the generic distinction
$2$ extends to all colimits in the category of sets.

The dual properties also show up in the respective partial orders (and
lattices). The images of injections $1\overset{u}{\longrightarrow}U$ are the
atoms $\left\{  u\right\}  $ in the inclusion partial order of subsets of $U$
and in the powerset Boolean algebra $\mathcal{P}(U)$. The inverse images of
surjections $U\overset{\alpha}{\longrightarrow}2$ are the atoms (binary
partitions) in the refinement partial order of partitions on $U$ and in the
partition lattice $\Pi(U)$ defined below.

Given a subset $S$ of $U$ and a partition $\pi$ on $U$, there is the
associated injection $S\longrightarrow U$ and the associated surjection
$U\longrightarrow\pi$ (taking $\pi$ as a set of blocks). The atom $\left\{
u\right\}  $ given by $1\overset{u}{\longrightarrow}U$ is contained in $S$,
iff $1\overset{u}{\longrightarrow}U$ uniquely factors through
$S\longrightarrow U$. Analogously, an atomic partition $U\overset{\alpha
}{\longrightarrow}2$ is refined by $\pi$ [$\operatorname*{dit}\left(
\alpha\right)  \subseteq\operatorname*{dit}\left(  \pi\right)  $] iff
$U\overset{\alpha}{\longrightarrow}2$ uniquely factors through
$U\longrightarrow\pi$.

\begin{center}
$%
\begin{array}
[c]{ccc}%
1 &  & \\
^{{}}\downarrow^{\exists!} & \searrow^{u} & \\
S & \longrightarrow & U
\end{array}
$ $\qquad%
\begin{array}
[c]{ccc}%
U & \longrightarrow & \pi\\
& \searrow^{\alpha} & ^{{}}\downarrow^{\exists!}\\
&  & 2
\end{array}
$

Analogous diagrams showing which atoms contained in an object (subset or partition)
\end{center}

The dual pullback and pushout constructions allow us to represent each
partition as a subset of a product and to represent each subset as a partition
on a coproduct.

Given a partition as a surjection $U\rightarrow\pi$, the pullback of the
surjection with itself, i.e., the \textit{kernel pair} \cite[p. 71]{mac:cwm}
of $U\rightarrow\pi$, gives the indit set $\operatorname*{indit}\left(
\pi\right)  $ as a subset of the product $U\times U$, i.e., as a binary
(equivalence) relation on $U$:

\begin{center}
$%
\begin{array}
[c]{ccc}%
\operatorname*{indit}\left(  \pi\right)  & \overset{p_{2}}{\longrightarrow} &
U\\
^{p_{1}}\downarrow^{{}} &  & \downarrow\\
U & \longrightarrow & \pi
\end{array}
$

Pullback for equivalence relation $\operatorname*{indit}\left(  \pi\right)  $.
\end{center}

Given a subset as an injection $S\rightarrow U$, the pushout of the injection
with itself, i.e., the \textit{cokernel pair} \cite[p. 66]{mac:cwm}\ of
$S\rightarrow U$, gives a partition $\Delta\left(  S\right)  $ on the
coproduct (disjoint union) $U%
{\textstyle\biguplus}
U$, i.e., a binary corelation which might be called a \textit{subset
corelation}:

\begin{center}
$%
\begin{array}
[c]{ccc}%
S & \longrightarrow & U\\
\downarrow &  & ^{{}}\downarrow^{\lbrack u^{\ast}]}\\
U & \overset{[u]}{\longrightarrow} & \Delta\left(  S\right)
\end{array}
$

Pushout for subset corelation $\Delta\left(  S\right)  $.
\end{center}

\noindent The disjoint union $U%
{\textstyle\biguplus}
U$ consists of the elements $u\in U $ and the copies $u^{\ast}$ of the
elements $u\in U$. The subset corelation $\Delta\left(  S\right)  $ is
constructed by identifying any $u$ and its copy $u^{\ast}$ for $u\in S$ so
$\Delta\left(  S\right)  $ is the partition on $U%
{\textstyle\biguplus}
U$ whose only non-singleton blocks are the pairs $\left\{  u,u^{\ast}\right\}
$ for $u\in S$.

The constructions can also be reversed by viewing the pullback square as a
pushout square, and by viewing the pushout square as a pullback square.
Equivalently, we can reconstruct $\pi$ as the coequalizer of the two
projection maps $p_{1},p_{2}$ from $\operatorname*{indit}\left(  \pi\right)
\subseteq U\times U$ to $U$ \cite[p. 89]{law:sfm}.

\begin{center}
$%
\begin{array}
[c]{ccccc}%
\operatorname*{indit}\left(  \pi\right)  & \underset{p_{2}}{\overset{p_{1}%
}{\rightrightarrows}} & U & \longrightarrow & \pi\\
&  &  & \searrow^{\alpha} & ^{{}}\downarrow^{\exists!}\\
&  &  &  & 2
\end{array}
$

Partition $\pi$ as coequalizer of $\operatorname*{indit}\left(  \pi\right)
\overset{p_{1}}{\longrightarrow}U$ and $\operatorname*{indit}\left(
\pi\right)  \overset{p_{2}}{\longrightarrow}U$.
\end{center}

Dually, we have the two maps $U\rightarrow\Delta\left(  S\right)  $ given by
$u\mapsto\left[  u\right]  _{\Delta\left(  S\right)  }$ and $u\mapsto\left[
u^{\ast}\right]  _{\Delta\left(  S\right)  }$, and the subset $S$ is
reconstructed as their equalizer:

\begin{center}
$%
\begin{array}
[c]{ccccc}%
1 &  &  &  & \\
^{\exists!}\downarrow^{{}} & \searrow^{u} &  &  & \\
S & \rightarrow & U & \rightrightarrows & \Delta\left(  S\right)
\end{array}
$

Subset $S$ as equalizer of $\left[  u\right]  :U\rightarrow\Delta\left(
S\right)  $ and $\left[  u^{\ast}\right]  :U\rightarrow\Delta\left(  S\right)
$.
\end{center}

In general, the equalizer (in the category of sets) of two set functions
$f,g:X\rightarrow Y$ is the largest subset $S$ of the domain $X$ so that no
element of $S$ goes via the functions to a distinction $\left(  f\left(
x\right)  ,g\left(  x\right)  \right)  $ of the codomain $Y$.

Dually, the coequalizer of two set functions $f,g:X\rightarrow Y$ is the
largest (most refined) partition $\pi$ on the codomain $Y$ so that no
distinction of $\pi$ comes via the functions from an element of the domain $X
$ (i.e., has the form $\left(  f\left(  x\right)  ,g\left(  x\right)  \right)
$ for some $x\in X$).

Then the functions $\left[  u\right]  ,\left[  u^{\ast}\right]
:U\rightrightarrows\Delta\left(  S\right)  $ are such that $S$ is the largest
subset of the domain $U$ so that no element of the subset goes via those
functions to a distinction of the codomain $\Delta\left(  S\right)  $.

The functions $p_{1},p_{2}:\operatorname*{indit}\left(  \pi\right)
\rightrightarrows U$ are such that $\pi$ is the largest partition on the
codomain $U$ so that no distinction of the partition comes via those functions
from an element of the domain $\operatorname*{indit}\left(  \pi\right)  $
\cite[p. 89]{law:sfm}.

\begin{center}%
\begin{tabular}
[c]{l|l|l|}\cline{2-3}%
\textit{Dualities} & Subsets & Partitions\\\hline\hline
\multicolumn{1}{|l|}{Generic element} & Generic element $1$ & Generic
distinction $2$\\\hline
\multicolumn{1}{|l|}{Basic property} & Each element $u\in U$ & Each
distinction $u\not =u^{\prime}$\\
\multicolumn{1}{|l|}{of generic element} & realized by some $1\rightarrow U$ &
realized by some $U\rightarrow2$\\\hline
\multicolumn{1}{|l|}{Separating functions} & $1$ is a separator & $2$ is a
coseparator\\\hline
\multicolumn{1}{|l|}{Sufficient test set} & $1$ is a test set for limits & $2
$ is a test set for colimits\\\hline
\multicolumn{1}{|l|}{Objects} & Subsets: monos $S\longrightarrow U$ &
Partitions: epis $U\longrightarrow\pi$\\\hline
\multicolumn{1}{|l|}{Atoms in partial orders} & Images of monos
$1\overset{u}{\longrightarrow}U$ & Inv. images of epis $U\overset{\alpha
}{\longrightarrow}2$\\\hline
\multicolumn{1}{|l|}{Inclusion of atoms} & $1\overset{u}{\longrightarrow}U$
uniquely factors & $U\overset{\alpha}{\longrightarrow}2$ uniquely factors\\
\multicolumn{1}{|l|}{} & through $S\longrightarrow U$ & through
$U\longrightarrow\pi$\\\hline
\multicolumn{1}{|l|}{Subsets $\leftrightarrow$ Partitions} & Partition
$\Delta\left(  S\right)  $ on $U%
{\textstyle\biguplus}
U$ & Subset $\operatorname*{indit}\left(  \pi\right)  $ of $U\times U$\\
\multicolumn{1}{|l|}{} & is cokernal pair of $S\rightarrow U$ & is kernel pair
of $U\rightarrow\pi$\\\hline
\multicolumn{1}{|l|}{Inverse operation} & Subset $S$ is equalizer & Partition
$\pi$ is coequalizer\\
\multicolumn{1}{|l|}{} & of $\left[  u\right]  ,\left[  u^{\ast}\right]
:U\rightrightarrows\Delta\left(  S\right)  $ & of $p_{1},p_{2}%
:\operatorname*{indit}\left(  \pi\right)  \rightrightarrows U$\\\hline
\end{tabular}

Summary of dual relationships
\end{center}

\subsection{Lattice of partitions}

Traditionally the "lattice of partitions," e.g., \cite{Birk:lt} or
\cite{grat:glt}, was defined as isomorphic to the lattice of equivalence
relations where the partial order was inclusion between the equivalence
relations as subsets of $U\times U$. But since equivalence relations and
partition relations are complementary subsets of the closure space $U\times U
$, we have two anti-isomorphic lattices with opposite partial orders.

Which lattice should be used in partition logic? For the purposes of comparing
formulas with ordinary logic (interpreted as applying to subsets of elements),
it is crucial to take the lattice of partitions as (isomorphic to) the lattice
$\mathcal{O}\left(  U\times U\right)  $ of partition relations (sets of
distinctions), the opposite of the lattice of equivalence relations.

The \textit{lattice of partitions} $\Pi(U)$ on $U$ adds the operations of join
and meet to the partial ordering of partitions on $U$ with the top $1$ and the
bottom $0$.\footnote{For a survey of what is known about partition lattices,
see \cite{grat:glt} where the usual opposite presentation is used.} There are
at least four ways that partitions and operations on partitions might be defined:

\begin{enumerate}
\item the basic set-of-blocks definition of partitions and their operations;

\item the closure space approach using open subsets or dit sets and the
interior operator on $U\times U$;

\item the graph-theoretic approach where the blocks of a partition on $U$ are
the nodes in the connected components of a simple (at most one arc between two
nodes and no loops at a node) undirected graph;\footnote{See any introduction
to graph theory such as Wilson
\citeyear{wil:gt}
for the basic notions.} and

\item the approach where the blocks of a partition on $U$ are the atoms of a
complete Boolean subalgebra of the powerset Boolean algebra $\mathcal{P}(U)$
of subsets of $U$ \cite{ore:ter}.
\end{enumerate}

The lattice of partitions $\Pi(U)$ is the partition analogue of the powerset
Boolean lattice $\mathcal{P}(U)$. In the powerset lattice, the partial order
is inclusion of elements, and in the partition lattice, it is inclusion of distinctions.

The join $\pi\vee\sigma$ in $\Pi(U)$ is the partition whose blocks are the
non-empty intersections $B\cap C$ of the blocks of the two partitions. The
equivalent dit-set definition in $\mathcal{O}\left(  U\times U\right)  $ is
simply the union: $\operatorname*{dit}\left(  \pi\vee\sigma\right)
=\operatorname*{dit}\left(  \pi\right)  \cup\operatorname*{dit}\left(
\sigma\right)  $.

Recall that the closure operator on the closure space was not topological in
the sense that the union of two closed sets is not necessarily closed and thus
the intersection of two open sets (i.e., two dit sets) is not necessarily
open. Hence the definition of the meet of two partitions requires some more
complication. The dit-set definition in $\mathcal{O}\left(  U\times U\right)
$ is the easiest: the dit set of the \textit{meet} of two partitions is the
interior of the intersection of the two dit sets, i.e.,

\begin{center}
$\operatorname*{dit}\left(  \sigma\wedge\pi\right)  =\operatorname*{int}%
\left(  \operatorname*{dit}\left(  \sigma\right)  \cap\operatorname*{dit}%
\left(  \pi\right)  \right)  $.
\end{center}

\noindent In the older literature, this meet of two partitions is what is
defined as the join\ of the two equivalence relations. Given the two
partitions as sets of blocks $\left\{  B\right\}  _{B\in\pi}$ and $\left\{
C\right\}  _{C\in\sigma}$ in $\Pi(U)$, two elements $u$ and $u^{\prime}$ are
directly equated, $u\sim u^{\prime}$ if $u$ and $u^{\prime}$ are in the same
block of $\pi$ or $\sigma$ so the set of directly equated pairs is:
$\operatorname*{indit}\left(  \sigma\right)  \cup\operatorname*{indit}\left(
\pi\right)  $. Then $u$ and $u^{\ast}$ are in the same block of the join in
$\Pi(U)$ if there is a finite sequence $u=u_{1}\sim u_{2}\sim...\sim
u_{n}=u^{\ast}$ that indirectly equates $u$ and $u^{\ast}$. The operation of
indirectly equating two elements is just the closure operation in the closure
space so the set of pairs indirectly equated, i.e., equated in the join
$\sigma\wedge\pi$ in $\Pi(U)$, is:

\begin{center}
$\operatorname*{indit}\left(  \sigma\wedge\pi\right)  =\overline{\left(
\operatorname*{indit}\left(  \sigma\right)  \cup\operatorname*{indit}\left(
\pi\right)  \right)  }$.
\end{center}

\noindent The complementary subset of $U\times U$ is the dit set of the meet
of the partitions in $\mathcal{O}\left(  U\times U\right)  $:

\begin{center}
$\operatorname*{dit}\left(  \sigma\wedge\pi\right)  =\operatorname*{indit}%
\left(  \sigma\wedge\pi\right)  ^{c}=\overline{\left(  \operatorname*{indit}%
\left(  \sigma\right)  \cup\operatorname*{indit}\left(  \pi\right)  \right)
}^{c}=\operatorname*{int}\left(  \operatorname*{dit}\left(  \sigma\right)
\cap\operatorname*{dit}\left(  \pi\right)  \right)  $.
\end{center}

This defines the lattice of partitions $\Pi(U)$ and the isomorphic lattice
$\mathcal{O}\left(  U\times U\right)  $ which represents the partitions as
open subsets of the product $U\times U$:

\begin{center}
\fbox{$\Pi(U)\cong\mathcal{O}\left(  U\times U\right)  $}.

Representation of the lattice of partitions $\Pi(U)$

as the lattice of open subsets $\mathcal{O}\left(  U\times U\right)  $.
\end{center}

The analogies between the lattice of subsets $\mathcal{P}(U)$ and the lattice
of partitions $\Pi(U)$ are summarized in the following table.

\begin{center}%
\begin{tabular}
[c]{l|l|l|}\cline{2-3}%
\textit{Analogies} & Boolean lattice of subsets & Lattice of
partitions\\\hline\hline
\multicolumn{1}{|l|}{Elements} & Elements of subsets & Distinctions of
partitions\\\hline
\multicolumn{1}{|l|}{Partial order} & Inclusion of elements & Inclusion of
distinctions\\\hline
\multicolumn{1}{|l|}{Join} & Elements of join are & Distinctions of join are\\
\multicolumn{1}{|l|}{} & union of elements & union of distinctions\\\hline
\multicolumn{1}{|l|}{Meet} & Largest subset & Largest partition\\
\multicolumn{1}{|l|}{} & of only common elements & of only common
distinctions\\\hline
\multicolumn{1}{|l|}{Top} & Subset $U$ with all elements & Partition $1$ with
all distinctions\\\hline
\multicolumn{1}{|l|}{Bottom} & Subset $\emptyset$ with no elements & Partition
$0$ with no distinctions\\\hline
\end{tabular}

Elements-distinctions analogies between the Boolean lattice of subsets and the
lattice of partitions
\end{center}

With this definition of the lattice of partitions $\Pi(U)$, the usual lattice
of equivalence relations is $\Pi(U)^{op}$ where the top is $\widehat{1}%
=U\times U=\operatorname*{indit}\left(  0\right)  $ and the bottom is
$\widehat{0}=\Delta=\operatorname*{indit}\left(  1\right)  $%
.\footnote{Inevitably notational conflicts arise for such common symbols as
"$0$" and "$1$" so where there is less risk of confusion, different uses of
these symbols will be clear from the context. In other cases, the symbols are
modified as in using $\widehat{1}$ and $\widehat{0}$ for the top and bottom of
the opposite lattice of equivalence relations.}

\subsection{Two other definitions of the partition meet operation}

Since the partition meet is the first non-trivial definition of a partition
operation, we might also give the equivalent definitions using the
graph-theoretic method and the complete-Boolean-subalgebras method.

The power of the dit-set approach to defining partition operations is that it
allows us to mimic subset operations using dit sets and the interior
operations as needed. The power of the graph-theoretic approach is that it
allows a very intuitive connection back to the truth tables of classical
propositional logic. The truth tables for the classical Boolean propositional
connectives can be stated in an abbreviated form using signed formulas such as
$T\left(  \pi\wedge\sigma\right)  $ or $F\sigma$. The truth table for the
Boolean meet $\pi\wedge\sigma$ is abbreviated by saying the \textit{Boolean
conditions} for $T\left(  \pi\wedge\sigma\right)  $ are "$T\pi$ and $T\sigma$"
while the \textit{Boolean conditions} for $F\left(  \pi\wedge\sigma\right)  $
are "$F\pi$ or $F\sigma$". Thus for the four Boolean operations of join
$\pi\vee\sigma$, meet $\pi\wedge\sigma$, implication $\sigma\Rightarrow\pi$,
and Sheffer stroke, not-and or nand $\sigma\mid\pi$, the table of Boolean
conditions is as follows:

\begin{center}%
\begin{tabular}
[c]{|c||c|c||c|c||}\hline
Signed Formula & $T\left(  \pi\vee\sigma\right)  $ & $F\left(  \pi\vee
\sigma\right)  $ & $T\left(  \sigma\Rightarrow\pi\right)  $ & $F\left(
\sigma\Rightarrow\pi\right)  $\\\hline
Boolean Conditions & $T\pi$ or $T\sigma$ & $F\pi$ and $F\sigma$ & $F\sigma$ or
$T\pi$ & $T\sigma$ and $F\pi$\\\hline
\end{tabular}

Boolean conditions for $\vee$ and $\Rightarrow$,

and%

\begin{tabular}
[c]{|c||c|c||c|c||}\hline
Signed Formula & $T\left(  \pi\wedge\sigma\right)  $ & $F\left(  \pi
\wedge\sigma\right)  $ & $T\left(  \sigma\mid\pi\right)  $ & $F\left(
\sigma\mid\pi\right)  $\\\hline
Boolean Conditions & $T\pi$ and $T\sigma$ & $F\pi$ or $F\sigma$ & $F\sigma$ or
$F\pi$ & $T\sigma$ and $T\pi$\\\hline
\end{tabular}

Boolean conditions for $\wedge$ and $\mid$.
\end{center}

Given any partition $\pi$ on $U$, and any pair of elements $\left(
u,u^{\prime}\right)  $, we say that $T\pi$ holds at $\left(  u,u^{\prime
}\right)  $ if $\left(  u,u^{\prime}\right)  $ is a distinction of $\pi$, and
that $F\pi$ holds at $\left(  u,u^{\prime}\right)  $ if $\left(  u,u^{\prime
}\right)  $ is not a distinction of $\pi$, i.e., if $u$ and $u^{\prime}$ are
in the same block of $\pi$. Given any two partitions $\pi$ and $\sigma$ on
$U$, we can define the partition version of \textit{any} Boolean connective
$\pi\ast\sigma$ by putting an arc between any two nodes $u$ and $u^{\prime}$
if the Boolean conditions for $F\left(  \pi\ast\sigma\right)  $ hold at
$\left(  u,u^{\prime}\right)  $. Then the blocks of the partition operation
$\pi\ast\sigma$ are taken as the nodes in the connected components of that
graph. Thus two elements $u$ and $u^{\prime}$ are in the same block of the
partition $\pi\ast\sigma$ if there is a chain or finite sequence
$u=u_{1},u\,_{2},...,u_{n-1},u_{n}=u^{\prime}$ such that for each
$i=1,...,n-1$, the Boolean conditions for $F\left(  \pi\ast\sigma\right)  $
hold at $\left(  u_{i},u_{i+1}\right)  $.

In order for $\pi\ast\sigma$ to distinguish $u$ and $u^{\prime}$, it has to
"cut" them apart in the sense of the graph-theoretic notion of a "cut" which
is the graph-theoretic dual to the notion of a chain \cite[p. 31]{rock:nfmo}.
A set of arcs in a graph form a \textit{cut between the nodes }$u$%
\textit{\ and }$u^{\prime}$ if every chain connecting $u$ and $u^{\prime} $
contains an arc from the set--so that the set of arcs cut every chain
connecting the two points. The complementation-duality between chains and cuts
is brought out by the fact that if we arbitrarily color the arcs of any simple
undirected graph by either black or white, then for any two nodes, there is
either a white cut between the nodes or a black chain connecting the nodes.
The above graph-theoretic definition of $\pi\ast\sigma$, i.e., two points are
not distinguished if there is "falsifying" chain connecting the points with
the Boolean conditions for $F\left(  \pi\ast\sigma\right)  $ holding at each
arc (i.e., a black chain), can be stated in an equivalent dual form. Two
points are distinguished in $\pi\ast\sigma$ if the set of arcs where the
Boolean conditions for $T\left(  \pi\ast\sigma\right)  $ hold form a (white)
cut between the two points. These \textit{falsifying chains} and
\textit{distinguishing cuts} results are the basis for the semantic tableaus
used below.

This graph-theoretic approach can be used to uniformly define all the
partition logical operations in terms of the corresponding Boolean logical
operations, but the case at hand is the meet. The graph constructed for the
meet would have an arc between $u$ and $u^{\prime}$ if the Boolean conditions
for $F\left(  \pi\wedge\sigma\right)  $ held at $\left(  u,u^{\prime}\right)
$, i.e., if $F\pi$ or $F\sigma$ held at $\left(  u,u^{\prime}\right)  $. But
this just means that $\left(  u,u^{\prime}\right)  \in\operatorname*{indit}%
\left(  \sigma\right)  \cup\operatorname*{indit}\left(  \pi\right)  $, and the
nodes in the connected components of that graph are the nodes $u$ and
$u^{\prime}$ connected by a finite sequence $u=u_{1},u\,_{2},...,u_{n-1}%
,u_{n}=u^{\prime}$ where for each $i=1,...,n-1$, $\left(  u_{i},u_{i+1}%
\right)  \in\operatorname*{indit}\left(  \sigma\right)  \cup
\operatorname*{indit}\left(  \pi\right)  $, which is the closure space
definition of the meet given above.

\begin{example}
Let $\sigma=\left\{  \left\{  a,b,c\right\}  ,\left\{  d,e\right\}  \right\}
$ and $\pi=\left\{  \left\{  a,b\right\}  ,\left\{  c,d,e\right\}  \right\}
$. In the graph below, all the arcs in the complete graph $K_{5}$ on five
nodes are labelled according to the status of the two endpoints in the two
partitions. The Boolean conditions for $F\left(  \sigma\wedge\pi\right)  $ are
"$F\sigma$ or $F\pi$" . The arcs where those conditions hold are the solid
lines. In the graph with only the solid arcs, there is only one connected
component so $\sigma\wedge\pi=\left\{  \left\{  a,b,c,d,e\right\}  \right\}
=0$. Equivalently, the set of arcs where the Boolean conditions for $T\left(
\sigma\wedge\pi\right)  $ hold, i.e., the dashed arcs, do not "cut" apart any
pair of points.
\end{example}

%

\begin{center}
\includegraphics[
natheight=182.005005bp,
natwidth=317.994995bp,
height=184.6875pt,
width=321.1875pt
]%
{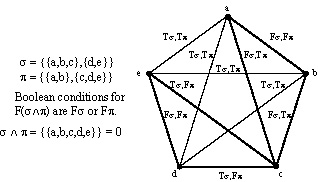}%
\\
Figure 1: Graph for meet $\sigma\wedge\pi$%
\end{center}

For the Boolean subalgebra approach, given a partition $\pi$ on $U$, define
$\mathcal{B}\left(  \pi\right)  \subseteq\mathcal{P}(U)$ as the complete
subalgebra generated by the blocks of $\pi$ as the atoms so that all the
elements of $\mathcal{B}\left(  \pi\right)  $ are formed as the arbitrary
unions and intersections of blocks of $\pi$. Conversely, given any complete
subalgebra $\mathcal{B}$ of $\mathcal{P}(U)$, the intersection of all elements
of $\mathcal{B}$ containing an element $u\in U$ will provide the atoms of
$\mathcal{B}$ which are the blocks in a partition $\pi$ on $U$ so that
$\mathcal{B}=\mathcal{B}\left(  \pi\right)  $. Thus an operation on complete
subalgebras of the powerset Boolean algebra will define a partition operation.
Since the blocks of the partition meet $\pi\wedge\sigma$ are minimal under the
property of being the exact union of $\pi$-blocks and also the exact union of
$\sigma$-blocks, a nice feature of this approach to partitions is that:

\begin{center}
$\mathcal{B}\left(  \pi\wedge\sigma\right)  =\mathcal{B}\left(  \pi\right)
\cap\mathcal{B}\left(  \sigma\right)  $.
\end{center}

The powerset Boolean algebra (BA) $\mathcal{P}(U)$ is not just a lattice; it
has additional structure which can be defined using the binary connective of
the set implication: $A\Rightarrow B=\left(  U-A\right)  \cup B=A^{c}\cup B$,
for $A,B\subseteq U$. The lattice structure on $\Pi(U)$ needs to be enriched
with other operations such as the new binary operation of implication on partitions.

\subsection{Partition implication operation}

Boolean algebras, or more generally, Heyting algebras are not just lattices;
there is another operation $A\Rightarrow B$, the implication operation. In a
Heyting algebra, the implication can be introduced by an adjunction (treating
the partial order as the morphisms in a category) that can be written in the
Gentzen style\footnote{Sometimes the Gentzen-style statement $%
\begin{tabular}
[c]{c}%
$x\rightarrow Gy$\\\hline
$Fx\rightarrow y$%
\end{tabular}
\ $ of an adjunction, $\operatorname*{Hom}_{Y}\left(  Fx,y\right)
\cong\operatorname*{Hom}_{X}(x,Gy)$, has the top and bottom reversed. But
there is a theory showing how adjoints arise out of representations of
heteromorphisms \cite{ell:af}, and that theory suggests that the Gentzen-style
statement should be written as above since there are "behind the scenes"
heteromorphisms (dashed arrows) as vertical downward maps $Gy\dashrightarrow
y$ and $x\dashrightarrow Fx$ so that this "adjoint square" commutes, i.e.,
$x\rightarrow Gy\dashrightarrow y=x\dashrightarrow Fx\rightarrow y$.} which in
this case is an "if and only if" statement:

\begin{center}
$%
\begin{tabular}
[c]{c}%
$C\leq A\Rightarrow B$\\\hline
$C\wedge A\leq B$%
\end{tabular}
$

Implication as the right adjoint to meet in a Heyting algebra.
\end{center}

\noindent In the standard model of the Heyting algebra of open sets of a
topological space, the implication is defined for open sets $A$ and $B$ as:

\begin{center}
$A\Rightarrow B=\operatorname*{int}(A^{c}\cup B)$.
\end{center}

A co-Heyting algebra is also a lattice with top and bottom but with the dual
adjunction where the \textit{difference} $B^{c}\backslash A^{c}$ is left
adjoint to the join:

\begin{center}
$%
\begin{tabular}
[c]{c}%
$B^{c}\leq A^{c}\vee C^{c}$\\\hline
$B^{c}\backslash A^{c}\leq C^{c}$%
\end{tabular}
$

Difference as the left adjoint to join in a co-Heyting algebra.
\end{center}

\noindent In the standard model of the co-Heyting algebra of closed sets of a
topological space, the difference is defined for closed sets $A^{c}$ and
$B^{c}$ (where $A$ and $B$ are open sets) as:

\begin{center}
$B^{c}\backslash A^{c}=\overline{\left(  B^{c}\cap A^{cc}\right)  }%
=\overline{\left(  B^{c}\cap A\right)  }=\left(  A\Rightarrow B\right)  ^{c}$.
\end{center}

Neither of these adjunctions holds in the lattice of partitions $\Pi(U)$ (or
its opposite). The adjunctions imply distributivity for Heyting and co-Heyting
algebras, and lattices of partitions (usually viewed in the opposite
presentation as the lattice of equivalence relations) are standard examples of
non-distributive lattices.

How might the implication partition $\sigma\Rightarrow\pi$ of two partitions
(or the difference between two equivalence relations) be defined? Some
motivation might be extracted from Heyting algebras, or, equivalently,
intuitionistic propositional logic. The subset version of intuitionistic
propositional logic is explicit in its topological interpretation where the
variables are interpreted as open subsets of a topological space $U$ and the
valid formulas are those that evaluate to the whole space $U$ regardless of
what open subsets are assigned to the atomic variables. The implication is
then defined as: $A\Rightarrow B=\operatorname*{int}(A^{c}\cup B)$ for open
subsets $A$ and $B$ which gives the classical definition if the topology is
discrete. Since we have an interior operator on the (non-topological) closure
space $U\times U$, this suggests that the implication partition $\sigma
\Rightarrow\pi$ might be defined by the dit-set definition:

\begin{center}
$\operatorname*{dit}\left(  \sigma\Rightarrow\pi\right)  =\operatorname*{int}%
\left(  \operatorname*{dit}\left(  \sigma\right)  ^{c}\cup\operatorname*{dit}%
\left(  \pi\right)  \right)  =\overline{\left(  \operatorname*{indit}\left(
\pi\right)  \cap\operatorname*{indit}\left(  \sigma\right)  ^{c}\right)  }%
^{c}$.
\end{center}

\noindent The equivalence relation that corresponds to a partition is its
indit set so the corresponding notion of the \textit{difference}
$\operatorname*{indit}\left(  \pi\right)  -\operatorname*{indit}\left(
\sigma\right)  $ between two equivalence relations would be the equivalence relation:

\begin{center}
$\operatorname*{indit}\left(  \pi\right)  -\operatorname*{indit}\left(
\sigma\right)  =\overline{\left(  \operatorname*{indit}\left(  \pi\right)
\cap\operatorname*{indit}\left(  \sigma\right)  ^{c}\right)  }%
=\operatorname*{indit}\left(  \sigma\Rightarrow\pi\right)
=\operatorname*{dit}\left(  \sigma\Rightarrow\pi\right)  ^{c}$.
\end{center}

\noindent The dit set $\operatorname*{dit}\left(  \sigma\Rightarrow\pi\right)
$ and its complement, the indit set $\operatorname*{indit}\left(
\sigma\Rightarrow\pi\right)  =\operatorname*{indit}\left(  \pi\right)
-\operatorname*{indit}\left(  \sigma\right)  $, define the \textit{same}
partition which is denoted $\sigma\Rightarrow\pi$ rather than say "$\pi
-\sigma$" since we have made the symmetry-breaking decision to define the
lattice of partitions to be isomorphic to the lattice of partition relations
rather than the opposite lattice of equivalence relations.

Since the dit-set definition of $\sigma\Rightarrow\pi$ involves the interior
operator on the closure space $U\times U$, it would be very convenient to have
a direct set-of-blocks definition of the implication partition $\sigma
\Rightarrow\pi$. From Boolean algebras and Heyting algebras, we can extract
one desideratum for the implication $\sigma\Rightarrow\pi$: if $\sigma\leq\pi$
in the partial order of the Boolean or Heyting algebra, then and only then
$\sigma\Rightarrow\pi=1$. Hence\ for any partitions $\sigma$ and $\pi$ on $U$,
if $\sigma$ is refined by $\pi$, i.e., $\sigma\preceq\pi$ in $\Pi(U)$, then
and only then we should have $\sigma\Rightarrow\pi=1$ (the discrete
partition).\footnote{The equality sign "$=$" is not a sign in the formal
language of partition logic so "$\sigma\Rightarrow\pi=1$" is not a formula in
that language. It simply says that the formulas "$\sigma\Rightarrow\pi$" \ and
"$1$" denote the same partitions in $\Pi(U)$.} The property is realized by the
simple set-of-blocks definition of the implication, temporarily denoted as
$\sigma\overset{\ast}{\Rightarrow}\pi$, that if a block $B\in\pi$ is contained
in a block $C\in\sigma$, then $B$ is "discretized," i.e., replaced by
singleton blocks $\left\{  u\right\}  $ for all $u\in B$, in the implication
$\sigma\overset{\ast}{\Rightarrow}\pi$ and otherwise the block $B$ remains the
same. The following proposition says that the dit-set definition is the same
as the set-of-blocks definition so that either may be used to define the
partition implication $\sigma\Rightarrow\pi$.

\begin{proposition}
$\sigma\Rightarrow\pi=\sigma\overset{\ast}{\Rightarrow}\pi$.
\end{proposition}

\noindent Proof: By the two definitions, $\operatorname*{dit}\left(
\pi\right)  \subseteq\operatorname*{dit}\left(  \sigma\Rightarrow\pi\right)  $
and $\operatorname*{dit}\left(  \pi\right)  \subseteq\operatorname*{dit}%
\left(  \sigma\overset{\ast}{\Rightarrow}\pi\right)  $ with the reverse
inclusions holding between the indit sets. We prove the proposition by showing
$\operatorname*{indit}\left(  \sigma\overset{\ast}{\Rightarrow}\pi\right)
\subseteq\operatorname*{indit}\left(  \sigma\Rightarrow\pi\right)  $ and
$\operatorname*{indit}\left(  \sigma\Rightarrow\pi\right)  \subseteq
\operatorname*{indit}\left(  \sigma\overset{\ast}{\Rightarrow}\pi\right)  $
where $\operatorname*{indit}\left(  \sigma\Rightarrow\pi\right)
=\overline{\left(  \operatorname*{indit}\left(  \pi\right)
-\operatorname*{indit}\left(  \sigma\right)  \right)  }=\overline{\left[
\operatorname*{dit}\left(  \sigma\right)  \cap\operatorname*{indit}\left(
\pi\right)  \right]  }$. Let $\left(  u,u^{\prime}\right)  \in
\operatorname*{indit}\left(  \sigma\overset{\ast}{\Rightarrow}\pi\right)  $
where $\operatorname*{indit}\left(  \sigma\overset{\ast}{\Rightarrow}%
\pi\right)  \subseteq\operatorname*{indit}\left(  \pi\right)  $ so that
$u,u^{\prime}\in B$ for some block $B\in\pi$. Moreover if $B$ were contained
in any block $C\in\sigma$, then $\left(  u,u^{\prime}\right)  \in
\operatorname*{dit}\left(  \sigma\overset{\ast}{\Rightarrow}\pi\right)
=\operatorname*{indit}\left(  \sigma\overset{\ast}{\Rightarrow}\pi\right)
^{c}$ contrary to assumption so $B$ is not contained in any $C\in\sigma$. If
$u$ and $u^{\prime}$ were in different blocks of $\sigma$ then $\left(
u,u^{\prime}\right)  \not \in \operatorname*{indit}\left(  \sigma\right)  $ so
that $\left(  u,u^{\prime}\right)  $ would not be subtracted off in the
formation of $\operatorname*{indit}\left(  \sigma\Rightarrow\pi\right)
=\overline{\left(  \operatorname*{indit}\left(  \pi\right)
-\operatorname*{indit}\left(  \sigma\right)  \right)  }$ and thus would be in
$\operatorname*{indit}\left(  \sigma\Rightarrow\pi\right)  $ which was to be
shown. Hence we may assume that $u$ and $u^{\prime}$ are in the same block
$C\in\sigma$. Thus $\left(  u,u^{\prime}\right)  $ was subtracted off in
$\operatorname*{indit}\left(  \pi\right)  -\operatorname*{indit}\left(
\sigma\right)  $ and we need to show that it is restored in the closure
$\overline{\left(  \operatorname*{indit}\left(  \pi\right)
-\operatorname*{indit}\left(  \sigma\right)  \right)  }$. Since $u,u^{\prime
}\in B\cap C$ but $B$ is not contained in any one block of $\sigma$, there is
another $\sigma$-block $C^{\prime}$ such that $B\cap C^{\prime}\not =%
\emptyset.$ Let $u^{\prime\prime}\in B\cap C^{\prime}$. Then $\left(
u,u^{\prime\prime}\right)  $ and $\left(  u^{\prime},u^{\prime\prime}\right)
$ are not in $\operatorname*{indit}\left(  \sigma\right)  $ since
$u,u^{\prime}\in C$ and $u^{\prime\prime}\in C^{\prime} $ but those two pairs
are in $\operatorname*{indit}\left(  \pi\right)  $ since $u,u^{\prime
},u^{\prime\prime}\in B$. Hence the pairs $\left(  u,u^{\prime\prime}\right)
,\left(  u^{\prime},u^{\prime\prime}\right)  \in\operatorname*{indit}\left(
\pi\right)  -\operatorname*{indit}\left(  \sigma\right)
=\operatorname*{indit}\left(  \pi\right)  \cap\operatorname*{dit}\left(
\sigma\right)  $ which implies that $\left(  u,u^{\prime}\right)  $ must be in
the closure $\operatorname*{indit}\left(  \sigma\Rightarrow\pi\right)
=\overline{\left(  \operatorname*{indit}\left(  \pi\right)
-\operatorname*{indit}\left(  \sigma\right)  \right)  }$. That establishes
$\operatorname*{indit}\left(  \sigma\overset{\ast}{\Rightarrow}\pi\right)
\subseteq\operatorname*{indit}\left(  \sigma\Rightarrow\pi\right)  $.

To prove the converse $\operatorname*{indit}\left(  \sigma\Rightarrow
\pi\right)  \subseteq\operatorname*{indit}\left(  \sigma\overset{\ast
}{\Rightarrow}\pi\right)  $, if $\left(  u,u^{\prime}\right)  \in
\operatorname*{indit}\left(  \sigma\Rightarrow\pi\right)  =\overline{\left[
\operatorname*{dit}\left(  \sigma\right)  \cap\operatorname*{indit}\left(
\pi\right)  \right]  }$, then there is a sequence $u=u_{1},u_{2}%
,...,u_{n}=u^{\prime}$ with every pair $\left(  u_{i},u_{i+1}\right)
\in\operatorname*{dit}\left(  \sigma\right)  \cap\operatorname*{indit}\left(
\pi\right)  $. Now $\left(  u_{i},u_{i+1}\right)  \in\operatorname*{indit}%
\left(  \pi\right)  $ implies there exists a block $B_{i}\in\pi$ with
$u_{i},u_{i+1}\in B_{i}$ for $i=1,...,n-1$. But $u_{i},u_{i+1}\in B_{i}$ and
$u_{i+1},u_{i+2}\in B_{i+1}$ implies $B_{i}=B_{i+1}$ so all the elements
$u_{i}$ belong to the same block $B\in\pi$ and in particular, $u,u^{\prime}\in
B$. Now if there was a $C\in\sigma$ with $B\subseteq C$, then, contrary to
assumption, we could not have any $\left(  u_{i},u_{i+1}\right)
\in\operatorname*{dit}\left(  \sigma\right)  $ since all the $u_{i}\in
B\subseteq C $. Hence there is no $C\in\sigma$ containing $B$ so $B$ would not
be discretized in $\sigma\overset{\ast}{\Rightarrow}\pi$ and thus $\left(
u,u^{\prime}\right)  \in\operatorname*{indit}\left(  \sigma\overset{\ast
}{\Rightarrow}\pi\right)  $. $\blacksquare$

Hence we may drop the temporary notation $\sigma\overset{\ast}{\Rightarrow}%
\pi$ and consider the partition implication $\sigma\Rightarrow\pi$ as
characterized by the set-of-blocks definition: form $\sigma\Rightarrow\pi$
from $\pi$ by discretizing any block $B\in\pi$ contained in a block
$C\in\sigma$.\footnote{For the analogy with subsets, the set difference
$X-Y=X\cap Y^{c}$ is obtained from $X$ by deleting any $u\in X$ that is
contained in $Y$, i.e., $\left\{  u\right\}  $ is locally replaced by the null
set, the zero element of the Boolean algebra of subsets of $U$. Similarly, in
the difference $\operatorname*{indit}\left(  \pi\right)
-\operatorname*{indit}(\sigma)$ of equivalence relations, any equivalence
class $B$ of $\pi$ contained in an equivalence class $C$ of $\sigma$ is
locally replaced by the zero in the lattice of equivalence relations, i.e., is
discretized.}

Another way to characterize the partition implication $\sigma\Rightarrow\pi$
is by using an adjunction.\footnote{See Mac Lane
\citeyear[p. 93]{mac:cwm}
for the notion of an adjunction between partial orders, i.e., a covariant
Galois connection.} In a Heyting algebra, the implication is characterized by
the adjunction which in our notation would be:

\begin{center}
$%
\begin{tabular}
[c]{c}%
$\tau\preceq\sigma\Rightarrow\pi$\\\hline
$\tau\wedge\sigma\preceq\pi$%
\end{tabular}
$.
\end{center}

\noindent For partitions, the top implies the bottom, but the bottom does not
imply the top. The simplest non-trivial partition algebra is that on the three
element set $U=\left\{  a,b,c\right\}  $ where we may take $\tau=\left\{
\left\{  a,b\right\}  ,\left\{  c\right\}  \right\}  $, $\sigma=\left\{
\left\{  a,c\right\}  ,\left\{  b\right\}  \right\}  $, and $\pi=\left\{
\left\{  a\right\}  ,\left\{  b,c\right\}  \right\}  $. Then $\tau\wedge
\sigma=0$ so the bottom $0\preceq\pi$ is true. But $\sigma\Rightarrow\pi=\pi$
(since no non-singleton block of $\pi$ is contained in a block of $\sigma$),
so the top is $\tau\preceq\pi$ which is false.

However, on the closure space $U\times U$, for any $S\subseteq U\times U$,
there is the usual adjunction $\mathcal{P}\left(  U\times U\right)
\rightleftarrows\mathcal{P}\left(  U\times U\right)  $ defining the set implication:

\begin{center}
$%
\begin{tabular}
[c]{c}%
$T\subseteq S\Rightarrow P$\\\hline
$T\cap S\subseteq P$%
\end{tabular}
$
\end{center}

\noindent(where $S\Rightarrow P$ is just $S^{c}\cup P$) for any subsets
$T,P\in\mathcal{P}\left(  U\times U\right)  $. Moreover, the dit-set
representation $\Pi\left(  U\right)  \rightarrow\mathcal{P}\left(  U\times
U\right)  $ where $\tau\longmapsto\operatorname*{dit}\left(  \tau\right)  $
has a right adjoint where $P\in\mathcal{P}\left(  U\times U\right)  $ is taken
to the partition $G\left(  P\right)  $ whose dit set is $\operatorname*{int}%
\left(  P\right)  $:

\begin{center}
$%
\begin{tabular}
[c]{c}%
$\tau\preceq G\left(  P\right)  $\\\hline
$\operatorname*{dit}\left(  \tau\right)  \subseteq P$%
\end{tabular}
$.
\end{center}

Composing the two right adjoints $\mathcal{P}\left(  U\times U\right)
\rightarrow\mathcal{P}\left(  U\times U\right)  \rightarrow\Pi(U)$ gives a
functor taking $P\in\mathcal{P}\left(  U\times U\right)  $ to $G_{S}\left(
P\right)  $ which is the partition whose dit set is $\operatorname*{int}%
\left(  S^{c}\cup P\right)  $. Its left adjoint is obtained by composing the
two left adjoints $\Pi(U)\rightarrow\mathcal{P}\left(  U\times U\right)
\rightarrow\mathcal{P}\left(  U\times U\right)  $ to obtain a functor taking a
partition $\tau$ to $F_{S}\left(  \tau\right)  =\operatorname*{dit}\left(
\tau\right)  \cap S$:

\begin{center}
$%
\begin{tabular}
[c]{c}%
$\tau\preceq G_{S}\left(  P\right)  $\\\hline
$F_{S}\left(  \tau\right)  \subseteq P$%
\end{tabular}
\ $.\footnote{Thanks to the referee for suggesting the simpler presentation of
this adjunction (as the composition of two adjunctions) as well as for other
helpful comments and proof-simplifying suggestions.}
\end{center}

Specializing $S=\operatorname*{dit}\left(  \sigma\right)  $ and
$P=\operatorname*{dit}\left(  \pi\right)  $ gives $G_{\operatorname*{dit}%
\left(  \sigma\right)  }\left(  \operatorname*{dit}\left(  \pi\right)
\right)  $ as the partition whose dit set is $\operatorname*{int}\left(
\operatorname*{dit}\left(  \sigma\right)  ^{c}\cup\operatorname*{dit}\left(
\pi\right)  \right)  $ which we know from above is the partition implication
$\sigma\Rightarrow\pi$, i.e., $G_{\operatorname*{dit}\left(  \sigma\right)
}\left(  \operatorname*{dit}\left(  \pi\right)  \right)  =\sigma\Rightarrow
\pi$. Using these restrictions, the adjunction gives the iff statement
characterizing the partition implication.

\begin{center}
$\operatorname*{dit}\left(  \tau\right)  \cap\operatorname*{dit}\left(
\sigma\right)  \subseteq\operatorname*{dit}\left(  \pi\right)  $ iff
$\tau\preceq\sigma\Rightarrow\pi$.

Characterization of $\sigma\Rightarrow\pi$
\end{center}

\noindent Thus $\sigma\Rightarrow\pi$ is the most refined partition $\tau$
such that $\operatorname*{dit}\left(  \tau\right)  \cap\operatorname*{dit}%
\left(  \sigma\right)  \subseteq\operatorname*{dit}\left(  \pi\right)  $. The
arbitrary intersection of equivalence relations (indit sets) is an equivalence
relation so the arbitrary union of dit sets is a dit set, i.e., the dit set of
the join of the partitions whose dit sets were in the union. Moreover,
distributivity in $\mathcal{P}\left(  U\times U\right)  $ implies that the
arbitrary union of dit sets $\operatorname*{dit}\left(  \tau\right)  $ such
that $\operatorname*{dit}\left(  \tau\right)  \cap\operatorname*{dit}\left(
\sigma\right)  \subseteq\operatorname*{dit}\left(  \pi\right)  $ will also
satisfy that same condition. Hence we may construct the most refined partition
$\tau$ such that $\operatorname*{dit}\left(  \tau\right)  \cap
\operatorname*{dit}\left(  \sigma\right)  \subseteq\operatorname*{dit}\left(
\pi\right)  $ by taking the join of those partitions:

\begin{center}
$\sigma\Rightarrow\pi=\bigvee\left\{  \tau:\operatorname*{dit}\left(
\tau\right)  \cap\operatorname*{dit}\left(  \sigma\right)  \subseteq
\operatorname*{dit}\left(  \pi\right)  \right\}  $.
\end{center}

The equivalent graph-theoretic definition of the partition implication can be
illustrated using the previous example.

\begin{example}
Let $\sigma=\left\{  \left\{  a,b,c\right\}  ,\left\{  d,e\right\}  \right\}
$ and $\pi=\left\{  \left\{  a,b\right\}  ,\left\{  c,d,e\right\}  \right\}  $
as before. In the graph below, all the arcs in the complete graph $K_{5}$ on
five nodes are again labelled according to the status of the two endpoints in
the two partitions. The Boolean conditions for $F\left(  \sigma\Rightarrow
\pi\right)  $ are "$T\sigma$ and $F\pi$" . The arcs where those conditions
hold are the solid lines. In the graph with only the solid arcs, there are
three connected components giving the blocks of the implication:
$\sigma\Rightarrow\pi=\left\{  \left\{  a\right\}  ,\{b\},\{c,d,e\}\right\}
$. Note that only the $\pi$-block $\left\{  a,b\right\}  $ is contained in a
$\sigma$-block so $\sigma\Rightarrow\pi$ is like $\pi$ except that $\left\{
a,b\right\}  $ is discretized.
\end{example}

%

\begin{center}
\includegraphics[
natheight=182.005005bp,
natwidth=310.024902bp,
height=185.6875pt,
width=314.875pt
]%
{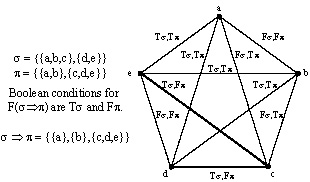}%
\\
Figure 2: Graph for implication $\sigma\Rightarrow\pi$%
\end{center}

\subsection{Partition negation operation}

In intuitionistic logic, the negation $\lnot\sigma$ would be defined as the
implication $\sigma\Rightarrow0$ with the consequent taken as the zero element
$0$, i.e., $\lnot\sigma=\sigma\Rightarrow0$. In the topological interpretation
using open subsets, $\sigma$ would be an open subset and $\lnot\sigma$ would
be the interior of its complement. Adapted to partitions, these give the same
dit-set definition of the partition negation (since $\operatorname*{dit}%
\left(  0\right)  =\emptyset$):

\begin{center}
$\operatorname*{dit}\left(  \lnot\sigma\right)  =\operatorname*{int}\left(
\operatorname*{dit}\left(  \sigma\right)  ^{c}\right)  =\operatorname*{dit}%
\left(  \sigma\Rightarrow0\right)  $.
\end{center}

\noindent It is a perhaps surprising fact that this dit set is always empty
(so that $\lnot\sigma=0$) except in the singular case where $\sigma=0$ in
which case we, of course, have $\lnot0=\left(  0\Rightarrow0\right)
=1$.\footnote{In graph theory, this is the result that given any disconnected
(simple) graph $G$, its complement $G^{c}$ (set of all links not in $G$) is
connected \cite[p. 30]{wil:gt}.} The key fact is that any two partitions
(aside from the blob) must have some dits in common.

\begin{theorem}
[Common-dits theorem]Any two non-empty dit sets have some dits in common.
\end{theorem}

\noindent Proof: Let $\pi$ and $\sigma$ be any two partitions on $U$ with
non-empty dit sets, i.e., $\pi\not =0\not =\sigma$. We need to show that
$\operatorname*{dit}\left(  \pi\right)  \cap\operatorname*{dit}\left(
\sigma\right)  \not =\emptyset$. Since $\pi$ is not the blob $0$, consider two
elements $u$ and $u^{\prime}$ distinguished by $\pi$ but identified by
$\sigma$ [otherwise $\left(  u,u^{\prime}\right)  \in\operatorname*{dit}%
\left(  \pi\right)  \cap\operatorname*{dit}\left(  \sigma\right)  $]. Since
$\sigma$ is also not the blob, there must be a third element $u^{\prime\prime
}$ not in the same block of $\sigma$ as $u$ and $u^{\prime}$.%

\begin{center}
\includegraphics[
natheight=83.997498bp,
natwidth=112.017403bp,
height=86.4259bp,
width=114.6326bp
]%
{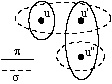}%
\end{center}

\begin{center}
Figure 3: $\left(  u,u^{\prime\prime}\right)  $ as a common dit
\end{center}

\noindent But since $u$ and $u^{\prime}$ are in different blocks of $\pi$, the
third element $u^{\prime\prime}$ must be distinguished from one or the other
or both in $\pi$. Hence $\left(  u,u^{\prime\prime}\right)  $ or $\left(
u^{\prime},u^{\prime\prime}\right)  $ must be distinguished by both partitions
and thus must be in $\operatorname*{dit}\left(  \pi\right)  \cap
\operatorname*{dit}\left(  \sigma\right)  $. $\blacksquare$

It should be noted that the interior of the intersection $\operatorname*{dit}%
\left(  \pi\right)  \cap\operatorname*{dit}\left(  \sigma\right)  $ could be
empty, i.e., $\sigma\wedge\pi=0$, even when the intersection is non-empty. It
might also be useful to consider the contrapositive form of the common-dits
theorem which is about equivalence relations. If the union of two equivalence
relations is the universal equivalence relation, i.e., $\operatorname*{indit}%
\left(  \pi\right)  \cup\operatorname*{indit}\left(  \sigma\right)  =U\times U
$, then one of the equivalence relations is the universal one, i.e.,
$\operatorname*{indit}\left(  \pi\right)  =U\times U$ or
$\operatorname*{indit}\left(  \sigma\right)  =U\times U$.

For any non-blob partition $\sigma$, $\operatorname*{dit}\left(  \lnot
\sigma\right)  =\operatorname*{int}\left(  \operatorname*{dit}\left(
\sigma\right)  ^{c}\right)  $ is a dit set disjoint from the non-empty
$\operatorname*{dit}\left(  \sigma\right)  $ so by the common-dits theorem, it
has to be empty and thus $\lnot\sigma=0 $. Negation becomes more useful if we
generalize by replacing the blob in the definition $\lnot\sigma=\sigma
\Rightarrow0$ by an arbitrary but fixed partition $\pi$. This leads to the
notion of the $\pi$\textit{-negation} of a partition $\sigma$ which is just
the implication $\sigma\Rightarrow\pi$ with the fixed partition $\pi$ as the
consequent. We added a $\pi$ to the negation symbol to represent this negation
relative to $\pi$:

\begin{center}
$\pi$-negation: $\overset{\pi}{\lnot}\sigma=\sigma\Rightarrow\pi$.
\end{center}

\noindent The unadorned negation $\lnot\sigma$ is the $0$-negation, i.e.,
$\lnot\sigma=\sigma\Rightarrow0$. Using this suggestive notation, the
partition tautology that internalizes modus ponens, $\left(  \sigma
\wedge\left(  \sigma\Rightarrow\pi\right)  \right)  \Rightarrow\pi$, is the
law of non-contradiction, $\overset{\pi}{\lnot}\left(  \sigma\wedge
\overset{\pi}{\lnot}\sigma\right)  $, for $\pi$-negation. While it is useful
to establish the notion of partition negation, it need not be taken as a
primitive operation.

\subsection{Partition stroke, not-and, or nand operation}

In addition to the lattice operations of the join and meet, and the
implication operation, we introduce the \textit{Sheffer stroke},
\textit{not-and}, or \textit{nand} operation $\sigma\mid\tau$, with the
dit-set definition:

\begin{center}
$\operatorname*{dit}\left(  \sigma\mid\tau\right)  =\operatorname*{int}\left[
\operatorname*{indit}\left(  \sigma\right)  \cup\operatorname*{indit}\left(
\tau\right)  \right]  $.
\end{center}

For a set-of-blocks definition consider a graph whose nodes are the elements
$u\in U$. Given $\sigma=\left\{  C\right\}  $ and $\tau=\left\{  D\right\}  $,
each element $u$ is in a unique block $C\cap D$ of the join $\sigma\vee\tau$.
Given elements $u\in C\cap D$ and $u^{\prime}\in C^{\prime}\cap D^{\prime}$,
$u$ is connected by an arc or link in the graph, i.e., $u\sim u^{\prime}$, if
$C\not =C^{\prime}$ and $D\not =D^{\prime}$, i.e., if $\left(  C\cap D\right)
\times\left(  C^{\prime}\cap D^{\prime}\right)  \subseteq\operatorname*{dit}%
\left(  \sigma\right)  \cap\operatorname*{dit}\left(  \tau\right)  =\left[
\operatorname*{indit}\left(  \sigma\right)  \cup\operatorname*{indit}\left(
\tau\right)  \right]  ^{c} $. Then the nodes in each connected component of
the graph are the blocks of $\sigma\mid\tau$. Two nodes $u,u^{\prime}$ are
connected in this graph if and only if the ordered pair $\left(  u,u^{\prime
}\right)  $ is in the closure $\overline{\left(  \operatorname*{dit}\left(
\sigma\right)  \cap\operatorname*{dit}\left(  \tau\right)  \right)
}=\overline{\left[  \operatorname*{indit}\left(  \sigma\right)  \cup
\operatorname*{indit}\left(  \tau\right)  \right]  ^{c}}$, and thus they are a
distinction if and only if they are in the complement of the closure which is
the interior: $\operatorname*{int}\left[  \operatorname*{indit}\left(
\sigma\right)  \cup\operatorname*{indit}\left(  \tau\right)  \right]  $. Hence
this graph-theoretic definition of the nand operation is the same as the
dit-set definition.

To turn it into a set-of-blocks definition, note that when $u\sim u^{\prime} $
because $C\not =C^{\prime}$ and $D\not =D^{\prime}$ then all the elements of
$C\cap D$ and $C^{\prime}\cap D^{\prime}$ are in the same block of the nand
$\sigma\mid\tau$. But if for any non-empty $C\cap D$, there is no other block
$C^{\prime}\cap D^{\prime}$ of the join with $C\not =C^{\prime} $ and
$D\not =D^{\prime}$, then the elements of $C\cap D$ would not even be
connected with each other so they would be singletons in the nand. Hence for
the set-of-blocks definition of the nand $\sigma\mid\tau$, the blocks of the
nand partition are formed by taking the unions of any join blocks $C\cap D$
and $C^{\prime}\cap D^{\prime}$ which differ in both "components" but by
taking as singletons the elements of any $C\cap D$ which does not differ from
any other join block in both components.

\begin{example}
Let $\sigma=\left\{  \left\{  a,b,c\right\}  ,\left\{  d,e\right\}  \right\}
$ and $\pi=\left\{  \left\{  a,b\right\}  ,\left\{  c,d,e\right\}  \right\}  $
as before. In the graph below, all the arcs in the complete graph $K_{5}$ on
five nodes are again labelled according to the status of the two endpoints in
the two partitions. The Boolean conditions for $F\left(  \sigma\mid\pi\right)
$ are "$T\sigma$ and $T\pi$" . The arcs where those conditions hold are the
solid lines. In the graph with only the solid arcs, there are two connected
components giving the blocks of the nand: $\sigma\mid\pi=\left\{  \left\{
a,b,d,e\right\}  ,\{c\}\right\}  $.
\end{example}

%

\begin{center}
\includegraphics[
natheight=182.005005bp,
natwidth=314.010010bp,
height=185.6875pt,
width=318.875pt
]%
{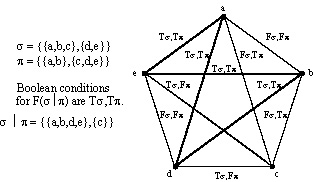}%
\\
Figure 4: Graph for nand $\sigma\mid\tau$%
\end{center}

\begin{example}
If $\sigma=\left\{  C,C^{\prime}\right\}  $ where $C=\left\{  u\right\}  $ and
$C^{\prime}=U-\left\{  u\right\}  $ and $\tau=\left\{  D,D^{\prime}\right\}  $
where $D=U-\left\{  u^{\prime}\right\}  $ and $D^{\prime}=\left\{  u^{\prime
}\right\}  $, then $\sigma\vee\tau=\left\{  \left\{  u\right\}  ,\left\{
u^{\prime}\right\}  ,U-\left\{  u,u^{\prime}\right\}  \right\}  $. Hence $u\in
C\cap D=\left\{  u\right\}  \cap\left(  U-\left\{  u^{\prime}\right\}
\right)  $ and $u^{\prime}\in C^{\prime}\cap D^{\prime}=\left(  U-\left\{
u\right\}  \right)  \cap\left\{  u^{\prime}\right\}  $ so $u\sim u^{\prime}$
in the graph for $\sigma\mid\tau$. But the elements $u^{\prime\prime}\in
C^{\prime}\cap D=U-\left\{  u,u^{\prime}\right\}  $ are not connected to any
other elements since $C^{\prime}\cup D=\left(  U-\left\{  u\right\}  \right)
\cup\left(  U-\left\{  u^{\prime}\right\}  \right)  =U$ so they are all
singletons in the nand. Hence $\sigma\mid\tau=\left\{  \left\{  u,u^{\prime
}\right\}  ,\left\{  u^{\prime\prime}\right\}  ,...\right\}  $.
\end{example}

This example can be stated in more general terms. A \textit{modular} partition
is a partition with at most one non-singleton block. A non-zero partition
$\varphi$ is an \textit{atom} in the lattice of partitions $\Pi(U)$ if
$0\preceq\pi\preceq\varphi$ implies $\pi=0$ or $\pi=\varphi$. A non-unitary
partition $\varphi$ is a \textit{coatom} if $\varphi\preceq\pi\preceq1$
implies $\pi=\varphi$ or $\pi=1$. All coatoms are modular where the
non-singleton block is some pair $\left\{  u,u^{\prime}\right\}  $. The
example then shows that the nand of any two distinct modular atoms is a coatom.

For subsets $S,T\subseteq U$, the nand subset $S\mid T=S^{c}\cup T^{c}=\left(
S\cap T\right)  ^{c}$ has as elements those elements $u\in U$ which are not
elements of both $S$ and $T$. Using the relationship between elements of a
subset and distinctions of a partition, the nand partition $\sigma\mid\tau$
has as distinctions those pairs $\left(  u,u^{\prime}\right)  \in U\times
U-\Delta$ which are, directly or indirectly, not distinctions of both $\sigma$
and $\tau$. In the example above, $\left(  u,u^{\prime}\right)  $ is a
distinction of both $\sigma$ and $\tau$ so it is not a distinction of
$\sigma\mid\tau$. For any third element $u^{\prime\prime}\in U$, then
$u^{\prime\prime}$ paired with any other element of $U$ is not a dit of both
$\sigma$ and $\tau$ so the pair is a distinction of $\sigma\mid\tau$, i.e.,
$\left\{  u^{\prime\prime}\right\}  $ is a singleton in the nand partition.

A number of the relations which we are accustomed to in subset logic also hold
in partition logic. For instance, negation can be defined using the nand:
$\sigma\mid\sigma=\lnot\sigma$. In fact, if $\sigma\preceq\tau$, then
$\sigma\mid\tau=\lnot\sigma$. For example, since $\sigma$ is always refined by
$\tau\Rightarrow\sigma$ for any $\tau$, $\sigma\mid\left(  \tau\Rightarrow
\sigma\right)  =\lnot\sigma$. The formula $\sigma\mid\sigma=\lnot\sigma$ is
also a special case of the formula $\left(  \sigma\mid\tau\right)
\wedge\left(  \sigma\Rightarrow\tau\right)  =\lnot\sigma$ derived in the next section.

In subset logic, the "and" and the nand subsets would be complements of one
another but the relationship is more subtle in partition logic. We say that
two partitions $\varphi$ and $\varphi^{\prime}$ which refine a partition $\pi
$, i.e., $\pi\preceq\varphi,\varphi^{\prime}$, are $\pi$-\textit{orthogonal}
if $\overset{\pi}{\lnot}\varphi\vee\overset{\pi}{\lnot}\varphi^{\prime}=1$.
Since all partitions refine $0$, two partitions $\varphi$ and $\varphi
^{\prime}$ are $0$\textit{-orthogonal} or, simply, \textit{orthogonal} if
$\lnot\varphi\vee\lnot\varphi^{\prime}=1$. This may look odd as a criterion
for orthogonality but it is classically equivalent to the more familiar
$\varphi\wedge\varphi^{\prime}=0$.

\begin{lemma}
$\varphi$ and $\varphi^{\prime}$ are orthogonal, i.e., $\lnot\varphi\vee
\lnot\varphi^{\prime}=1$, iff $\varphi\mid\varphi^{\prime}=1$.
\end{lemma}

\noindent Proof: If $\lnot\varphi\vee\lnot\varphi^{\prime}=1$, then
$\operatorname*{int}\left(  \operatorname*{indit}\left(  \varphi\right)
\right)  \cup\operatorname*{int}\left(  \operatorname*{indit}\left(
\varphi^{\prime}\right)  \right)  =\operatorname*{dit}\left(  1\right)
=U^{2}-\Delta$. By the monotonicity of the interior operator,
$\operatorname*{int}\left(  \operatorname*{indit}\left(  \varphi\right)
\right)  \cup\operatorname*{int}\left(  \operatorname*{indit}\left(
\varphi^{\prime}\right)  \right)  \subseteq\operatorname*{int}\left(
\operatorname*{indit}\left(  \varphi\right)  \cup\operatorname*{indit}\left(
\varphi^{\prime}\right)  \right)  =\operatorname*{dit}\left(  \varphi
\mid\varphi^{\prime}\right)  $ so $\varphi\mid\varphi^{\prime}=1$. Conversely
if $\varphi\mid\varphi^{\prime}=1$, then $\operatorname*{int}\left(
\operatorname*{indit}\left(  \varphi\right)  \cup\operatorname*{indit}\left(
\varphi^{\prime}\right)  \right)  =\operatorname*{dit}\left(  1\right)
=U^{2}-\Delta$. Since $\Delta\subseteq\operatorname*{indit}\left(
\varphi\right)  ,\operatorname*{indit}\left(  \varphi^{\prime}\right)  $ (so
that only $\Delta$ is removed by the interior operator),
$\operatorname*{indit}\left(  \varphi\right)  \cup\operatorname*{indit}\left(
\varphi^{\prime}\right)  =U^{2}$. It was previously noted that if the union of
two equivalence relations is the universal equivalence relation $U^{2}$, then
one of the equivalence relations must be the universal one. Hence either
$\varphi=0$ or $\varphi^{\prime}=0$ and since $\lnot0=1$, we have either way,
$\lnot\varphi\vee\lnot\varphi^{\prime}=1$. $\blacksquare$

Just as the unary negation operation $\lnot\varphi$ is usefully generalized by
the binary operation $\overset{\pi}{\lnot}\varphi=\varphi\Rightarrow\pi$, so
the binary nand operation $\sigma\mid\tau$ is usefully generalized by the
ternary operation of $\pi$\textit{-nand }defined by:

\begin{center}
$\operatorname*{dit}\left(  \sigma\mid_{\pi}\tau\right)  =\operatorname*{int}%
\left(  \operatorname*{indit}\left(  \sigma\right)  \cup\operatorname*{indit}%
\left(  \tau\right)  \cup\operatorname*{dit}\left(  \pi\right)  \right)  $.
\end{center}

\noindent Then a similar argument shows that for $\pi\preceq\varphi
,\varphi^{\prime}$:

\begin{center}
$\varphi$ and $\varphi^{\prime}$ are $\pi$-orthogonal iff $\varphi\mid_{\pi
}\varphi^{\prime}=1$.
\end{center}

Thus two partitions are orthogonal when if one of the partitions is non-zero,
then the other partition must be zero (i.e., at least one is zero). If
$\varphi$ and $\varphi^{\prime}$ are orthogonal, i.e., $\varphi\mid
\varphi^{\prime}=1$, then $\varphi\wedge\varphi^{\prime}=0$ follows but not
vice-versa (see next example).

Every partition $\sigma$ and its $0$-negation $\lnot\sigma$ are orthogonal
since $\lnot\sigma\vee\lnot\lnot\sigma=1$. In the example above, the meet of
$\sigma=\left\{  \left\{  u\right\}  ,U-\left\{  u\right\}  \right\}  $ and
$\tau=\left\{  \left\{  u^{\prime}\right\}  ,U-\left\{  u^{\prime}\right\}
\right\}  $ is $\sigma\wedge\tau=0$ and $\lnot0=1$ but $\sigma\mid\tau\not =1$
so the negation $\lnot\left(  \sigma\wedge\tau\right)  $ operation is not
necessarily the same as the nand $\sigma\mid\tau$ operation. However, the
"and" or meet $\sigma\wedge\tau$ and the "not-and" or nand $\sigma\mid\tau$
are orthogonal; if one is non-zero, the other must be zero. Thus no pair
$\left(  u,u^{\prime}\right)  $ can be a dit of both and hence $\left(
\sigma\mid\tau\right)  \mid\left(  \sigma\wedge\tau\right)  =1$ is a partition
tautology. The same example above shows that the nand $\sigma\mid\tau$ is also
not the same as $\lnot\sigma\vee\lnot\tau$ (which equals $0$ in the example).
Although the three formulas are equal in subset logic, in partition logic we
only have the following refinement relations holding in general:

\begin{center}
$\lnot\sigma\vee\lnot\tau\preceq\sigma\mid\tau\preceq\lnot\left(  \sigma
\wedge\tau\right)  $.
\end{center}

\noindent Thus only one direction $\lnot\sigma\vee\lnot\tau\preceq\lnot\left(
\sigma\wedge\tau\right)  $ holds in general so the "strong" DeMorgan law
$\lnot\sigma\vee\lnot\tau=\lnot\left(  \sigma\wedge\tau\right)  $ does not
hold in partition logic. However, the other "weak" DeMorgan law holds in
partition logic even for $\pi$-negation, i.e., $\overset{\pi}{\lnot}\left(
\sigma\vee\tau\right)  =\overset{\pi}{\lnot}\sigma\wedge\overset{\pi}{\lnot
}\tau$.

\begin{example}
The universe set $U=\left\{  Tom,John,Jim\right\}  $ consists of three people
and there are two partitions: $\alpha$ which distinguishes people according to
the first letter of their name so that $\alpha=\left\{  \left\{  Tom\right\}
,\left\{  John,Jim\right\}  \right\}  $, and $\omega$ which distinguishes
people according to the last letter of their name so that $\omega=\left\{
\left\{  Tom,Jim\right\}  ,\left\{  John\right\}  \right\}  $. Then the meet
$\alpha\wedge\omega$ would identify people who are directly and indirectly
identified by the two partitions. Tom and John are not directly identified but
are indirectly identified: $Tom\overset{\omega}{\sim}Jim\overset{\alpha}{\sim
}John$ so that $\sigma\wedge\omega=0$. But since the meet is $0$,
the\ $0$-orthogonal nand of the two partitions could be non-zero, and in fact
$\alpha\mid\omega=\left\{  \left\{  Tom,John\right\}  ,\left\{  Jim\right\}
\right\}  $. Thus the fact that Tom and John are directly distinguished by
both the first and last letters of their names results in them not being
distinguished by the not-and partition.
\end{example}

In any dit-set definition of a partition $\varphi$ as $\operatorname*{dit}%
\left(  \varphi\right)  =\operatorname*{int}\left(  P\right)  $ for some
$P\subseteq U\times U $, two elements $u$ and $u^{\prime}$ will be in the same
block of $\varphi$ if and only if they are in the closure $\overline{\left(
P^{c}\right)  }$, i.e., if there is a finite sequence of links $\left(
u_{i},u_{i+1}\right)  \in P^{c}$ connecting $u$ and $u^{\prime}$. The question
arises of there being an upper bound on the number of links required to put
two elements in the same block. In the simple case of the join $\sigma\vee
\tau$ where $\operatorname*{dit}\left(  \sigma\vee\tau\right)
=\operatorname*{dit}\left(  \sigma\right)  \cup\operatorname*{dit}\left(
\tau\right)  $, no interior operator is needed since the union of open subsets
of the closure space $U\times U$ is open. Thus the complement $\left(
\operatorname*{dit}\left(  \sigma\right)  \cup\operatorname*{dit}\left(
\tau\right)  \right)  ^{c}=\operatorname*{indit}\left(  \sigma\right)
\cap\operatorname*{indit}\left(  \tau\right)  $ is already closed (i.e., the
intersection of two equivalence relations is an equivalence relation) so one
link $\left(  u,u^{\prime}\right)  \in\operatorname*{indit}\left(
\sigma\right)  \cap\operatorname*{indit}\left(  \tau\right)  $ suffices to put
$u$ and $u^{\prime}$ into the same block of the join $\sigma\vee\tau$. Thus
for the join, one link suffices.

For the implication $\sigma\Rightarrow\tau$, $\left(  u,u^{\prime}\right)
\in\operatorname*{indit}\left(  \sigma\Rightarrow\tau\right)  $ if and only if
$\left(  u,u^{\prime}\right)  \in\operatorname*{indit}\left(  \tau\right)  $,
say, $u,u^{\prime}\in D\in\tau$, and there is no $C\in\sigma$ such that
$D\subseteq C$ so the block $D$ remains whole in the implication
$\sigma\Rightarrow\tau$. But that means there is another block $C^{\prime}%
\in\sigma$ such that $D\cap C^{\prime}\not =\emptyset$, i.e., there is an
$a\in D\cap C^{\prime}$ such that $\left(  u,a\right)  $ and $\left(
a,u^{\prime}\right)  $ are both dits of $\sigma$ but indits of $\tau.$ Thus
there is at most a two link chain connecting $u$ and $u^{\prime}$ where each
link is in $\operatorname*{dit}\left(  \sigma\right)  \cap
\operatorname*{indit}\left(  \tau\right)  =\left(  \operatorname*{dit}\left(
\sigma\right)  ^{c}\cup\operatorname*{dit}\left(  \tau\right)  \right)  ^{c}$.
Thus for the implication, two links suffice.

For the meet of two partitions, it is well-known that there is no upper bound
on the finite number of links needed to connect two elements which are in the
same block. For instance on the natural numbers, take $\sigma=\left\{
\left\{  0,1\right\}  ,\left\{  2,3\right\}  ,...\right\}  $ and
$\tau=\left\{  \left\{  0\right\}  ,\left\{  1,2\right\}  ,\left\{
3,4\right\}  ,...\right\}  $ so that $\sigma\wedge\tau$ is the blob and thus
any two elements are connected. But clearly there is no upper bound on the
number of links needed to connect any two elements.

For the nand operation, it is perhaps interesting that four links suffice. To
show this, we first exhibit an example where four links are required, i.e., no
shorter set of links would suffice. Then we show that in general, longer
chains can always be shortened to four or fewer links..

For an example where four links are required, consider the four-link chain
$u,a,b,c,u^{\prime}$ connecting $u$ and $u^{\prime}$ in the nand $\sigma
\mid\tau$ where $\sigma=\left\{  \left\{  u,u^{\prime},b\right\}  ,\left\{
a,c\right\}  \right\}  $ and $\tau=\left\{  \left\{  u,c\right\}  ,\left\{
u^{\prime},a\right\}  ,\left\{  b\right\}  \right\}  $. Each link $\left(
u,a\right)  $, $\left(  a,b\right)  $, $\left(  b,c\right)  $, and $\left(
c,u^{\prime}\right)  $ in the four-link chain is in the set $\left(
\operatorname*{indit}\left(  \sigma\right)  \cup\operatorname*{indit}\left(
\tau\right)  \right)  ^{c}=\operatorname*{dit}\left(  \sigma\right)
\cap\operatorname*{dit}\left(  \tau\right)  $ so $\left(  u,u^{\prime}\right)
$ is in its closure, i.e., $u$ and $u^{\prime} $ are in the same block of
$\sigma\mid\tau=0$. And there are no short-cuts. By placing the five points on
the vertices of a pentagon, then it is easy to see that none of the
short-cutting chords are in $\operatorname*{dit}\left(  \sigma\right)
\cap\operatorname*{dit}\left(  \tau\right)  $.

\begin{lemma}
Four links suffice to put any two elements in the same block of any nand
$\sigma\mid\tau$.
\end{lemma}

\noindent Proof: The proof can be formulated abstractly using sequences of
ordered pairs which can be pictured as points on the plane. Suppose we have a
chain of ordered pairs $\left(  x_{1},y_{1}\right)  $, $\left(  x_{2}%
,y_{2}\right)  $,..., $\left(  x_{n},y_{n}\right)  $ where each pair differs
from the previous one on both coordinates. However if any pair differs on both
coordinates with a previous pair, then all intermediate pairs could be cut out
thus shortening the chain. We want to construct a subchain with four or less
links. Since we cannot just directly connect the end points they must agree on
one coordinate such as the $x$ coordinate. Then $\left(  x_{2},y_{2}\right)  $
must also agree on one coordinate with $\left(  x_{n},y_{n}\right)  $ or we
would just connect them and be finished with a two-link chain. But they cannot
agree on the $x$ coordinate since it has to differ on both coordinates from
the first point $\left(  x_{1},y_{1}\right)  $. Hence it has to agree on the
$y$ coordinate with $\left(  x_{n},y_{n}\right)  $.%

\begin{center}
\includegraphics[
natheight=100.000000bp,
natwidth=114.010002bp,
height=102.5529bp,
width=116.6252bp
]%
{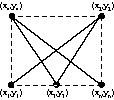}%
\end{center}

\begin{center}
Figure 5: Four links suffice
\end{center}

The third point $\left(  x_{3},y_{3}\right)  $ must agree with the first and
last which means on the $x$ coordinate as pictured above. Then the fourth
point $\left(  x_{4},y_{4}\right)  $ must differ from $\left(  x_{3}%
,y_{3}\right)  $ on both coordinates but must agree with the first and second
points on some coordinates. Thus it must agree with the first point on the $y$
coordinate and with the second point on the $x$ coordinate. But then it will
differ from the last point $\left(  x_{n},y_{n}\right)  $ on both coordinates
so it can be directly connected giving a four link subchain where each
successive pair differs on both coordinates.

To map this abstract proof into the case at hand, recall that
$\operatorname*{indit}\left(  \sigma\mid\tau\right)  =\overline{\left[
\operatorname*{dit}\left(  \sigma\right)  \cap\operatorname*{dit}\left(
\tau\right)  \right]  }$ so that $\left(  u,u^{\prime}\right)  $ is an indit
of $\sigma\mid\tau$ if there is a finite sequence $u=u_{1},u_{2}%
,...,u_{n}=u^{\prime}$ with each pair $\left(  u_{i},u_{i+1}\right)
\in\operatorname*{dit}\left(  \sigma\right)  \cap\operatorname*{dit}\left(
\tau\right)  $. Different horizontal coordinates correspond to different
$\sigma$ blocks and different vertical coordinates correspond to different
$\tau$ blocks where only a finite number of points are needed to model the
finite sequence. The first link $\left(  u_{1},u_{2}\right)  $ then maps to
the first line segment from $\left(  x_{1},y_{1}\right)  $ (the pair of
coordinates representing the blocks $C_{1}\in\sigma$ and $D_{1}\in\tau$
containing $u_{1}=u$) to $\left(  x_{2},y_{2}\right)  $ (the pair representing
the different blocks $C_{2}\in\sigma$ and $D_{2}\in\tau$ containing $u_{2}$).
The second link $\left(  u_{2},u_{3}\right)  $ maps to the second line segment
from $\left(  x_{2},y_{2}\right)  $ to $\left(  x_{3},y_{3}\right)  $, and so
forth. $\blacksquare$

\subsection{Sixteen binary operations on partitions}

What other partition operations might be defined? For binary operations
$\sigma\ast\tau$ on Boolean $0,1$ variables $\sigma$ and $\tau$, there are
four combinations of values for $\sigma$ and $\tau$, and thus there are
$2^{4}=16$ possible binary Boolean operations: $2\times2\rightarrow2$.
Thinking in terms of subsets $S,T\subseteq U$ instead of Boolean propositional
variables, there are the four basic disjoint regions in the general position
Venn diagram for $S$ and $T$, namely $S\cap T$, $S\cap T^{c} $, $S^{c}\cap T$,
and $S^{c}\cap T^{c}$. Then there are again $2^{4}=16 $ subsets of $U$ defined
by including or not including each of these four basic regions. That defines
the $16$ binary \textit{logical} operations on subsets of $U$.

Now take $S=\operatorname*{dit}\left(  \sigma\right)  $ and
$T=\operatorname*{dit}\left(  \tau\right)  $ as subsets of $U\times U$ and
define the $16$ subsets of $U\times U$ in the same way. Some of these such as
$S\cup T=\operatorname*{dit}\left(  \sigma\right)  \cup\operatorname*{dit}%
\left(  \tau\right)  =\operatorname*{dit}\left(  \sigma\vee\tau\right)  $ will
be open and thus will be the dit sets of partitions on $U$. For those which
are not already open, we must apply the interior operator to get the dit set
of a partition on $U$. This gives $16$ binary operations on partitions that
would naturally be called \textit{logical} since they are immediately paired
with the corresponding $16$ binary logical operations on subsets. We will use
the same notation for the partition operations. For instance, for subsets
$S,T\subseteq U$, the conditional or implication subset is $S^{c}\cup
T=S\Rightarrow T$. When $S=\operatorname*{dit}\left(  \sigma\right)  $ and
$T=\operatorname*{dit}\left(  \tau\right)  $ as subsets of $U\times U$, the
subset $S^{c}\cup T$ is not necessarily open so we must apply the interior
operator to get the dit set defining the corresponding implication operation
on partitions, i.e., $\operatorname*{int}\left[  \operatorname*{dit}\left(
\sigma\right)  ^{c}\cup\operatorname*{dit}\left(  \tau\right)  \right]
=\operatorname*{dit}\left(  \sigma\Rightarrow\tau\right)  $.

Alternatively, the sixteen logical partition operations could be defined using
the graph-theoretic approach rather than the dit-set approach. Given two
partitions $\sigma$ and $\tau$ on $U$, we label each arc in the complete
undirected simple graph on the node set $U$ according to whether or not the
end points of the arc were distinctions (with labels $T\sigma$ or $T\tau$) or
indistinctions ($F\sigma$ or $F\tau$) of the partitions. Thus each arc is
associated with a row in the Boolean truth table for the two atomic variables
$\sigma$ and $\tau$. For any of the sixteen Boolean logical operations
$\sigma\ast\tau$, we keep the arcs where $\sigma\ast\tau$ was assigned a $0$
or $F\left(  \sigma\ast\tau\right)  $ in the truth table for that Boolean
operation and discard the other arcs. Then the blocks in the partition
$\sigma\ast\tau$ are the nodes in the connected components of that graph. This
gives the same sixteen logical binary operations on partitions as the dit-set approach.

In both subset and partition logic, there are only two nullary operations
(constants), $0$ and $1$. With unary operations, the situation is still
straightforward. There are only four subset logical unary operations: identity
and negation (or complementation) in addition to the two nullary operations
(seen as constant unary operations). These immediately yield the partition
operations of identity $\sigma$ and negation $\lnot\sigma$ in addition to the
two partition constant operations $0$ and $1$. If these partition operations
are compounded using the logical operations such as negation, implication,
join, meet, and nand, then two other distinct unary operations are generated:
the double negation $\lnot\lnot\sigma$ and the excluded middle operation
$\sigma\vee\lnot\sigma$ (which is also equal to $\lnot\lnot\sigma
\Rightarrow\sigma$, the direction of the usual law of double negation that is
not a partition tautology)--to make six logical unary partition operations.

The situation for binary partition operations is considerably more
complicated. If the sixteen binary operations on subsets are compounded, then
the result is always one of the sixteen binary subset operations, e.g.,
$S\cap\left(  S\Rightarrow T\right)  =S\cap T$. But the presence of the
interior operator significantly changes the partition case. Compounding gives
many new binary operations on partitions, e.g., $\lnot\left(  \sigma\wedge
\tau\right)  $ and $\lnot\sigma\vee\lnot\tau$ (noted in the analysis of
$\sigma\mid\tau$), and they could also be called "logical"
operations.\footnote{In Boolean logic, formulas using only the implication
$\Rightarrow$ and $0$ suffice to define\ all of the $16$ binary logical
operations on subsets. Although beyond the scope of this paper, there are, by
the author's calculations, $134$ binary logical operations on partitions
definable just with formulas using only the implication $\Rightarrow$ and
$0$.} For our purposes here, we will settle for being able to define the
sixteen binary logical operations on partitions that correspond to the sixteen
logical binary subset operations. But which binary operations form a complete
set that suffices to define all those sixteen operations?

\subsection{Conjunctive normal form in partition logic}

The four operations, the join, meet, implication, and nand, suffice to define
the sixteen binary logical partition operations by using the partition version
of conjunctive normal form--which, in turn, is based on the following result.

\begin{lemma}
For any subsets $A,B\subseteq U\times U$, $\operatorname*{int}\left[  A\cap
B\right]  =\operatorname*{int}\left[  \operatorname*{int}\left(  A\right)
\cap\operatorname*{int}\left(  B\right)  \right]  $.
\end{lemma}

\noindent Proof: Since $\operatorname*{int}\left(  A\right)  \subseteq A$ and
$\operatorname*{int}\left(  B\right)  \subseteq B$, $\operatorname*{int}%
\left(  A\right)  \cap\operatorname*{int}\left(  B\right)  \subseteq A\cap B$
and thus $\operatorname*{int}\left[  \operatorname*{int}\left(  A\right)
\cap\operatorname*{int}\left(  B\right)  \right]  \subseteq\operatorname*{int}%
\left[  A\cap B\right]  $. Conversely, $A\cap B\subseteq A,B$ so
$\operatorname*{int}\left(  A\cap B\right)  \subseteq\operatorname*{int}%
\left(  A\right)  \cap\operatorname*{int}\left(  B\right)  $ and since
$\operatorname*{int}\left(  A\cap B\right)  $ is open, $\operatorname*{int}%
\left[  A\cap B\right]  \subseteq\operatorname*{int}\left[
\operatorname*{int}\left(  A\right)  \cap\operatorname*{int}\left(  B\right)
\right]  $. $\blacksquare$

In the treatment of the $16$ subsets defined from four basic regions $S\cap T
$, $S\cap T^{c}$, $S^{c}\cap T$, and $S^{c}\cap T^{c}$, we were in effect
using disjunctive normal form to define the $15$ non-empty subsets by taking
the unions of the $15$ combinations of those four basic regions. But the above
lemma shows that the conjunctive normal form will be more useful in partition
logic (since the corresponding result for the union and the interior operator
does not hold).

In the subset version of the conjunctive normal form, the $15$ non-universal
subsets are obtained by taking the intersections of $15$ combinations of the
four regions: $S\cup T$, $S\cup T^{c}$, $S^{c}\cup T$, and $S^{c}\cup T^{c}$.
Taking $S=\operatorname*{dit}\left(  \sigma\right)  $ and
$T=\operatorname*{dit}\left(  \tau\right)  $, the interiors of these four
basic "conjuncts" are, respectively, the dit sets of: $\sigma\vee\tau$,
$\tau\Rightarrow\sigma$, $\sigma\Rightarrow\tau$, and $\sigma\mid\tau$. By
expressing each of the $15$ non-universal subsets of $U\times U$ in
conjunctive normal form, applying the interior operator, and then using the
lemma to distribute the interior operator across the intersections, we express
each of the $15$ partition operations (aside from the constant $1$) as a meet
of some combination of the join $\sigma\vee\tau$, the implications
$\tau\Rightarrow\sigma$ and $\sigma\Rightarrow\tau$, and the nand $\sigma
\mid\tau$. The constant operation $1$ can be obtained using just the
implication $\sigma\Rightarrow\sigma$ or $\tau\Rightarrow\tau$. These results
and some other easy reductions are given in the following tables. In the first
table, the interior of the subset of $U\times U$ in the first column yields
the dit set of the binary operation given in the second column.\footnote{For
notation, we have followed, for the most part, Church
\citeyear{church:ML}%
.}

\begin{center}%
\begin{tabular}
[c]{|c|c|}\hline
$15$ regions Conjunctive Normal Form & Binary operation on
partitions\\\hline\hline
$\left(  S\cup T\right)  \cap\left(  S^{c}\cup T\right)  \cap\left(  S\cup
T^{c}\right)  \cap\left(  S^{c}\cup T^{c}\right)  $ & $0$\\\hline
$\left(  S^{c}\cup T\right)  \cap\left(  S\cup T^{c}\right)  \cap\left(
S^{c}\cup T^{c}\right)  $ & $\sigma\overline{\vee}\tau=\lnot\sigma\wedge
\lnot\tau$\\\hline
$\left(  S\cup T\right)  \cap\left(  S\cup T^{c}\right)  \cap\left(  S^{c}\cup
T^{c}\right)  $ & $\tau\nLeftarrow\sigma=\sigma\wedge\lnot\tau$\\\hline
$\left(  S\cup T^{c}\right)  \cap\left(  S^{c}\cup T^{c}\right)  $ &
$\lnot\tau=\tau\Rightarrow0$\\\hline
$\left(  S\cup T\right)  \cap\left(  S^{c}\cup T\right)  \cap\left(  S^{c}\cup
T^{c}\right)  $ & $\sigma\nLeftarrow\tau=\lnot\sigma\wedge\tau$\\\hline
$\left(  S^{c}\cup T\right)  \cap\left(  S^{c}\cup T^{c}\right)  $ &
$\lnot\sigma=\sigma\Rightarrow0$\\\hline
$\left(  S\cup T\right)  \cap\left(  S^{c}\cup T^{c}\right)  $ &
$\sigma\not \equiv \tau$\\\hline
$S^{c}\cup T^{c}$ & $\sigma\mid\tau$\\\hline
$\left(  S\cup T\right)  \cap\left(  S^{c}\cup T\right)  \cap\left(  S\cup
T^{c}\right)  $ & $\sigma\wedge\tau$\\\hline
$\left(  S^{c}\cup T\right)  \cap\left(  T^{c}\cup S\right)  $ & $\sigma
\equiv\tau$\\\hline
$\left(  S\cup T\right)  \cap\left(  S\cup T^{c}\right)  $ & $\sigma$\\\hline
$S\cup T^{c}$ & $\tau\Rightarrow\sigma$\\\hline
$\left(  S\cup T\right)  \cap\left(  S^{c}\cup T\right)  $ & $\tau$\\\hline
$S^{c}\cup T$ & $\sigma\Rightarrow\tau$\\\hline
$S\cup T$ & $\sigma\vee\tau$\\\hline
\end{tabular}

Interior of column 1 gives partition operation in column 2
\end{center}

Using the lemma, the interior is distributed across the intersections of the
subset CNF to give the partition CNF in the following table.

\begin{center}%
\begin{tabular}
[c]{|c|c|}\hline
Binary operation & Partition CNF for $15$ binary operations\\\hline\hline
$0$ & \multicolumn{1}{|l|}{$=\left(  \sigma\vee\tau\right)  \wedge\left(
\sigma\Rightarrow\tau\right)  \wedge\left(  \tau\Rightarrow\sigma\right)
\wedge\left(  \sigma\mid\tau\right)  $}\\\hline
$\sigma\overline{\vee}\tau=\lnot\sigma\wedge\lnot\tau$ &
\multicolumn{1}{|l|}{$=\left(  \sigma\Rightarrow\tau\right)  \wedge\left(
\tau\Rightarrow\sigma\right)  \wedge\left(  \sigma\mid\tau\right)  $}\\\hline
$\tau\nLeftarrow\sigma=\sigma\wedge\lnot\tau$ & \multicolumn{1}{|l|}{$=\left(
\sigma\vee\tau\right)  \wedge\left(  \tau\Rightarrow\sigma\right)
\wedge\left(  \sigma\mid\tau\right)  $}\\\hline
$\lnot\tau=\tau\Rightarrow0$ & \multicolumn{1}{|l|}{$=\left(  \tau
\Rightarrow\sigma\right)  \wedge\left(  \sigma\mid\tau\right)  $}\\\hline
$\sigma\nLeftarrow\tau=\lnot\sigma\wedge\tau$ & \multicolumn{1}{|l|}{$=\left(
\sigma\vee\tau\right)  \wedge\left(  \sigma\Rightarrow\tau\right)
\wedge\left(  \sigma\mid\tau\right)  $}\\\hline
$\lnot\sigma=\sigma\Rightarrow0$ & \multicolumn{1}{|l|}{$=\left(
\sigma\Rightarrow\tau\right)  \wedge\left(  \sigma\mid\tau\right)  $}\\\hline
$\sigma\not \equiv \tau$ & \multicolumn{1}{|l|}{$=\left(  \sigma\vee
\tau\right)  \wedge\left(  \sigma\mid\tau\right)  $}\\\hline
$\sigma\mid\tau$ & \multicolumn{1}{|l|}{$=\sigma\mid\tau$}\\\hline
$\sigma\wedge\tau$ & \multicolumn{1}{|l|}{$=\left(  \sigma\vee\tau\right)
\wedge\left(  \sigma\Rightarrow\tau\right)  \wedge\left(  \tau\Rightarrow
\sigma\right)  $}\\\hline
$\sigma\equiv\tau$ & \multicolumn{1}{|l|}{$=\left(  \sigma\Rightarrow
\tau\right)  \wedge\left(  \tau\Rightarrow\sigma\right)  $}\\\hline
$\sigma$ & \multicolumn{1}{|l|}{$=\left(  \sigma\vee\tau\right)  \wedge\left(
\tau\Rightarrow\sigma\right)  $}\\\hline
$\tau\Rightarrow\sigma$ & \multicolumn{1}{|l|}{$=\tau\Rightarrow\sigma$%
}\\\hline
$\tau$ & \multicolumn{1}{|l|}{$=\left(  \sigma\vee\tau\right)  \wedge\left(
\sigma\Rightarrow\tau\right)  $}\\\hline
$\sigma\Rightarrow\tau$ & \multicolumn{1}{|l|}{$=\sigma\Rightarrow\tau$%
}\\\hline
$\sigma\vee\tau$ & \multicolumn{1}{|l|}{$=\sigma\vee\tau$}\\\hline
\end{tabular}

Distributing interior across intersections gives partition CNF
\end{center}

The CNF identities shows that the $15$ functions, and thus all their further
combinations, could be defined in terms of the four primitive operations of
join, meet, implication, and nand.\footnote{There are other combinations which
can be taken as primitive since the \textit{inequivalence}, \textit{symmetric
difference, exclusive-or}, or \textit{xor} $\sigma\not \equiv \tau$ can be
used to define the nand operation: $\left(  \left(  \sigma\vee\tau\right)
\Rightarrow\left(  \sigma\not \equiv \tau\right)  \right)  =\sigma\mid\tau$.}

The fourteen non-zero operations occur in natural pairs: $\Rightarrow$ and
$\nRightarrow$, $\Leftarrow$ and $\nLeftarrow$, $\equiv$ and $\not \equiv $,
$\vee$ and $\overline{\vee}$, and $\wedge$ and $\mid$ in addition to $\sigma$
and $\lnot\sigma$, and $\tau$ and $\lnot\tau$. Except in the case of the join
$\vee$ (and, of course, $\sigma$ and $\tau$), the second operation in the pair
is not the negation of the first. The relationship is not negation but
$0$-orthogonality. The pairs of formulas $\sigma\Rightarrow\tau$ and
$\sigma\nRightarrow\tau$ (and similarly for the other pairs) are
$0$-orthogonal; if one is non-zero, the other must be zero. Later we see a
different pairing of the operations by duality.

\subsection{Partition algebra $\Pi(U)$ on $U$ and its dual $\Pi(U)^{op}$}

The partition lattice of all partitions on $U$ with the top $1$ and bottom $0
$ enriched with the binary operations of implication and nand is the
\textit{partition algebra} $\Pi\left(  U\right)  $ of $U$. It plays the role
for partition logic that the Boolean algebra $\mathcal{P}(U)$ of all subsets
of $U$ plays in ordinary subset logic. Dualization in classical propositional
logic--when expressed in terms of subsets--amounts to reformulating the
operations as operations on subset complements. But since the complements are
in the same Boolean algebra, Boolean duality can be expressed as a theorem
about a Boolean algebra. We have defined the lattice of partitions (sets of
disjoint and mutually exhaustive non-empty subsets of a set) as being
isomorphic to the lattice of partition relations $\mathcal{O}\left(  U\times
U\right)  $ on $U\times U$ (anti-reflexive, symmetric, and anti-transitive
relations). The complement of a partition relation is an equivalence relation
(reflexive, symmetric, and transitive relations) which is not an element in
the same lattice. Hence in partition logic, duality is naturally expressed as
a relationship between the partition algebra $\Pi(U)$ (seen as the algebra of
partition relations) and the dual algebra $\Pi(U)^{op}$ of equivalence relations.

Given a formula $\varphi$ in Boolean logic, the dual formula $\varphi^{d}$ is
obtained by interchanging $0$ and $1$, and by interchanging each of the
following pairs of operations: $\Rightarrow$ and $\nLeftarrow$, $\vee$ and
$\wedge$, $\equiv$ and $\not \equiv $, $\Leftarrow$ and $\nRightarrow$, and
$\overline{\vee}$ and $\mid$, while leaving the atomic variables and negation
$\lnot$ unchanged \cite[p. 106]{church:ML}. In partition logic, we use exactly
the same definition of dualization except that the atomic variables (and
constants) will now stand for equivalence relations rather than partitions so
we will indicate this by adding the superscript "$d$" to the atomic variables.
However the partition formulas may be assumed to involve only $\vee$, $\wedge
$, $\Rightarrow$, and $\mid$ along with $0$ and $1$. Hence the dual of modus
ponens $\varphi=\left(  \sigma\wedge\left(  \sigma\Rightarrow\tau\right)
\right)  \Rightarrow\tau$ is $\varphi^{d}=\left(  \sigma^{d}\vee\left(
\sigma^{d}\nLeftarrow\tau^{d}\right)  \nLeftarrow\tau^{d}\right)  $. The
\textit{converse non-implication} $\nLeftarrow$ (to use Church's terminology)
is the difference operation \cite[p. 201]{law:sfm}, i.e., $\sigma
^{d}\nLeftarrow\tau^{d}$ is the result of subtracting $\sigma^{d}$ from
$\tau^{d}$ so it might otherwise be symbolized as $\tau^{d}-\sigma^{d}$ (or
$\tau^{d}\backslash\sigma^{d}$). Then the dual to the modus ponens formula
would be: $\varphi^{d}=\tau^{d}-\left(  \sigma^{d}\vee\left(  \tau^{d}%
-\sigma^{d}\right)  \right)  $. This, incidentally, is the formula that would
have been compared to modus ponens $\left(  \sigma\wedge\left(  \sigma
\Rightarrow\tau\right)  \right)  \Rightarrow\tau$ in subset and intuitionistic
logic if the lattice of partitions had been written upside down instead of
just comparing the \textit{same} formulas in subset, intuitionistic, and
partition logic (a benefit of writing that lattice right side up). Similarly
the \textit{non-implication} $\sigma^{d}\nRightarrow\tau^{d}$, dual to the
reverse implication $\sigma\Leftarrow\tau$, might otherwise by symbolized as
the difference $\sigma^{d}-\tau^{d}$ (or $\sigma^{d}\backslash\tau^{d}$). The
difference $\tau^{d}-\sigma^{d}$ and nor $\sigma^{d}\overline{\vee}\tau^{d}$
will be taken as primitive operations on equivalence relations. The difference
and nor operations on partitions (as opposed to equivalence relations) are not
primitive: $\sigma\nLeftarrow\tau=\tau\wedge\lnot\sigma$ and $\sigma
\overline{\vee}\tau=\lnot\sigma\wedge\lnot\tau$. The equivalence and
inequivalence operations on partitions are also not taken as primitive:
$\sigma\equiv\tau=\left(  \sigma\Rightarrow\tau\right)  \wedge\left(
\tau\Rightarrow\sigma\right)  $ and $\sigma\not \equiv \tau=\left(  \sigma
\vee\tau\right)  \wedge\left(  \sigma\mid\tau\right)  $.

The process of dualization is reversible. Starting with a formula $\varphi
^{d}$ with superscript "$d$" on all atomic variables (to indicate they refer
to equivalence relations instead of partitions), dualizing means making the
same interchanges of operation symbols and constants, and erasing the "$d$"
superscripts so that the dual of the dual is the original formula.

We have used the lower case Greek letters $\pi$, $\sigma$, $...$ to stand for
set-of-blocks partitions while the corresponding partition relations were the
dit sets $\operatorname*{dit}\left(  \pi\right)  $, $\operatorname*{dit}%
\left(  \sigma\right)  $, $...$. The Greek letters with the superscript "$d$"
stand for equivalence relations which take the form $\operatorname*{indit}%
\left(  \pi\right)  $, $\operatorname*{indit}\left(  \sigma\right)  $, $...$.
Thus atomic variables such as $\pi$ dualize to $\pi^{d}$ and would be
interpreted as denoting indit sets $\operatorname*{indit}\left(  \pi\right)  $.

The operations of the dual algebra $\Pi(U)^{op}$ of equivalence relations on
$U$ could be defined directly but it is more convenient to define them using
duality from the partition operations.

\begin{enumerate}
\item The top of the dual algebra, usually denoted $\widehat{1}$, is
$0^{d}=\operatorname*{indit}\left(  0\right)  =U\times U$, the universal
equivalence relation that identifies everything. The bottom of the dual
algebra, usually denoted $\widehat{0}$, is $1^{d}=\operatorname*{indit}\left(
1\right)  =\Delta$, the diagonal where each element of $U$ is only identified
with itself.

\item Given any equivalence relations $\operatorname*{indit}\left(
\pi\right)  $ and $\operatorname*{indit}\left(  \sigma\right)  $ on $U$, their
meet $\wedge$ is defined via duality as the indit set of the join of the two
corresponding partitions: $\operatorname*{indit}\left(  \pi\right)
\wedge\operatorname*{indit}\left(  \sigma\right)  =\operatorname*{indit}%
\left(  \pi\vee\sigma\right)  =\operatorname*{indit}\left(  \pi\right)
\cap\operatorname*{indit}\left(  \sigma\right)  $. Using the superscript-$d$
notation, this is: $\pi^{d}\wedge\sigma^{d}=\left(  \pi\vee\sigma\right)
^{d}=\operatorname*{indit}\left(  \pi\vee\sigma\right)  $.

\item Similarly the join of two equivalence relations is defined via duality
as: $\operatorname*{indit}\left(  \pi\right)  \vee\operatorname*{indit}\left(
\sigma\right)  =\operatorname*{indit}\left(  \pi\wedge\sigma\right)
=\overline{\left\{  \operatorname*{indit}\left(  \pi\right)  \cup
\operatorname*{indit}\left(  \sigma\right)  \right\}  }$, so that using the
superscript-$d$ notation: $\pi^{d}\vee\sigma^{d}=\left(  \pi\wedge
\sigma\right)  ^{d}=\operatorname*{indit}\left(  \pi\wedge\sigma\right)  $.

\item The difference of two equivalence relations is defined via duality as:
$\operatorname*{indit}\left(  \pi\right)  -\operatorname*{indit}\left(
\sigma\right)  =\operatorname*{indit}\left(  \sigma\Rightarrow\pi\right)
=\overline{\left\{  \operatorname*{dit}\left(  \sigma\right)  \cap
\operatorname*{indit}\left(  \pi\right)  \right\}  }$, which in the other
notation is: $\pi^{d}-\sigma^{d}=\left(  \sigma\Rightarrow\pi\right)
^{d}=\operatorname*{indit}\left(  \sigma\Rightarrow\pi\right)  $.

\item The not-or or nor operation on equivalence relations is defined via
duality as: $\operatorname*{indit}\left(  \pi\right)  \overline{\vee
}\operatorname*{indit}\left(  \sigma\right)  =\operatorname*{indit}\left(
\pi\mid\sigma\right)  =\overline{\left\{  \left(  \operatorname*{indit}\left(
\pi\right)  \cup\operatorname*{indit}\left(  \sigma\right)  \right)
^{c}\right\}  }$, which gives: $\pi^{d}\overline{\vee}\sigma^{d}=\left(
\pi\mid\sigma\right)  ^{d}=\operatorname*{indit}\left(  \pi\mid\sigma\right)
$.
\end{enumerate}

\noindent That completes the definition of the dual algebra $\Pi(U)^{op}$ of
equivalence relations on $U$ with the top $\widehat{1}$, bottom $\widehat{0}$,
and the four primitive operations of meet, join, difference, and nor.

The dualization operation $\varphi\longmapsto\varphi^{d}$ is a purely
syntactic operation on formulas, but in the partition algebra $\Pi(U)$ and
equivalence relation algebra $\Pi(U)^{op}$ we reason semantically about
partitions and equivalence relations on $U$. Given a compound formula
$\varphi$ in the language of the partition algebra, it would be interpreted by
interpreting its atomic variables as denoting partitions on $U$ and then
applying the partition operations (join, meet, implication, and nand) to
arrive at an interpretation of $\varphi$. Such an interpretation automatically
supplies an interpretation of the dual formula $\varphi^{d}$. If $\alpha$ was
an atomic variable of $\varphi$ and was interpreted as denoting a partition on
$U$, then $\alpha^{d}$ is interpreted as denoting the equivalence relation
$\operatorname*{indit}\left(  \alpha\right)  $. Then the equivalence relation
operations (meet, join, difference, and nor) are applied to arrive at an
equivalence relation interpretation of the formula $\varphi^{d}$. The
relationship between the two interpretations is very simple.

\begin{proposition}
$\varphi^{d}=\operatorname*{indit}\left(  \varphi\right)  $.
\end{proposition}

\noindent Proof: The proof uses induction over the complexity of the formulas
[where complexity is defined in the standard way as in propositional logic
\cite{fit:il}]. If $\varphi$ is one of the constants $0$ or $1$, then the
proposition holds since: $0^{d}=\widehat{1}=\operatorname*{indit}\left(
0\right)  $ and $1^{d}=\widehat{0}=\operatorname*{indit}\left(  1\right)  $.
If $\varphi=\alpha$ is atomic, then it is true by the definition: $\sigma
^{d}=\operatorname*{indit}\left(  \sigma\right)  $. If $\varphi$ is a compound
formula then the main connective in $\varphi$ is one of the four primitive
partition operations and the main connective in $\varphi^{d}$ is one of the
four primitive equivalence relation operations. Consider the case:
$\varphi=\pi\wedge\sigma$ so that $\varphi^{d}=\pi^{d}\vee\sigma^{d}$. By the
induction hypothesis, $\pi^{d}=\operatorname*{indit}\left(  \pi\right)  $ and
$\sigma^{d}=\operatorname*{indit}\left(  \sigma\right)  $, and by the
definition of the equivalence relation join: $\varphi^{d}=\pi^{d}\vee
\sigma^{d}=\operatorname*{indit}\left(  \pi\right)  \vee\operatorname*{indit}%
\left(  \sigma\right)  =\overline{\left\{  \operatorname*{indit}\left(
\pi\right)  \cup\operatorname*{indit}\left(  \sigma\right)  \right\}
}=\operatorname*{indit}\left(  \varphi\right)  $. The other three cases
proceed in a similar manner. $\blacksquare$

\begin{corollary}
The map $\varphi\longmapsto\operatorname*{indit}\left(  \varphi\right)  $ is a
dual-isomorphism: $\Pi(U)\rightarrow\Pi(U)^{op}$ between the partition algebra
and the dual equivalence relation algebra.
\end{corollary}

\noindent Proof: Clearly the mapping is a set isomorphism since each partition
$\varphi$ on $U$ is uniquely determined by its dit set $\operatorname*{dit}%
\left(  \varphi\right)  $, and thus by its complement $\operatorname*{indit}%
\left(  \varphi\right)  $. By "dual-isomorphism," we mean that each operation
in the partition algebra is mapped to the dual operation in the equivalence
relation algebra. Suppose $\varphi=\sigma\Rightarrow\pi$ so that $\varphi
^{d}=\pi^{d}-\sigma^{d}$. By the proposition, this means that
$\operatorname*{indit}\left(  \varphi\right)  =\operatorname*{indit}\left(
\pi\right)  -\operatorname*{indit}\left(  \sigma\right)  $ (where we must be
careful to note that "$-$" is the difference operation on equivalence
relations which is the closure of the set-difference operation
$\operatorname*{indit}\left(  \pi\right)  \cap\operatorname*{indit}\left(
\sigma\right)  ^{c}$ on subsets of $U\times U$) so that $\varphi
\longmapsto\operatorname*{indit}\left(  \varphi\right)  $ maps the partition
operation of implication to the equivalence relation operation of difference.
The other operations are treated in a similar manner. $\blacksquare$

The previous result $\operatorname*{int}\left[  A\cap B\right]
=\operatorname*{int}\left[  \operatorname*{int}\left(  A\right)
\cap\operatorname*{int}\left(  B\right)  \right]  $ for $A,B\subseteq U\times
U$ could also be expressed using the closure operation as $\overline{\left[
A\cup B\right]  }=\overline{\left[  \overline{A}\cup\overline{B}\right]  }$
and thus the conjunctive normal form treatment of the $15$ binary operations
on partitions in terms of the operations of $\vee$, $\wedge$, $\Rightarrow$,
and $\mid$ dualizes to the disjunctive normal form treatment of the $15$
(dual) binary operations on equivalence relations in terms of the respective
dual operations $\wedge$, $\vee$, $-$, and $\overline{\vee}$, which are the
primitive operations in the algebra of equivalence relations $\Pi(U)^{op}$.

The previous two tables giving the CNF treatment of the $15$ partition
operations dualize to give two similar tables for the DNF treatment of the
$15$ non-zero operations on equivalence relations. In the following table, let
$S^{\prime}=\operatorname*{indit}\left(  \sigma\right)  $ and $T^{\prime
}=\operatorname*{indit}\left(  \tau\right)  $ where $\left(  {}\right)  ^{c}$
is complementation in $U\times U$. In the table, we have also taken the
liberty of writing the "converse non-implication" operation as the difference
operation on both equivalence relations and partitions: $\tau^{d}-\sigma
^{d}=\sigma^{d}\nLeftarrow\tau^{d}$ and $\tau-\sigma=\sigma\nLeftarrow
\tau=\tau\wedge\lnot\sigma$.

\begin{center}%
\begin{tabular}
[c]{|c|c|c|}\hline
$15$ regions Disjunctive Normal Form & Bin. op. on eq. rel. & Dual
to\\\hline\hline
$S^{\prime c}\cap T^{\prime c}$ & $\sigma^{d}\overline{\vee}\tau^{d}$ &
$\sigma\mid\tau$\\\hline
$S^{\prime}\cap T^{\prime c}$ & $\sigma^{d}-\tau^{d}$ & $\tau\Rightarrow
\sigma$\\\hline
$\left(  S^{\prime}\cap T^{\prime c}\right)  \cup\left(  S^{\prime c}\cap
T^{\prime c}\right)  $ & $\lnot\tau^{d}$ & $\lnot\tau$\\\hline
$S^{\prime c}\cap T^{\prime}$ & $\tau^{d}-\sigma^{d}$ & $\sigma\Rightarrow
\tau$\\\hline
$\left(  S^{\prime c}\cap T^{\prime}\right)  \cup\left(  S^{\prime c}\cap
T^{\prime c}\right)  $ & $\lnot\sigma^{d}$ & $\lnot\sigma$\\\hline
$\left(  S^{\prime c}\cap T^{\prime}\right)  \cup\left(  S^{\prime}\cap
T^{\prime c}\right)  $ & $\sigma^{d}\not \equiv \tau^{d}$ & $\sigma\equiv\tau
$\\\hline
$\left(  S^{\prime c}\cap T^{\prime}\right)  \cup\left(  S^{\prime c}\cap
T^{\prime c}\right)  \cup\left(  S^{\prime}\cap T^{\prime c}\right)  $ &
$\sigma^{d}\mid\tau^{d}$ & $\sigma\overline{\vee}\tau$\\\hline
$S^{\prime}\cap T^{\prime}$ & $\sigma^{d}\wedge\tau^{d}$ & $\sigma\vee\tau
$\\\hline
$\left(  S^{\prime}\cap T^{\prime}\right)  \cup\left(  S^{\prime c}\cap
T^{\prime c}\right)  $ & $\sigma^{d}\equiv\tau^{d}$ & $\sigma\not \equiv \tau
$\\\hline
$\left(  S^{\prime}\cap T^{\prime}\right)  \cup\left(  S^{\prime}\cap
T^{\prime c}\right)  $ & $\sigma^{d}$ & $\sigma$\\\hline
$\left(  S^{\prime}\cap T^{\prime}\right)  \cup\left(  S^{\prime}\cap
T^{\prime c}\right)  \cup\left(  S^{\prime c}\cap T^{\prime c}\right)  $ &
$\tau^{d}\Rightarrow\sigma^{d}$ & $\sigma-\tau$\\\hline
$\left(  S^{\prime}\cap T^{\prime}\right)  \cup\left(  S^{\prime c}\cap
T^{\prime}\right)  $ & $\tau^{d}$ & $\tau$\\\hline
$\left(  S^{\prime c}\cap T^{\prime}\right)  \cup\left(  S^{\prime c}\cap
T^{\prime c}\right)  \cup\left(  S^{\prime}\cap T^{\prime}\right)  $ &
$\sigma^{d}\Rightarrow\tau^{d}$ & $\tau-\sigma$\\\hline
$\left(  S^{\prime}\cap T^{\prime}\right)  \cup\left(  S^{\prime}\cap
T^{\prime c}\right)  \cup\left(  S^{\prime c}\cap T^{\prime}\right)  $ &
$\sigma^{d}\vee\tau^{d}$ & $\sigma\wedge\tau$\\\hline
$\left(  S^{\prime}\cap T^{\prime}\right)  \cup\left(  S^{\prime}\cap
T^{\prime c}\right)  \cup\left(  S^{\prime c}\cap T^{\prime}\right)
\cup\left(  S^{\prime c}\cap T^{\prime c}\right)  $ & $\widehat{1}$ &
$0$\\\hline
\end{tabular}

Closure of column $1$ gives equivalence relation binary operation in column
$2$
\end{center}

For instance, the CNF expression for the partition inequivalence or symmetric
difference is: $\sigma\not \equiv \tau=\left(  \sigma\vee\tau\right)
\wedge\left(  \sigma\mid\tau\right)  $ so that:%

\begin{align*}
\operatorname*{dit}\left(  \sigma\not \equiv \tau\right)   &
=\operatorname*{int}\left[  \operatorname*{int}\left(  \operatorname*{dit}%
\left(  \sigma\right)  \cup\operatorname*{dit}\left(  \tau\right)  \right)
\cap\operatorname*{int}\left(  \operatorname*{dit}\left(  \sigma\right)
^{c}\cup\operatorname*{dit}\left(  \tau\right)  ^{c}\right)  \right] \\
&  =\operatorname*{int}\left[  \left(  \operatorname*{dit}\left(
\sigma\right)  \cup\operatorname*{dit}\left(  \tau\right)  \right)
\cap\left(  \operatorname*{dit}\left(  \sigma\right)  ^{c}\cup
\operatorname*{dit}\left(  \tau\right)  ^{c}\right)  \right]  \text{.}%
\end{align*}

\noindent Taking complements yields:%

\begin{align*}
\operatorname*{indit}\left(  \sigma\not \equiv \tau\right)   &  =\overline
{\left[  \left(  \operatorname*{indit}\left(  \sigma\right)  \cap
\operatorname*{indit}\left(  \tau\right)  \right)  \cup\left(
\operatorname*{indit}\left(  \sigma\right)  ^{c}\cap\operatorname*{indit}%
\left(  \tau\right)  ^{c}\right)  \right]  }\\
&  =\overline{\left[  \overline{\left(  \operatorname*{indit}\left(
\sigma\right)  \cap\operatorname*{indit}\left(  \tau\right)  \right)  }%
\cup\overline{\left(  \operatorname*{indit}\left(  \sigma\right)  ^{c}%
\cap\operatorname*{indit}\left(  \tau\right)  ^{c}\right)  }\right]  }\\
&  =\overline{\left[  \left(  \sigma^{d}\wedge\tau^{d}\right)  \cup\left(
\sigma^{d}\overline{\vee}\tau^{d}\right)  \right]  }\\
&  =\left(  \sigma^{d}\wedge\tau^{d}\right)  \vee\left(  \sigma^{d}%
\overline{\vee}\tau^{d}\right) \\
&  =\sigma^{d}\equiv\tau^{d}\text{.}%
\end{align*}

\noindent Thus the equivalence $\sigma^{d}\equiv\tau^{d}$ of equivalence
relations has the disjunctive normal form: $\sigma^{d}\equiv\tau^{d}=\left(
\sigma^{d}\wedge\tau^{d}\right)  \vee\left(  \sigma^{d}\overline{\vee}\tau
^{d}\right)  $ in the "dual" logic of equivalence relations. The disjunctive
normal forms for the $15$ operations on equivalence relations is given in the
following table.

\begin{center}%
\begin{tabular}
[c]{|c|c|}\hline
Binary operation & Equivalence relation DNF for $15$ binary
operations\\\hline\hline
\multicolumn{1}{|c|}{$\sigma^{d}\overline{\vee}\tau^{d}$} &
\multicolumn{1}{|l|}{$=\sigma^{d}\overline{\vee}\tau^{d}$}\\\hline
\multicolumn{1}{|c|}{$\sigma^{d}-\tau^{d}$} & \multicolumn{1}{|l|}{$=\sigma
^{d}-\tau^{d}$}\\\hline
\multicolumn{1}{|c|}{$\lnot\tau^{d}$} & \multicolumn{1}{|l|}{$=\left(
\sigma^{d}-\tau^{d}\right)  \vee\left(  \sigma^{d}\overline{\vee}\tau
^{d}\right)  $}\\\hline
\multicolumn{1}{|c|}{$\tau^{d}-\sigma^{d}$} & \multicolumn{1}{|l|}{$=\tau
^{d}-\sigma^{d}$}\\\hline
\multicolumn{1}{|c|}{$\lnot\sigma^{d}$} & \multicolumn{1}{|l|}{$=\left(
\tau^{d}-\sigma^{d}\right)  \vee\left(  \sigma^{d}\overline{\vee}\tau
^{d}\right)  $}\\\hline
\multicolumn{1}{|c|}{$\sigma^{d}\not \equiv \tau^{d}$} &
\multicolumn{1}{|l|}{$=\left(  \tau^{d}-\sigma^{d}\right)  \vee\left(
\sigma^{d}-\tau^{d}\right)  $}\\\hline
\multicolumn{1}{|c|}{$\sigma^{d}\mid\tau^{d}$} &
\multicolumn{1}{|l|}{$=\left(  \tau^{d}-\sigma^{d}\right)  \vee\left(
\sigma^{d}\overline{\vee}\tau^{d}\right)  \vee\left(  \sigma^{d}-\tau
^{d}\right)  $}\\\hline
\multicolumn{1}{|c|}{$\sigma^{d}\wedge\tau^{d}$} &
\multicolumn{1}{|l|}{$=\sigma^{d}\wedge\tau^{d}$}\\\hline
\multicolumn{1}{|c|}{$\sigma^{d}\equiv\tau^{d}$} &
\multicolumn{1}{|l|}{$=\left(  \sigma^{d}\wedge\tau^{d}\right)  \vee\left(
\sigma^{d}\overline{\vee}\tau^{d}\right)  $}\\\hline
\multicolumn{1}{|c|}{$\sigma^{d}$} & \multicolumn{1}{|l|}{$=\left(  \sigma
^{d}\wedge\tau^{d}\right)  \vee\left(  \sigma^{d}-\tau^{d}\right)  $}\\\hline
\multicolumn{1}{|c|}{$\tau^{d}\Rightarrow\sigma^{d}$} &
\multicolumn{1}{|l|}{$=\left(  \sigma^{d}\wedge\tau^{d}\right)  \vee\left(
\sigma^{d}-\tau^{d}\right)  \vee\left(  \sigma^{d}\overline{\vee}\tau
^{d}\right)  $}\\\hline
\multicolumn{1}{|c|}{$\tau^{d}$} & \multicolumn{1}{|l|}{$=\left(  \sigma
^{d}\wedge\tau^{d}\right)  \vee\left(  \tau^{d}-\sigma^{d}\right)  $}\\\hline
\multicolumn{1}{|c|}{$\sigma^{d}\Rightarrow\tau^{d}$} &
\multicolumn{1}{|l|}{$=\left(  \sigma^{d}\wedge\tau^{d}\right)  \vee\left(
\tau^{d}-\sigma^{d}\right)  \vee\left(  \sigma^{d}\overline{\vee}\tau
^{d}\right)  $}\\\hline
\multicolumn{1}{|c|}{$\sigma^{d}\vee\tau^{d}$} &
\multicolumn{1}{|l|}{$=\left(  \sigma^{d}\wedge\tau^{d}\right)  \vee\left(
\sigma^{d}-\tau^{d}\right)  \vee\left(  \tau^{d}-\sigma^{d}\right)  $}\\\hline
\multicolumn{1}{|c|}{$\widehat{1}$} & \multicolumn{1}{|l|}{$=\left(
\sigma^{d}\wedge\tau^{d}\right)  \vee\left(  \sigma^{d}-\tau^{d}\right)
\vee\left(  \tau^{d}-\sigma^{d}\right)  \vee\left(  \sigma^{d}\overline{\vee
}\tau^{d}\right)  $}\\\hline
\end{tabular}

Distributing closure across unions gives equivalence relation DNF
\end{center}

These DNF identities give the expression of the non-primitive binary
operations on equivalence relations, e.g., $\equiv$, $\not \equiv $, $\mid
$,and $\Rightarrow$, in terms of the primitive operations. The constant
$\widehat{0}$ may be defined as $\sigma^{d}-\sigma^{d}$ dual to the definition
of $1$ as $\sigma\Rightarrow\sigma$.

In referring to the dual logic of equivalence relations, we must keep distinct
the different notions of duality. Partition logic is dual to subset logic in
the sense of the category-theoretic duality between monomorphisms and
epimorphisms (or between subsets and quotient sets). But equivalence relation
logic is only dual to partition logic in the sense of
complementation--analogous to the duality between Heyting algebras and
co-Heyting algebras, or between open subsets and closed subsets of a
topological space. Since the complement of an open set is a closed set that is
not necessarily open, complementation-duality for partition logic and
intuitionistic propositional logic is a duality between two types of algebras
(partition algebras and equivalence relation algebras in the one case and
Heyting and co-Heyting algebras in the other case). But the complement of a
subset is another subset so complementation-duality for subset logic is a
duality within a Boolean algebra.

\subsection{Subset and partition tautologies}

For present purposes, we may take the formulas of classical propositional
logic (i.e., subset logic) as using the binary operations of $\vee$, $\wedge$,
$\Rightarrow$, and $\mid$ along with the constants $0$ and $1$ so that we have
exactly the same well-formed formulas in subset logic and partition logic. A
\textit{truth-table tautology} is a formula that always evaluates to $1$ in
the Boolean algebra $\mathcal{P}(1)$ regardless of the assignments of $0$ and
$1$ to the atomic variables. A \textit{subset tautology} is a formula that
always evaluates to $1$ (the universe set $U$) in the Boolean algebra
$\mathcal{P}(U)$ regardless of the subsets assigned to the atomic variables.

\begin{proposition}
Subset tautologies = truth-table tautologies.
\end{proposition}

\noindent Proof: Since $1$ is a special case of $U$, all subset tautologies
are truth-table tautologies. To see that all truth-table tautologies are
subset tautologies, we reinterpret the columns in a truth table as giving the
truth and falsity of the set-membership statements for the logical subset
operations. For instance, if $\sigma$ and $\tau$ are interpreted as subsets of
a universe $U$, then $\sigma\Rightarrow\tau$ is the subset $\sigma^{c}\cup
\tau$ so that $u\in\left(  \sigma\Rightarrow\tau\right)  $ is assigned a $1$
(representing truth) iff $u\not \in \sigma$ or $u\in\tau$, i.e., $u\in\sigma$
is assigned a $0$ or $u\in\tau$ is assigned a $1$. The truth table for the
modus ponens tautology would then be rewritten with the new column labels.

\begin{center}%
\begin{tabular}
[c]{|c|c|c|c|c|}\hline
$u\in\sigma$ & $u\in\tau$ & $u\in\left(  \sigma\Rightarrow\tau\right)  $ &
$u\in\left(  \sigma\wedge\left(  \sigma\Rightarrow\tau\right)  \right)  $ &
$u\in\left(  \left(  \sigma\wedge\left(  \sigma\Rightarrow\tau\right)
\right)  \Rightarrow\tau\right)  $\\\hline\hline
$1$ & $1$ & $1$ & $1$ & $1$\\\hline
$1$ & $0$ & $0$ & $0$ & $1$\\\hline
$0$ & $1$ & $1$ & $0$ & $1$\\\hline
$0$ & $0$ & $1$ & $0$ & $1$\\\hline
\end{tabular}

Truth table reinterpreted as giving set-membership conditions
\end{center}

\noindent Given a formula that is a truth-table tautology, we assign any
subsets of any given $U$ to the atomic variables. For any specific element
$u\in U$, it will be either a member or non-member of each of those subsets so
one of the cases (rows) in the reinterpreted truth table will apply. Since the
truth table operations give the set-membership conditions for the
corresponding logical subset operations (illustrated in the example above),
and since, by assumption, the final column under the formula is all $1$'s, the
arbitrary element $u$ of $U$ is a member of the subset which interprets the
truth-table tautologous formula. Hence that subset must be $U$ and thus the
formula is a subset tautology. $\blacksquare$

A \textit{partition tautology} is a formula that always evaluates to $1$ (the
discrete partition) in the partition algebra $\Pi(U)$ regardless of the
partitions assigned to the atomic variables.\footnote{Needless to say, the
constants $0$ and $1$ are always assigned the bottom and top, respectively, in
any evaluation or interpretation of a formula in either $\mathcal{P}(U)$ or
$\Pi(U)$.} It is also useful to define a \textit{weak partition tautology} as
a formula that never evaluates to $0$ (the indiscrete partition) regardless of
the partitions assigned to the atomic variables. Of course, any partition
tautology is a weak partition tautology. Moreover, it is easily seen that:

\begin{proposition}
$\varphi$ is a weak partition tautology iff $\lnot\lnot\varphi$ is a partition tautology.
\end{proposition}

An immediate question is the relationship of partition tautologies and weak
partition tautologies to the classical subset tautologies as well as to the
valid formulas of intuitionistic propositional logic (where formulas are
assumed to be written in the same language).

There is a sense in which results in partition logic can be trivially seen as
a generalization of results in ordinary subset logic. This\textit{\ reduction
principle} is based on the observation that any partition logic result holding
for all $U$ will hold when restricted to any two element universe $\left\vert
U\right\vert =2$. There is an isomorphism between the partition algebra
$\Pi\left(  2\right)  $ on the two-element set and the Boolean algebra
$\mathcal{P}(1)$ on the one-element set. There are only two partitions, the
bottom $0$ and top $1$ on $U$ where $\left\vert U\right\vert =2$. Moreover,
the partition operations of join, meet, implication, and nand in this special
case satisfy the truth tables for the corresponding Boolean operations on
subsets (using $0$ and $1$ in the usual manner in the truth tables). For
instance, in $\Pi(U)$ where $\left\vert U\right\vert =2$, we can only
substitute $0$ or $1$ for the atomic variables in $\sigma\Rightarrow\tau$. The
result is $0$ in the case where $\sigma=1$ and $\tau=0$, and the result in $1$
in the other three cases. But that is just the truth table for the Boolean
implication operation in $\mathcal{P}(1)$. Similarly for the other operations
so there is a BA isomorphism: $\Pi(2)\cong\mathcal{P}\left(  1\right)  $.
Hence if a partition logic result holds for all $U$, then it holds for a
two-element $U$ where the partition operations on on the partitions $0$ and
$1$ are isomorphic to the Boolean operations on the subsets $0$ and $1$ (where
$0$ and $1$ in the Boolean case stand for the null subset and the universe set
of a one-element universe). But if a formula always evaluates to $1$ on the
one-element universe, i.e., if it is a truth-table tautology, then it is a
subset tautology. Thus we might say that partition logic restricted to a
two-element universe is Boolean logic:

\begin{center}
$\Pi(2)\cong\mathcal{P}\left(  1\right)  $.

Reduction Principle
\end{center}

For instance, if $\varphi$ is a weak partition tautology, e.g., $\varphi
=\sigma\vee\lnot\sigma$, then it will never evaluate to $0$ in any $\Pi(U)$
where it is always assumed $\left\vert U\right\vert \geq2$. For $\left\vert
U\right\vert =2$, there are only two partitions $0$ and $1$, so never
evaluating to $0$ means always evaluating to the partition $1$. By the
reduction principle, the Boolean operations in $\mathcal{P}\left(  1\right)  $
would always evaluate to the subset $1$ so the formula is a truth-table
tautology and thus a subset tautology. This proves the following proposition.

\begin{proposition}
All weak partition tautologies are subset tautologies. $\blacksquare$
\end{proposition}

\begin{corollary}
All partition tautologies are subset tautologies. $\blacksquare$
\end{corollary}

The converse is not true with Peirce's law, $\left(  \left(  \sigma
\Rightarrow\pi\right)  \Rightarrow\sigma\right)  \Rightarrow\sigma$,
accumulation, $\sigma\Rightarrow\left(  \pi\Rightarrow\left(  \sigma\wedge
\pi\right)  \right)  $, and distributivity, $\left(  \left(  \pi\vee
\sigma\right)  \wedge\left(  \pi\vee\tau\right)  \right)  \Rightarrow\left(
\pi\vee\left(  \sigma\wedge\tau\right)  \right)  $, being examples of subset
tautologies that are not partition tautologies.

There is no inclusion either way between partition tautologies and the valid
formulas of intuitionistic propositional logic. In view of the complex nature
of the partition meet, it is not surprising that a formula such as the
accumulation formula, $\sigma\Rightarrow\left(  \pi\Rightarrow\left(
\pi\wedge\sigma\right)  \right)  $, is valid in both Boolean and
intuitionistic logic but not in partition logic. The ("non-weak") law of
excluded middle, $\sigma\vee\lnot\sigma$, is a weak partition tautology, and
the weak law of excluded middle, $\lnot\sigma\vee\lnot\lnot\sigma$, is a
("non-weak") partition tautology that is not intuitionistically valid.

In the dual algebra $\Pi(U)^{op}$ of equivalence relations, the bottom is the
smallest equivalence relation $\widehat{0}=\Delta=\operatorname*{indit}\left(
1\right)  $ containing only the diagonal pairs $\left(  u,u\right)  $. Dual to
the notion of a partition tautology is the notion of an \textit{equivalence
relation contradiction} which is a formula (with the atomic variables written
with the "$d$" superscript) that always evaluates to the bottom $\widehat{0}%
=\Delta=1^{d}$ of $\Pi(U)^{op}$ regardless of the equivalence relations
substituted for the atomic variables. Similarly, a formula (with the atomic
variables written with the "$d$" superscript) is a \textit{weak equivalence
relation contradiction} if it never evaluates to the top $\widehat{1}=U\times
U=0^{d}$ of $\Pi(U)^{op}$. We then have the following duality theorem.

\begin{proposition}
[Principle of duality for partition logic]Given a formula $\varphi$, $\varphi$
is a (weak) partition tautology iff $\varphi^{d}$ is a (resp. weak)
equivalence relation contradiction.
\end{proposition}

\noindent Proof: Using the dual-isomorphism $\Pi(U)\rightarrow\Pi(U)^{op}$, a
partition formula $\varphi$ evaluates to the top of $\Pi(U)$, i.e.,
$\operatorname*{dit}\left(  \varphi\right)  =\operatorname*{dit}\left(
1\right)  =U\times U-\Delta$ when any partitions are substituted for the
atomic variables of $\varphi$ iff $\varphi^{d}$ evaluates to the bottom of
$\Pi(U)^{op}$, i.e., $1^{d}=\operatorname*{dit}\left(  1\right)
^{c}=\operatorname*{indit}\left(  1\right)  =\widehat{0}=\Delta$, when any
equivalence relations are substituted for the atomic variables of $\varphi
^{d}$. Similarly for the weak notions. $\blacksquare$

Using the reduction principle, restricting the above proposition and its
related concepts to $\left\vert U\right\vert =2$ would yield the usual Boolean
duality principle \cite[p. 107]{church:ML} that $\varphi$ is a truth-table
tautology iff $\varphi^{d}$ is a truth-table contradiction (where the weak or
"non-weak" notions coincide in the Boolean case and where $\Pi(2)\cong%
\mathcal{P}\left(  1\right)  \cong\Pi(2)^{op}$).

In the Boolean case, if a formula $\varphi$ is not a subset tautology, then
there is a non-empty universe set $U$ and an assignment of subsets of $U$ to
the atomic variables of $\varphi$ so that $\varphi$ does not evaluate to $1 $
(the universe set $U$). Such a model showing that $\varphi$ is not a tautology
is called a \textit{countermodel} for $\varphi$. It is a remarkable aspect of
the Boolean logic of subsets that to determine subset tautologies, it suffices
to restrict the universe set $U$ to a one-element set and to consider only
truth-table tautologies. If $\varphi$ has a subset countermodel, then it has a
countermodel using the subsets of a one-element set.\footnote{The drawback of
the coextensiveness of truth-table tautologies and general subset tautologies
is that it is easier to take the specific propositional interpretation of
Boolean logic as the whole subject matter of the logic.}

The analogous question can be posed for partition logic. Is there a finite
number $n$ so that if $\varphi$ always evaluates to $1$ for any partitions on
$U$ with $\left\vert U\right\vert \leq n$, then $\varphi$ is a partition
tautology? For instance, if $\varphi$ is not a partition tautology and is also
not a subset tautology, then it suffices to take $n=2$ since $\Pi
(2)\cong\mathcal{P}\left(  1\right)  $ so a truth-table countermodel in
$\mathcal{P}\left(  1\right)  $ also provides a partition countermodel in
$\Pi(2)$. Hence the question is only open for formulas $\varphi$ which are
subset tautologies but not partition tautologies. A standard device answers
this question in the negative.

\begin{proposition}
There is no fixed $n$ such that if any $\varphi$ has no partition countermodel
on any universe $U$ with $\left\vert U\right\vert \leq n$, then $\varphi$ has
no partition countermodel, i.e., is a partition tautology.
\end{proposition}

\noindent Proof: Consider any fixed $n\geq2$. We use the standard device of a
"universal disjunction of equations" \cite[p. 316]{grat:glt} to construct a
formula $\omega_{n}$ that evaluates to $1$ for any substitutions of partitions
on $U$ with $\left\vert U\right\vert \leq n$ and yet the formula is not a
partition tautology. Let $B_{n}$ be the Bell number, the number of partitions
on a set $U$ with $\left\vert U\right\vert =n$. Take the atomic variables to
be $\pi_{i}$ for $i=0,1,...,B_{n}$ so that there are $B_{n}+1$ atomic
variables. Let $\omega_{n}$ be the join of all the equivalences between
distinct atomic variables:

\begin{center}
$\omega_{n}=%
{\textstyle\bigvee}
\left\{  \pi_{i}\equiv\pi_{j}:0\leq i<j\leq B_{n}\right\}  $.
\end{center}

\noindent Then for any substitution of partitions on $U$ where $\left\vert
U\right\vert \leq n$ for the atomic variables, there is, by the pigeonhole
principle, some "disjunct" $\pi_{i}\equiv\pi_{j}=\left(  \pi_{i}\Rightarrow
\pi_{j}\right)  \wedge\left(  \pi_{j}\Rightarrow\pi_{i}\right)  $ which has
the same partition substituted for the two variables so the disjunct evaluates
to $1$ and thus the join $\omega_{n}$ evaluates to $1$. Thus $\omega_{n}$
evaluates to $1$ for any substitutions of partitions on any $U$ where
$\left\vert U\right\vert \leq n$. To see that $\omega_{n}$ is not a partition
tautology, take $U=\left\{  0,1,...,B_{n}\right\}  $ and let $\pi_{i}$ be the
atomic partition which has $i$ as a singleton and all the other elements of
$U$ as a block, i.e., $\pi_{i}=\left\{  \left\{  0,1,...,i-1,i+1,...,B_{n}%
\right\}  ,\left\{  i\right\}  \right\}  $. Then $\pi_{i}\Rightarrow\pi
_{j}=\pi_{j}$ and $\pi_{j}\wedge\pi_{i}=0$ so that $\omega_{n}=0$ for that
substitution and thus $\omega_{n}$ is not even a weak partition tautology.
$\blacksquare$

To see that the $\omega_{n}$ are subset tautologies, consider $n=2$, so that
$B_{2}=2$ and $\omega_{2}=\left(  \pi_{0}\equiv\pi_{1}\right)  \vee\left(
\pi_{0}\equiv\pi_{2}\right)  \vee\left(  \pi_{1}\equiv\pi_{2}\right)  $. Thus
$\omega_{2}$ is a truth-table tautology and hence any larger join $\omega_{n}$
for $n>2$ is also a truth-table tautology and thus a subset tautology.

There is no upper bound $n$ so that if any formula has a countermodel, then it
has a countermodel with $\left\vert U\right\vert \leq n$. However, it seems
likely to the author that if a partition formula has a countermodel, then it
has a finite countermodel (i.e., the finite model property) but that question
remains open.

\subsection{Boolean subalgebras $\mathcal{B}_{\pi}$ of $\Pi(U)$ for any
partition $\pi$}

In any Heyting algebra, the elements of the form $\lnot\sigma=\sigma
\Rightarrow0$ for some $\sigma$ are the \textit{regular} elements. They form a
Boolean algebra but it is not a subalgebra since the join of two regular
elements is not necessarily regular (so one must take the double negation of
the join to have the Boolean algebra join). In the topological interpretation,
the regular elements of the Heyting algebra of open subsets are the regular
open sets (the regular open sets are obtained as the interior of the closure
of a subset) and the union of two regular open subsets is open but not
necessarily regular open.

Following the analogy, we define a partition as being $\pi$\textit{-regular}
if it can be obtained as the implication $\sigma\Rightarrow\pi$ for some
partitions $\sigma$ and $\pi$. Intuitively, a $\pi$-regular partition is like
$\pi$ except that some blocks may have been discretized. Let

\begin{center}
$\mathcal{B}_{\pi}=\left\{  \sigma\Rightarrow\pi:\text{ for some }\sigma\in
\Pi(U)\right\}  $
\end{center}

\noindent be the subset of $\pi$-regular partitions with the induced partial
ordering of refinement. The top is still $1$ but the bottom is $\pi
=1\Rightarrow\pi$ itself. The implication partition $\sigma\Rightarrow\pi$ can
be interpreted as a Boolean probe for containment between blocks. If
$B\subseteq C$ for some $C\in\sigma$, then the probe finds containment and
this is indicated by setting the $\pi$-block $B$ locally equal to $1$, i.e.,
by discretizing $B$, and otherwise $B$ stays locally like $0$, i.e., stays as
a whole block (or "mini-blob") $B$. Whenever the refinement relation
$\sigma\preceq\pi$ holds, then all the non-singleton blocks $B\in\pi$ are
discretized in $\sigma\Rightarrow\pi$ (and the singleton blocks are already
discrete) so that $\sigma\Rightarrow\pi=1$ (and vice-versa).%

\begin{center}
\includegraphics[
natheight=150.996307bp,
natwidth=326.027405bp,
height=154.375pt,
width=331pt
]%
{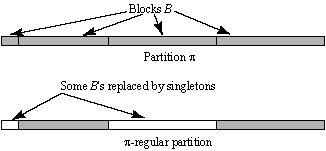}%
\\
Figure 6: $B$-slots in $\pi$-regular partition
\end{center}

The partition operations of meet and join operate on the blocks of $\pi
$-regular partitions in a completely Boolean manner. Since every $\pi$-regular
partition is like $\pi$ except that some blocks may be set locally to $1$
while the others remain locally like $0$, the meet of two $\pi$-regular
partitions, say $\sigma\Rightarrow\pi$ and $\tau\Rightarrow\pi$, will have no
interaction between distinct $\pi$-blocks. Each block of the meet will be
"truth-functionally" determined by whatever is in the $B$-slot of the two
constituents. If either of the $B$'s remains locally equal to $0$, then the
whole block $B$ fills the $B$-slot of the meet, i.e., $B$ is locally equal to
$0$ in the meet $\left(  \sigma\Rightarrow\pi\right)  \wedge\left(
\tau\Rightarrow\pi\right)  =\overset{\pi}{\lnot}\tau\wedge\overset{\pi}{\lnot
}\sigma$. But if both $B$'s were discretized in the constituents, i.e., both
are set locally to $1$, then the blocks in that $B$-slot of the meet are the
singletons from $B$, i.e., the discretized $B$ or $B$ set locally to $1$. That
local pattern of $0$'s and $1$'s is precisely the truth table for the Boolean meet.

If $\pi_{ns}$ is the set of non-singleton blocks of the partition $\pi$, then
the $\pi$-regular partitions are in one-to-one correspondence with the subsets
of $\pi_{ns}$, each of which can be represented by its characteristic function
$\chi:\pi_{ns}\rightarrow2=\left\{  0,1\right\}  $ which takes each
non-singleton block of $\pi$ to its local assignment. Thus for a $\pi$-regular
partition with the form $\sigma\Rightarrow\pi$, $\chi\left(  \sigma
\Rightarrow\pi\right)  :\pi_{ns}\rightarrow2$ takes a non-singleton block
$B\in\pi$ to $1$ if $B$ is discretized in $\sigma\Rightarrow\pi$ and otherwise
to $0$.

\noindent The argument just given shows that the characteristic function for
the meet of two $\pi$-regular partitions is obtained by the component-wise
Boolean meets of "conjuncts":

\begin{center}
$\chi\left(  \left(  \sigma\Rightarrow\pi\right)  \wedge\left(  \tau
\Rightarrow\pi\right)  \right)  =\chi\left(  \sigma\Rightarrow\pi\right)
\wedge\chi\left(  \tau\Rightarrow\pi\right)  $.
\end{center}

In a similar manner, the blocks in the join of two $\pi$-regular partitions,
$\sigma\Rightarrow\pi$ and $\tau\Rightarrow\pi$, would be the intersections of
what is in the $B$-slots. If $B$ was discretized (set locally to $1$) in
either of the constituents, then $B$ would be discretized in the join $\left(
\tau\Rightarrow\pi\right)  \vee\left(  \sigma\Rightarrow\pi\right)
=\overset{\pi}{\lnot}\tau\vee\overset{\pi}{\lnot}\sigma$ (since the
intersection of a discretized $B$ with a whole $B$ is still the discretized
$B$). But if both $B$'s were still whole (set locally to $0$) then their
intersection would still be the whole block $B$. This pattern of $0$'s and $1
$'s is precisely the truth table for the Boolean join or disjunction. In terms
of the characteristic functions of local assignments:

\begin{center}
$\chi\left(  \left(  \tau\Rightarrow\pi\right)  \vee\left(  \sigma
\Rightarrow\pi\right)  \right)  =\chi\left(  \tau\Rightarrow\pi\right)
\vee\chi\left(  \sigma\Rightarrow\pi\right)  $.
\end{center}

For the implication $\left(  \sigma\Rightarrow\pi\right)  \Rightarrow\left(
\tau\Rightarrow\pi\right)  $ between two $\pi$-regular partitions, the result
would have $B$ remaining whole, i.e., being set to $0$, only in the case where
$B$ was whole in the consequent partition $\tau\Rightarrow\pi$ but discretized
in the antecedent partition $\sigma\Rightarrow\pi$; otherwise $B$ is
discretized, i.e., set to $1$. This pattern of $0$'s and $1$'s is precisely
the truth table for the ordinary Boolean implication. In terms of the
characteristic functions:

\begin{center}
$\chi\left(  \left(  \sigma\Rightarrow\pi\right)  \Rightarrow\left(
\tau\Rightarrow\pi\right)  \right)  =\chi\left(  \sigma\Rightarrow\pi\right)
\Rightarrow\chi\left(  \tau\Rightarrow\pi\right)  $.
\end{center}

To show that $\mathcal{B}_{\pi}$ is a Boolean algebra, we must define negation
inside of $\mathcal{B}_{\pi}$. The negation of a $\pi$-regular element
$\sigma\Rightarrow\pi$ would be its implication to the bottom element which in
$\mathcal{B}_{\pi}$ is $\pi$ itself. Thus the negation of $\sigma
\Rightarrow\pi=\overset{\pi}{\lnot}\sigma$ is just the iterated implication:
$\left(  \sigma\Rightarrow\pi\right)  \Rightarrow\pi=\overset{\pi}{\lnot
}\overset{\pi}{\lnot}\sigma$, the double $\pi$-negation. It is easily seen
that this just "flips" the $B$-slots to the opposite state. The $B$'s set
(locally) to $1$ in $\sigma\Rightarrow\pi$ are flipped back to (locally) $0$
in $\left(  \sigma\Rightarrow\pi\right)  \Rightarrow\pi$, and the $B$'s left
whole in $\sigma\Rightarrow\pi$ are flipped to $1$ or discretized in $\left(
\sigma\Rightarrow\pi\right)  \Rightarrow\pi$. This pattern of $0$'s and $1$'s
is just the truth table for the Boolean negation. In terms of the
characteristic functions,

\begin{center}
$\chi\left(  \left(  \sigma\Rightarrow\pi\right)  \Rightarrow\pi\right)
=\lnot\chi\left(  \sigma\Rightarrow\pi\right)  $.
\end{center}

Thus it is easily seen that the set of $\pi$-regular elements $\mathcal{B}%
_{\pi}$ is a Boolean algebra, called the \textit{Boolean core} of the upper
interval $\left[  \pi,1\right]  =\left\{  \sigma\in\Pi(U):\pi\preceq
\sigma\preceq1\right\}  $, since it is isomorphic to the powerset Boolean
algebra $\mathcal{P}\left(  \pi_{ns}\right)  $ of the set $\pi_{ns}$ (when the
subsets are represented by their characteristic functions).

\begin{proposition}
$\mathcal{B}_{\pi}\cong\mathcal{P}\left(  \pi_{ns}\right)  $. $\blacksquare$
\end{proposition}

We previously saw that the partition lattice $\Pi\left(  U\right)  $ could be
represented by the lattice of open subsets $\operatorname*{dit}\left(
\pi\right)  $ of the product $U\times U$ (when taken as a closure space). The
representation of the partition lattice by the open subsets of the closure
space $U\times U$ continues to hold when the lattice is enriched with the
implication and nand operations.

\begin{center}
\fbox{$\Pi(U)\cong O\left(  U\times U\right)  $}

Representation of algebra of partitions $\Pi(U)$

as the algebra of open subsets $O\left(  U\times U\right)  $
\end{center}

\noindent Now we can see the dual representation of the Boolean algebra of
subsets $\mathcal{P}(U)$ by a certain Boolean algebra modeled using partition
operations. Start with the dual constructions of subsets
$\operatorname*{indit}\left(  \pi\right)  $ of the product $U\times U$ and the
partitions $\Delta\left(  S\right)  $ on the coproduct $U%
{\textstyle\biguplus}
U$. For the dual representation of $\mathcal{P}(U)$ we consider the partition
algebra $\Pi(U%
{\textstyle\biguplus}
U)$ on the coproduct and the Boolean core $\mathcal{B}_{\Delta}$, or
$\mathcal{B}_{\Delta}\left(  U%
{\textstyle\biguplus}
U\right)  $ to make the underlying universe explicit, associated with the
diagonal partition $\Delta\left(  U\right)  =\Delta\left(  \emptyset
^{c}\right)  $ consisting of all the pairs $\left\{  u,u^{\ast}\right\}  $ for
$u\in U$. Just as we previous took the complement of $\operatorname*{indit}%
\left(  \pi\right)  $ to arrive at the partition relations
$\operatorname*{dit}\left(  \pi\right)  =\operatorname*{indit}\left(
\pi\right)  ^{c}$, and we now consider the $\Delta$-complements
$\overset{\Delta}{\lnot}\Delta\left(  S\right)  =\Delta\left(  S^{c}\right)  $
which are the \textit{subset corelations}. The $\Delta$-regular partitions of
$\Pi(U%
{\textstyle\biguplus}
U)$ are precisely the subset corelations $\Delta\left(  S^{c}\right)  $. The
subset corelation $\Delta\left(  S^{c}\right)  $ locally assigns $\left\{
u,u^{\ast}\right\}  \in\Delta$ to $1$ (i.e., discretizes it) if $u\in S$ and
locally assigns $\left\{  u,u^{\ast}\right\}  \in\Delta$\ to $0 $ (i.e.,
leaves it whole) if $u\in S^{c}$. Rather than associate each partition $\pi$
with the partition relation $\operatorname*{dit}\left(  \pi\right)  $ on the
product $U\times U$, we now associate each subset $S\in\mathcal{P}\left(
U\right)  $ with the subset corelation $\Delta\left(  S^{c}\right)  $ on the
coproduct $U%
{\textstyle\biguplus}
U$ to get the dual representation:

\begin{center}
\fbox{$\mathcal{P}(U)\cong\mathcal{B}_{\Delta}\left(  U%
{\textstyle\biguplus}
U\right)  $}

Dual representation of the Boolean algebra of subsets $\mathcal{P}(U)$

as the BA of subset corelations $\mathcal{B}_{\Delta}\left(  U%
{\textstyle\biguplus}
U\right)  $.
\end{center}

The universe sets $U$ are assumed to have two or more elements to avoid the
degenerate case of a singleton universe where $0=1$, i.e., the indiscrete and
discrete partitions are the same. But in partitions $\pi$, singleton blocks
cannot be avoided and the same problem emerges locally. For a singleton block
$B$, being locally like $0$ (i.e., remaining whole) and being locally like $1$
(being discretized) are the same. Hence the singletons of $\pi$ play no role
in the Boolean algebras $\mathcal{B}_{\pi} $.

We previously saw another Boolean algebra $\mathcal{B}\left(  \pi\right)  $
associated with every partition $\pi$ on a set $U$, and the singletons will
play a role in connecting the two BAs. For each partition $\pi$ on $U$,
$\mathcal{B}\left(  \pi\right)  \subseteq\mathcal{P}(U)$ is the complete
subalgebra generated by the blocks of $\pi$ as the atoms so that all the
elements of $\mathcal{B}\left(  \pi\right)  $ are formed as the arbitrary
unions and intersections of blocks of $\pi$. Since each element of
$\mathcal{B}\left(  \pi\right)  $ is the union of a set of blocks of $\pi$, it
is isomorphic to the powerset BA of the set of blocks that make up $\pi$,
i.e., $\mathcal{B}\left(  \pi\right)  \cong\mathcal{P}\left(  \pi\right)  $.
Since $\mathcal{B}_{\pi}\cong\mathcal{P}\left(  \pi_{ns}\right)  $ is
isomorphic to the powerset BA of the set of non-singleton blocks of $\pi$, and
since the introduction of each singleton $\left\{  u\right\}  $ will have the
effect of doubling the elements of $\mathcal{P}\left(  \pi_{ns}\right)  $
(with or without the singleton), we can reach $\mathcal{P}\left(  \pi\right)
$ from $\mathcal{P}\left(  \pi_{ns}\right)  $ by taking the direct product
with the two element BA $2$ for each singleton in $\pi$. Thus we have the
following result which relates the two BAs associated with each partition
$\pi$.

\begin{proposition}
$\mathcal{B}\left(  \pi\right)  \cong\mathcal{B}_{\pi}\times\prod
\limits_{\left\{  u\right\}  \in\pi}2$. $\blacksquare$
\end{proposition}

\subsection{Transforming subset tautologies into partition tautologies}

Unlike the case of the Boolean algebra of regular elements in a Heyting
algebra, the Boolean core $\mathcal{B}_{\pi}$ is a \textit{sub}algebra of the
partition algebra $\Pi(U)$ for the "Boolean" operations of join, meet, and
implication (but not nand), i.e., the Boolean operations in $\mathcal{B}_{\pi
}$ are the partition operations from the partition algebra $\Pi(U)$. The BA
$\mathcal{B}_{\pi}$ even has the same top $1$ as the partition algebra; only
the bottoms are different, i.e., $\pi$ in $\mathcal{B}_{\pi}$ and $0$ in
$\Pi(U)$.

Since the Boolean core $\mathcal{B}_{\pi}$ of the interval $\left[
\pi,1\right]  $ and the whole partition algebra $\Pi(U)$ have the same top $1$
and the same operations of join, meet, and implication, we immediately have a
way to transform any subset tautology into a partition tautology. But we must
be careful about the connectives used in the subset tautology. The partition
operations of the join, meet, and implication are the same as the Boolean
operations in the Boolean core $\mathcal{B}_{\pi}$. But the negation in that
BA is not the partition negation $\lnot$ but the $\pi$-negation $\overset{\pi
}{\lnot}$. Similarly, the nand operation in the Boolean algebra $\mathcal{B}%
_{\pi}$ is not the partition nand $\mid$ but the $\pi$\textit{-nand} defined
by the ternary partition operation: $\operatorname*{dit}\left(  \sigma
\mid_{\pi}\tau\right)  =\operatorname*{int}\left[  \operatorname*{indit}%
\left(  \sigma\right)  \cup\operatorname*{indit}\left(  \tau\right)
\cup\operatorname*{dit}\left(  \pi\right)  \right]  $ which agrees with the
usual nand when $\pi=0$. But the nand operation in the BA $\mathcal{B}_{\pi}$
can be defined in terms of the other BA operations so we may assume that the
subset tautology is written without a nand operation $\mid$. Similarly we may
assume that negations $\lnot\sigma$ are written as $\sigma\Rightarrow0$ so
that no negation signs $\lnot$ occur in the partition tautology.

Given any propositional formula using the connectives of $\vee$, $\wedge$,
$\Rightarrow$ and the constants of $0$ and $1$, its \textit{single }$\pi
$\textit{-negation transform} is obtained by replacing each atomic variable
$\sigma$ by its single $\pi$-negation $\overset{\pi}{\lnot}\sigma
=\sigma\Rightarrow\pi$ and by replacing the constant $0$ by $\pi$. The binary
operations $\vee$, $\wedge$, and $\Rightarrow$ as well as the constant $1$ all
remain the same. For instance, the single $\pi$-negation transform of the
excluded middle formula $\sigma\vee\lnot\sigma=\sigma\vee\left(
\sigma\Rightarrow0\right)  $ is the weak excluded middle formula for $\pi$-negation:

\begin{center}
$\left(  \sigma\Rightarrow\pi\right)  \vee\left(  \left(  \sigma\Rightarrow
\pi\right)  \Rightarrow\pi\right)  =\overset{\pi}{\lnot}\sigma\vee
\overset{\pi}{\lnot}\overset{\pi}{\lnot}\sigma$.
\end{center}

A formula that is a subset tautology will always evaluate to $1$ in a Boolean
algebra regardless of what elements of the Boolean algebra are assigned to the
atomic variables. The single $\pi$-negation transformation maps any formula
into a formula for an element of the Boolean core $\mathcal{B}_{\pi}$. If the
original formula with the atomic variables $\sigma$, $\tau$,... was a subset
tautology, then the single $\pi$-negation transform of the formula will
evaluate to $1$ in $\mathcal{B}_{\pi}$ for any partitions ($\pi$-regular or
not) assigned to the original atomic variables $\sigma$, $\tau$, ... with
$\pi$ fixed. But this is true for any $\pi$ so the single $\pi$-negation
transform of any subset tautology will evaluate to $1$ for any partitions
assigned to the atomic variables $\pi$, $\sigma$, $\tau$,.... Thus it is a
partition tautology.

\begin{proposition}
The single $\pi$-negation transform of any subset tautology is a partition
tautology. $\blacksquare$
\end{proposition}

For example, since the law of excluded middle, $\sigma\vee\lnot\sigma$, is a
subset tautology, its single $\pi$-negation transform, $\overset{\pi}{\lnot
}\sigma\vee\overset{\pi}{\lnot}\overset{\pi}{\lnot}\sigma$, is a partition
tautology. This particular example is also intuitively obvious since the
blocks $B$ that were not discretized in $\overset{\pi}{\lnot}\sigma$ are
discretized in the double $\pi$-negation $\overset{\pi}{\lnot}\overset{\pi
}{\lnot}\sigma$ so all the non-singleton blocks are discretized in
$\overset{\pi}{\lnot}\sigma\vee\overset{\pi}{\lnot}\overset{\pi}{\lnot}\sigma$
(and the singleton blocks were already "discretized") so it is a partition
tautology. This formula is also an example of a partition tautology that is
not a valid formula of intuitionistic logic (either for $\pi=0$ or in general).

We can similarly define the \textit{double }$\pi$\textit{-negation transform}
of a formula as the formula where each atomic variable $\sigma$ is replaced by
its double $\pi$-negation $\overset{\pi}{\lnot}\overset{\pi}{\lnot}\sigma$ and
by replacing the constant $0$ by $\pi$. By the same argument, the double $\pi
$-negation transform of any subset tautology is a partition tautology so there
are at least two ways to transform any classical subset tautology into a
partition tautology.

\begin{proposition}
The double $\pi$-negation transform of any subset tautology is a partition
tautology. $\blacksquare$
\end{proposition}

The double $\pi$-negation transform of excluded middle, $\sigma\vee\lnot
\sigma$, is the partition tautology $\overset{\pi}{\lnot}\overset{\pi}{\lnot
}\sigma\vee\overset{\pi}{\lnot}\overset{\pi}{\lnot}\overset{\pi}{\lnot}\sigma
$. Since the $\pi$-negation has the effect of flipping the $\pi$-blocks $B$
back and forth being locally equal to $0$ or $1$ (i.e., from being whole to
being discretized), it is clear that $\overset{\pi}{\lnot}\sigma
=\overset{\pi}{\lnot}\overset{\pi}{\lnot}\overset{\pi}{\lnot}\sigma$ so the
formula $\overset{\pi}{\lnot}\overset{\pi}{\lnot}\sigma\vee\overset{\pi
}{\lnot}\overset{\pi}{\lnot}\overset{\pi}{\lnot}\sigma$ simplifies to
$\overset{\pi}{\lnot}\overset{\pi}{\lnot}\sigma\vee\overset{\pi}{\lnot}\sigma$.

There is also a partition analogue of the G\"{o}del transform \cite{god:int}
that produces an intuitionistic validity from each subset tautology. For any
classical formula $\varphi$ in the language of $\vee$, $\wedge$, and
$\Rightarrow$ as well as $0$ and $1$, we define the \textit{G\"{o}del }$\pi
$\textit{-transform} $\varphi_{\pi}^{g}$ of the formula as follows:

\begin{itemize}
\item If $\varphi$ is atomic, then $\varphi_{\pi}^{g}=\varphi\vee\pi$; if
$\varphi=0$, then $\varphi_{\pi}^{g}=\pi$, and if $\varphi=1$, then
$\varphi_{\pi}^{g}=1$;

\item If $\varphi=\sigma\vee\tau$, then $\varphi_{\pi}^{g}=\sigma_{\pi}%
^{g}\vee\tau_{\pi}^{g}$;

\item If $\varphi=\sigma\Rightarrow\tau$, then $\varphi_{\pi}^{g}=\sigma_{\pi
}^{g}\Rightarrow\tau_{\pi}^{g}$; and

\item if $\varphi=\sigma\wedge\tau$, then $\varphi_{\pi}^{g}=\overset{\pi
}{\lnot}\overset{\pi}{\lnot}\sigma_{\pi}^{g}\wedge\overset{\pi}{\lnot
}\overset{\pi}{\lnot}\tau_{\pi}^{g}$.
\end{itemize}

When $\pi=0$, then we write $\varphi_{0}^{g}=\varphi^{g}$.

\begin{lemma}
$\varphi$ is a subset tautology iff $\varphi^{g}$ is a weak partition
tautology iff $\lnot\lnot\varphi^{g}$ is a partition tautology.
\end{lemma}

\noindent Proof: The idea of the proof is that the partition operations on the
G\"{o}del $0$-transform $\varphi^{g}$ mimic the Boolean $0,1$-operations on
$\varphi$ if we associate the partition interpretation $\sigma^{g}=0$ with the
Boolean $\sigma=0$ and $\sigma^{g}\not =0$ with the Boolean $\sigma=1$. We
proceed by induction over the complexity of the formula $\varphi$ where the
induction hypothesis is that: $\varphi=1$ in the Boolean case iff $\varphi
^{g}\not =0$ in the partition case, which could also be stated as: $\varphi=0$
in the Boolean case iff $\varphi^{g}=0$ in the partition case.

\begin{enumerate}
\item If $\varphi$ is atomic, the Boolean assignment $\varphi=0$ (the Boolean
truth value $0$) is associated with the partition assignment of $\varphi=0$
(the indiscrete partition) and for atomic $\varphi$, $\varphi=\varphi
\vee0=\varphi^{g}$ so the hypothesis holds in the base case.

\item For the join in the Boolean case, $\varphi=\sigma\vee\tau=1$ iff
$\sigma=1$ or $\tau=1$. In the partition case, $\varphi^{g}=\sigma^{g}\vee
\tau^{g}\not =0$ iff $\sigma^{g}\not =0$ or $\tau^{g}\not =0$, so by the
induction hypothesis, $\varphi=\sigma\vee\tau=1$ iff $\sigma=1$ or $\tau=1$
iff $\sigma^{g}\not =0$ or $\tau^{g}\not =0$ iff $\varphi^{g}=\sigma^{g}%
\vee\tau^{g}\not =0$.

\item For the implication in the Boolean case, $\varphi=\sigma\Rightarrow
\tau=0$ iff $\sigma=1$ and $\tau=0$. In the partition case, $\varphi
^{g}=\sigma^{g}\Rightarrow\tau^{g}=0$ iff $\sigma^{g}\not =0$ and $\tau^{g}%
=0$. Hence using the induction hypothesis, $\varphi=\sigma\Rightarrow\tau=1$
iff $\sigma=0$ or $\tau=1$ iff $\sigma^{g}=0$ or $\tau^{g}\not =0 $ iff
$\varphi^{g}=\sigma^{g}\Rightarrow\tau^{g}\not =0$.

\item For the meet in the Boolean case, $\varphi=\sigma\wedge\tau=1$ iff
$\sigma=1=\tau$. In the partition case, $\varphi^{g}=\lnot\lnot\sigma
^{g}\wedge\lnot\lnot\tau^{g}=1$ iff $\lnot\lnot\sigma^{g}=1=\lnot\lnot\tau
^{g}$ iff $\sigma^{g}\not =0\not =\tau^{g}$. By the induction hypothesis,
$\varphi=\sigma\wedge\tau=1$ iff $\sigma=1=\tau$ iff $\sigma^{g}\not =%
0\not =\tau^{g}$ iff $\varphi^{g}=\lnot\lnot\sigma^{g}\wedge\lnot\lnot\tau
^{g}=1$ iff $\varphi^{g}=\lnot\lnot\sigma^{g}\wedge\lnot\lnot\tau^{g}\not =0$.
\end{enumerate}

\noindent Thus $\varphi$ is a subset tautology iff under any Boolean
interpretation, $\varphi=1$ iff for any partition interpretation, $\varphi
^{g}\not =0$ iff $\varphi^{g}$ is a weak partition tautology iff $\lnot
\lnot\varphi^{g}$ is a partition tautology. $\blacksquare$

In this case of $\pi=0$, the negation $\lnot\sigma=\sigma\Rightarrow0$ is
unchanged and, for atomic variables $\varphi$, $\varphi\vee0=\varphi$ so
atomic variables are left unchanged in the G\"{o}del $0$-transform. Hence any
classical formula $\varphi$ expressed in the language of $\lnot$, $\vee$, and
$\Rightarrow$ (excluding the meet $\wedge$) would be unchanged by the
G\"{o}del $0$-transform.

\begin{corollary}
For any formula $\varphi$ in the language of $\lnot$, $\vee$, and
$\Rightarrow$ along with $0$ and $1$, $\varphi$ is a subset tautology iff
$\varphi$ is a weak partition tautology iff $\lnot\lnot\varphi$ is a partition tautology.
\end{corollary}

\noindent For instance, the excluded middle subset tautology $\sigma\vee
\lnot\sigma$ is a weak partition tautology and $\lnot\lnot\left(  \sigma
\vee\lnot\sigma\right)  $ is a partition tautology.

The lemma generalizes to any $\pi$ in the following form.

\begin{proposition}
$\varphi$ is a subset tautology iff $\overset{\pi}{\lnot}\overset{\pi}{\lnot
}\varphi_{\pi}^{g}$ is a partition tautology.
\end{proposition}

\noindent Proof: For any fixed partition $\pi$ on a universe set $U$, the
interpretation of the G\"{o}del $\pi$-transform $\varphi_{\pi}^{g}$ is in the
upper interval $\left[  \pi,1\right]  \subseteq\Pi(U)$. The key to the
generalization is the standard result that the upper interval $\left[
\pi,1\right]  $ can be represented as the product of the sets $\Pi(B)$ where
$B$ is a non-singleton block of $\pi$:

\begin{center}
$\left[  \pi,1\right]  \cong%
{\textstyle\prod}
\left\{  \Pi(B):B\in\pi\text{, }B\text{ non-singleton}\right\}  $%
.\footnote{Since the partition lattice is conventionally written upside down,
the usual result is stated in terms of the interval below $\pi$ \cite[p.
252]{grat:glt}.}
\end{center}

\noindent Once we establish that the G\"{o}del $\pi$-transform $\varphi_{\pi
}^{g}$ can be obtained, using the isomorphism, by computing the G\"{o}del
$0$-transform $\varphi^{g}$ "component-wise" in $\Pi(B)$, then we can apply
the lemma component-wise to obtain the result.

We use induction over the complexity of $\varphi$ with the hypothesis:
$\varphi_{\pi}^{g}$ has a block $B\in\pi$ iff $\varphi_{0}^{g}=\varphi^{g}$ is
equal to the zero $0_{B}$ of $\Pi(B)$.

\begin{enumerate}
\item Given a partition $\pi$ on $U$, any interpretation of an atomic
$\varphi$ as a partition on $U$ can be cut down to each non-singleton block
$B\in\pi$ to yield a partition on $B$. Then $\varphi_{\pi}^{g}=\varphi\vee\pi$
has a block $B\in\pi$ iff $\varphi_{0}^{g}=\varphi^{g}$ is equal to the zero
$0_{B}$ of $\Pi(B)$.

\item If $\varphi=\sigma\vee\tau$, then a block of $\varphi_{\pi}^{g}%
=\sigma_{\pi}^{g}\vee\tau_{\pi}^{g}$ is $B$ iff $B$ is a block of both
$\sigma_{\pi}^{g}$ and $\tau_{\pi}^{g}$ iff $\sigma^{g}=0_{B}=\tau^{g}$ in
$\Pi\left(  B\right)  $ iff $\varphi^{g}=\sigma^{g}\vee\tau^{g}=0_{B}$ in
$\Pi\left(  B\right)  $.

\item If $\varphi=\sigma\Rightarrow\tau$, then $\varphi_{\pi}^{g}=\sigma_{\pi
}^{g}\Rightarrow\tau_{\pi}^{g}$ has a block $B\in\pi$ iff $\sigma_{\pi}^{g}$
does not have the block $B$ and $\tau_{\pi}^{g}$ has the block $B$ iff
$\sigma^{g}$ is not equal to $0_{B}$ and $\tau^{g}$ is equal to $0_{B}$ in
$\Pi\left(  B\right)  $ iff $\varphi^{g}=\sigma^{g}\Rightarrow\tau^{g}=0_{B}$
in $\Pi\left(  B\right)  $.

\item If $\varphi=\sigma\wedge\tau$, then $\varphi_{\pi}^{g}=\overset{\pi
}{\lnot}\overset{\pi}{\lnot}\sigma_{\pi}^{g}\wedge\overset{\pi}{\lnot
}\overset{\pi}{\lnot}\tau_{\pi}^{g}$ has a block $B\in\pi$ iff both
$\sigma_{\pi}^{g}$ and $\tau_{\pi}^{g}$ have a block $B$ iff $\sigma^{g}%
=0_{B}=\tau^{g}$ in $\Pi(B)$ iff $\varphi^{g}=\lnot\lnot\sigma^{g}\wedge
\lnot\lnot\tau^{g}=0_{B}$ in $\Pi(B)$.
\end{enumerate}

\noindent Hence applying the lemma component-wise, $\varphi$ is a subset
tautology iff $\varphi^{g}$ never evaluates to $0_{B}$ in $\Pi(B)$ iff $B$ is
never a block of $\varphi_{\pi}^{g}$ iff every block $B\in\pi$ is discretized
in $\overset{\pi}{\lnot}\overset{\pi}{\lnot}\varphi_{\pi}^{g}$, i.e.,
$\overset{\pi}{\lnot}\overset{\pi}{\lnot}\varphi_{\pi}^{g}$ is a partition
tautology. $\blacksquare$

Thus the G\"{o}del $\pi$-transform of excluded middle $\varphi=\sigma
\vee\left(  \sigma\Rightarrow0\right)  =\sigma\vee\lnot\sigma$ is
$\varphi_{\pi}^{g}=\left(  \sigma\vee\pi\right)  \vee\left(  \left(
\sigma\vee\pi\right)  \Rightarrow\pi\right)  =\left(  \sigma\vee\pi\right)
\vee\overset{\pi}{\lnot}\left(  \sigma\vee\pi\right)  $ and $\overset{\pi
}{\lnot}\overset{\pi}{\lnot}\left[  \left(  \sigma\vee\pi\right)
\vee\overset{\pi}{\lnot}\left(  \sigma\vee\pi\right)  \right]  $ is a
partition tautology. Note that the single $\pi$-negation transform, the double
$\pi$-negation transform, and the G\"{o}del $\pi$-transform all gave different
formulas starting with the classical excluded middle tautology.

\subsection{Some partition results}

Before turning to the proof theory of partition logic, we might mention a few
interesting results. For many purposes, the structure of the partition
algebras $\Pi(U)$ is best analyzed by analyzing the upper intervals $\left[
\pi,1\right]  $ for any partition $\pi$. Partition lattices are the "standard"
examples of non-distributive lattices, but one can do much better than simply
say a partition lattice is non-distributive. The Boolean core of each interval
$\left[  \pi,1\right]  $ is, of course, distributive since it is a Boolean
algebra using the meet and join operations of the partition lattice. Moreover,
each partition in the interval $\left[  \pi,1\right]  $ distributes across the
Boolean core. To see this, note that one of these distributivity results is
essentially due to Oystein Ore. Ore
\citeyear{ore:ter}
did much of the path-breaking work on partitions. He defined two partitions as
being \textit{associable} if each block in their meet is a block in one (or
both) of the partitions.\footnote{Ore actually dealt with the join of
equivalence relations but we are using the opposite presentation.} Although
Ore did not consider $\pi$-regular partitions, any two $\pi$-regular
partitions are associable. He showed that any partition joined with the meet
of two associable partitions will distribute across the meet \cite[p.
585]{ore:ter}. Hence we have the following result for any partitions $\varphi
$, $\sigma$, $\tau$, and $\pi$.

\begin{lemma}
[Ore's distributivity theorem]$\varphi\vee\left(  \overset{\pi}{\lnot}%
\sigma\wedge\overset{\pi}{\lnot}\tau\right)  =\left(  \varphi\vee
\overset{\pi}{\lnot}\sigma\right)  \wedge\left(  \varphi\vee\overset{\pi
}{\lnot}\tau\right)  $.
\end{lemma}

Ore's theorem does not assume that $\varphi$ is in the interval $\left[
\pi,1\right]  $ but we can interchange join and meet if we restrict $\varphi$
to the interval.

\begin{lemma}
["Dual" to Ore's theorem]If $\varphi\in\left[  \pi,1\right]  $, then
$\varphi\wedge\left(  \overset{\pi}{\lnot}\sigma\vee\overset{\pi}{\lnot}%
\tau\right)  =\left(  \varphi\wedge\overset{\pi}{\lnot}\sigma\right)
\vee\left(  \varphi\wedge\overset{\pi}{\lnot}\tau\right)  $.
\end{lemma}

\begin{proposition}
[Distributivity over the Boolean core]If $\pi\preceq\varphi$,
\end{proposition}

\begin{align*}
\varphi\vee\left(  \overset{\pi}{\lnot}\sigma\wedge\overset{\pi}{\lnot}%
\tau\right)   &  =\left(  \varphi\vee\overset{\pi}{\lnot}\sigma\right)
\wedge\left(  \varphi\vee\overset{\pi}{\lnot}\tau\right) \\
\varphi\wedge\left(  \overset{\pi}{\lnot}\sigma\vee\overset{\pi}{\lnot}%
\tau\right)   &  =\left(  \varphi\wedge\overset{\pi}{\lnot}\sigma\right)
\vee\left(  \varphi\wedge\overset{\pi}{\lnot}\tau\right)  \text{.}%
\end{align*}

Distributivity over the Boolean core allows certain results based on subset
distributivity to be "transferred" to partitions. For instance, for arbitrary
subsets $S,T$ of $U$, the conjunctive normal form expression for the null set
$\emptyset$ is:

\begin{center}
$\left(  S\cup T\right)  \cap\left(  S\cup T^{c}\right)  \cap\left(  S^{c}\cup
T\right)  \cap$ $\left(  S^{c}\cup T^{c}\right)  =\emptyset$
\end{center}

\noindent so for any subset $P$ of $U$, we can distribute the union
$P=P\cup\emptyset$ across the intersections to obtain the \textit{CNF
decomposition of }$P$ using subsets $S,T$ of $U$:

\begin{center}
$P=\left(  S\cup T\cup P\right)  \cap\left(  S\cup T^{c}\cup P\right)
\cap\left(  S^{c}\cup T\cup P\right)  \cap$ $\left(  S^{c}\cup T^{c}\cup
P\right)  $.
\end{center}

\noindent Similarly for arbitrary partitions $\sigma,\tau$ on $U$, the
conjunctive normal form expression for the bottom $\pi$ of the Boolean core
$\mathcal{B}_{\pi}$ is:

\begin{center}
$\left(  \overset{\pi}{\lnot}\overset{\pi}{\lnot}\sigma\vee\overset{\pi
}{\lnot}\overset{\pi}{\lnot}\tau\right)  \wedge\left(  \overset{\pi}{\lnot
}\overset{\pi}{\lnot}\sigma\vee\overset{\pi}{\lnot}\tau\right)  \wedge\left(
\overset{\pi}{\lnot}\sigma\vee\overset{\pi}{\lnot}\overset{\pi}{\lnot}%
\tau\right)  \wedge\left(  \overset{\pi}{\lnot}\sigma\vee\overset{\pi}{\lnot
}\tau\right)  =\pi$
\end{center}

\noindent so for any $\varphi\in\left[  \pi,1\right]  $, we can distribute the
join $\varphi=\varphi\vee\pi$ across the meets to obtain the:

\begin{corollary}
[CNF decomposition identity]For any $\varphi\in\left[  \pi,1\right]  $ and any
partitions $\sigma,\tau$ on $U$:
\end{corollary}

\begin{center}
$\varphi=\left(  \overset{\pi}{\lnot}\overset{\pi}{\lnot}\sigma\vee
\overset{\pi}{\lnot}\overset{\pi}{\lnot}\tau\vee\varphi\right)  \wedge\left(
\overset{\pi}{\lnot}\overset{\pi}{\lnot}\sigma\vee\overset{\pi}{\lnot}\tau
\vee\varphi\right)  \wedge\left(  \overset{\pi}{\lnot}\sigma\vee
\overset{\pi}{\lnot}\overset{\pi}{\lnot}\tau\vee\varphi\right)  \wedge\left(
\overset{\pi}{\lnot}\sigma\vee\overset{\pi}{\lnot}\tau\vee\varphi\right)  $.
\end{center}

Dually, for arbitrary subsets $S,T$ of $U$, the disjunctive normal form
expression for the universe set $U$ is:

\begin{center}
$\left(  S\cap T\right)  \cup\left(  S\cap T^{c}\right)  \cup\left(  S^{c}\cap
T\right)  \cup\left(  S^{c}\cap T^{c}\right)  =U$
\end{center}

\noindent so for any subset $P$ of $U$, we can distribute the intersection
$P=P\cap U$ across the unions to obtain the \textit{DNF decomposition of }$P$
using subsets $S,T$ of $U$:

\begin{center}
$P=\left(  S\cap T\cap P\right)  \cup\left(  S\cap T^{c}\cap P\right)
\cup\left(  S^{c}\cap T\cap P\right)  \cup\left(  S^{c}\cap T^{c}\cap
P\right)  $.
\end{center}

\noindent Similarly for arbitrary partitions $\sigma,\tau$ on $U$, the
disjunctive normal form expression for the top $1$ of the Boolean core
$\mathcal{B}_{\pi}$ is:

\begin{center}
$\left(  \overset{\pi}{\lnot}\overset{\pi}{\lnot}\sigma\wedge\overset{\pi
}{\lnot}\overset{\pi}{\lnot}\tau\right)  \vee\left(  \overset{\pi}{\lnot
}\overset{\pi}{\lnot}\sigma\wedge\overset{\pi}{\lnot}\tau\right)  \vee\left(
\overset{\pi}{\lnot}\sigma\wedge\overset{\pi}{\lnot}\overset{\pi}{\lnot}%
\tau\right)  \vee\left(  \overset{\pi}{\lnot}\sigma\wedge\overset{\pi}{\lnot
}\tau\right)  =1$
\end{center}

\noindent so for any $\varphi\in\left[  \pi,1\right]  $, we can distribute the
meet $\varphi=\varphi\wedge1$ across the joins to obtain the:

\begin{corollary}
[DNF decomposition identity]For any $\varphi\in\left[  \pi,1\right]  $ and any
partitions $\sigma,\tau$ on $U$:
\end{corollary}

\begin{center}
$\varphi=\left(  \overset{\pi}{\lnot}\overset{\pi}{\lnot}\sigma\wedge
\overset{\pi}{\lnot}\overset{\pi}{\lnot}\tau\wedge\varphi\right)  \vee\left(
\overset{\pi}{\lnot}\overset{\pi}{\lnot}\sigma\wedge\overset{\pi}{\lnot}%
\tau\wedge\varphi\right)  \vee\left(  \overset{\pi}{\lnot}\sigma
\wedge\overset{\pi}{\lnot}\overset{\pi}{\lnot}\tau\wedge\varphi\right)
\vee\left(  \overset{\pi}{\lnot}\sigma\wedge\overset{\pi}{\lnot}\tau
\wedge\varphi\right)  $.
\end{center}

Since the CNF and DNF decomposition identities for any $\varphi\in\left[
\pi,1\right]  $ hold for arbitrary partitions $\sigma,\tau$ on $U$, they can
be used as engines to produce other identities by the clever choice of
$\sigma$ and $\tau$ so that simplifications will apply.

For instance, Lawvere,
\citeyear{law:intro}
and
\citeyear{law:coh}%
, has explored two interesting formulas in the context of co-Heyting algebras
(e.g., the closed subsets of a topological space) but both formulas are also
true in the partition algebras $\Pi\left(  U\right)  $. Since Lawvere was
working in a co-Heyting algebra, his suggestive terminology would be more
fitting in the algebra of equivalence relations (represented by the closed
subsets in the non-topological closure space $U\times U$). Lawvere uses the
"difference from $1$" negation that in the algebra of equivalence relations
would be $\left(  \lnot\sigma\right)  ^{d}=\left(  \sigma\Rightarrow0\right)
^{d}=0^{d}-\sigma^{d}$ where $0^{d}=\widehat{1}$ is the top or "one" of that
algebra, and that is dual to the "implication to $0$," i.e., $\lnot
\sigma=\sigma\Rightarrow0$, in the partition algebra. Moreover, we will
relativize the negation using an arbitrary $\pi$ in place of $0$.

Lawvere defines the "boundary" of an element as its meet with its negation, so
dualizing and using $\pi$-negation, we define the $\pi$\textit{-coboundary} of
a partition as the partition obtained from the excluded middle formula using
$\pi$-negation:

\begin{center}
$\partial^{\pi}\sigma=\sigma\vee\overset{\pi}{\lnot}\sigma$

$\pi$-coboundary of a partition $\sigma$
\end{center}

\noindent Lawvere's boundary was "nowhere dense" in the sense that its double
negation was the zero element. In the dual, the $\pi$-coboundary is $\pi
$\textit{-dense} in the sense that its double $\pi$-negation is $1$. The
G\"{o}del $\pi$-transform of the excluded middle tautology is $\overset{\pi
}{\lnot}\overset{\pi}{\lnot}\left[  \left(  \sigma\vee\pi\right)
\vee\overset{\pi}{\lnot}\sigma\right]  $, and since $\left(  \sigma\vee
\pi\right)  \vee\overset{\pi}{\lnot}\sigma=\sigma\vee\overset{\pi}{\lnot
}\sigma$, we have:

\begin{center}
$\overset{\pi}{\lnot}\overset{\pi}{\lnot}\partial^{\pi}\sigma=\overset{\pi
}{\lnot}\overset{\pi}{\lnot}\left(  \sigma\vee\overset{\pi}{\lnot}%
\sigma\right)  =1$.
\end{center}

Lawvere defined the "core" of an element as its double negation but we may
extend this to the double $\pi$-negation $\overset{\pi}{\lnot}\overset{\pi
}{\lnot}\sigma$. Lawvere then shows that each element is equal to its boundary
joined with its core. In the opposite presentation, this result is: $\left(
\sigma\vee\lnot\sigma\right)  \wedge\lnot\lnot\sigma=\sigma$. Generalizing
from $0$ to any $\pi$ then gives the following result in any $\Pi(U)$.

\begin{proposition}
[Lawvere's boundary + core law for partitions]$\partial^{\pi}\sigma
\wedge\overset{\pi}{\lnot}\overset{\pi}{\lnot}\sigma=\sigma\vee\pi$.
\end{proposition}

\noindent Proof: This is easily proved directly from Ore's distributivity
theorem using some basic identities such as: $\overset{\pi}{\lnot}\sigma
\wedge\overset{\pi}{\lnot}\overset{\pi}{\lnot}\sigma=\pi$ and $\sigma
\preceq\overset{\pi}{\lnot}\overset{\pi}{\lnot}\sigma$ so that $\sigma
\vee\overset{\pi}{\lnot}\overset{\pi}{\lnot}\sigma=\overset{\pi}{\lnot
}\overset{\pi}{\lnot}\sigma$. Then using Ore's theorem:

\begin{center}
$\sigma\vee\pi=\sigma\vee\left(  \overset{\pi}{\lnot}\sigma\wedge
\overset{\pi}{\lnot}\overset{\pi}{\lnot}\sigma\right)  =\left(  \sigma
\vee\overset{\pi}{\lnot}\sigma\right)  \wedge\left(  \sigma\vee\overset{\pi
}{\lnot}\overset{\pi}{\lnot}\sigma\right)  =\partial^{\pi}\sigma
\wedge\overset{\pi}{\lnot}\overset{\pi}{\lnot}\sigma$.
\end{center}

\noindent Alternatively, one could take $\varphi=\sigma\vee\pi$ and
$\tau=\sigma$ in the CNF decomposition of $\varphi$ and simplify using the
identities along with $\pi\preceq\overset{\pi}{\lnot}\sigma,\overset{\pi
}{\lnot}\overset{\pi}{\lnot}\sigma$ so that $\pi\vee\overset{\pi}{\lnot}%
\sigma=\overset{\pi}{\lnot}\sigma$ and $\pi\vee\overset{\pi}{\lnot
}\overset{\pi}{\lnot}\sigma=\overset{\pi}{\lnot}\overset{\pi}{\lnot}\sigma$.
$\blacksquare$

Lawvere also shows that the Leibniz rule for taking the derivative of the
product of functions, i.e., $\left(  fg\right)  ^{\prime}=f\left(  g^{\prime
}\right)  +\left(  f^{\prime}\right)  g$, applies in, say, the co-Heyting
algebra of closed subsets of a topological space using the notion of boundary
in place of the derivative. The Leibniz rule holds in the dual algebra of
equivalence relations using the dual notion of $\pi$-boundary, and the dual of
the Leibniz rule holds in the partition algebras $\Pi(U)$ using the notion of
$\pi$-coboundary.

\begin{proposition}
[co-Leibniz rule for partitions]$\partial^{\pi}(\sigma\vee\tau)=\left(
\partial^{\pi}\sigma\vee\tau\right)  \wedge\left(  \sigma\vee\partial^{\pi
}\tau\right)  $.
\end{proposition}

\noindent Proof: This result is an easy consequence of Ore's theorem and the
(weak) DeMorgan law for $\pi$-negation, $\overset{\pi}{\lnot}\left(
\sigma\vee\tau\right)  =\overset{\pi}{\lnot}\sigma\wedge\overset{\pi}{\lnot
}\tau$, that holds in partition logic. Alternatively, the result can be
obtained by taking $\varphi=\partial^{\pi}(\sigma\vee\tau)$ in the CNF
decomposition identity and simplifying. $\blacksquare$

\section{Correctness and completeness for partition logic}

\subsection{Beth-style tableaus for partition logic}

\subsubsection{Classical, intuitionistic, and partition "forcing" models}

It is a familiar fact from classical and intuitionistic logic that logics
might be syntactically presented in a number of ways: Hilbert-style axiom
systems, Beth-style tableaus, natural deduction systems, or Gentzen-style
sequent systems. For partition logic, it seems that the Beth-style tableaus
provide the easiest and most transparent approach so they will be exclusively
used here.

Beth-style tableaus are often called "semantic" since the rules, in effect,
try to construct a model for a formula at the syntactic level. For each of the
connectives, it will be useful to consider the corresponding classical and
intuitionistic tableaus for purposes of comparison. This requires presenting
an appropriate form of the classical and intuitionistic tableaus adapted to
the subset interpretation. As remarked before, classical and intuitionistic
logic are to be interpreted as being about subsets (open subsets in the
intuitionistic case). The rules for the connectives govern when the subsets
contain or don't contain a generic element $u$. Then the partition case is
motivated by elements-distinctions analogy with the generic element $u$
replaced by a generic pair $\left(  u,u^{\prime}\right)  $ of distinct
elements. The conditions governing when subsets contain elements are replaced
by the conditions governing when partitions make distinctions.

Tableaus with signed formulas $T\sigma$ or $F\sigma$ will be used
\cite{smul:fol}. But each signed formula must be accompanied by a generic
element or generic pair as in "$u:T\sigma$" or "$\left(  u,u^{\prime}\right)
:T\sigma$." In the classical or intuitionistic case, $u:T\sigma$ would be
interpreted as saying that the subsets represented by $\sigma$ contains the
element $u$ while $u:F\sigma$ would mean that $\sigma$ (i.e., the subset it
represents) does not contain $u$. Similarly, $\left(  u,u^{\prime}\right)
:T\sigma$ means that the partition represented by $\sigma$ makes the
distinction $\left(  u,u^{\prime}\right)  $, i.e., $u$ and $u^{\prime}$ are in
distinct blocks of $\sigma$, and $\left(  u,u^{\prime}\right)  :F\sigma$ would
mean that $u$ and $u^{\prime}$ are in the same block of $\sigma$.

For classical "propositional" or Boolean logic, the subsets in the intended
interpretation are the subsets of any non-empty universe set $U$. For
intuitionistic "propositional" logic, the intended interpretation is known as
a \textit{Kripke structure} or \textit{intuitionistic forcing model}
\cite{fit:il}. The universe $U$ is endowed with a partial ordering $\leq$ and
the relevant subsets are the up-closed subsets where $S\subseteq U$ is
\textit{up-closed} if $u\in S$ and $u\leq u^{\prime}$ implies that $u^{\prime
}\in S$. These subsets satisfy the conditions for being the open sets of a
topology on $U$. Ordinarily one has a forcing relation ($\models$) between the
points of $U$ and the unsigned formulas. However, signed formulas will be used
here to facilitate the connection to the tableaus using signed formulas:

\begin{center}
$u\models\varphi$ will be written $u:T\varphi$ and $u\not \models \varphi$ is
written $u:F\varphi$.
\end{center}

\noindent A Kripke structure satisfies the structural rule for any $T$-formula
$\varphi$, $\forall u^{\prime}\geq u$, if $u:T\varphi$ then $u^{\prime
}:T\varphi$ so that all $T$\textit{-sets} $T_{\varphi}=\left\{  u|\text{
\ }u:T\varphi\right\}  $ are up-closed (i.e., open). The $T$-conditions for
the connectives are given below while the $F$-conditions are obtained by contraposition.

$u:T\left(  \pi\vee\sigma\right)  $ iff $u:T\pi$ or $u:T\sigma$;

$u:T\left(  \pi\wedge\sigma\right)  $ iff $u:T\pi$ and $u:T\sigma$;

$u:T\left(  \sigma\Rightarrow\pi\right)  $ iff $\forall u^{\prime}\geq u$,
$u^{\prime}:F\sigma$ or $u^{\prime}:T\pi$; and

$u:T\left(  \sigma\mid\pi\right)  $ iff $\forall u^{\prime}\geq u$,
$u^{\prime}:F\sigma$ or $u^{\prime}:F\pi$.

Ordinarily, Kripke structures are defined using negation as a primitive
connective but we can define $\lnot\sigma=\sigma\mid\sigma$ so that setting
$\pi=\sigma$ in the condition for $T\left(  \sigma\mid\pi\right)  $ gives the
derived condition for the negation:

$u:T\left(  \lnot\sigma\right)  $ iff $\forall u^{\prime}\geq u$, $u^{\prime
}:F\sigma$.

Kripke structures make explicit certain features which are left implicit in
classical logic but which must be explicit in partition logic so they are
useful as an expository bridge. It was emphasized from the outset that
classical "propositional" logic should be seen as being about the subsets of a
universe set $U$ and that the "truth table" rules for the connectives are
really the subset membership conditions for a generic element $u$. Since the
classical operations on subsets do not require ever "leaving" the base point
$u$, say, to some other point $u^{\prime}$, all explicit reference to $u$ is
dropped. The definitions can all be interpreted as being about the subsets $0
$ and $1$ of a one point set $\left\{  u\right\}  $ which, in turn, can be
interpreted as falsity and truth for propositions. But once we have the notion
of Kripke structures, then we can see that classical propositional logic
arises, as it were, when the partial ordering on $U$ is discrete which gives
the discrete topology where all subsets are open subsets.

In a Kripke structure, the atomic variables are, in effect, interpreted as $T
$-sets (open subsets) and the conditions for the Kripke structure just give
the membership conditions for the $T$-sets of compound formulas since: $u\in
T_{\varphi}$ iff $u:T\varphi$. Thus a \textit{classical model} for classical
"propositional" logic would be a discrete Kripke structure, i.e., a non-empty
universe set $U$ together with the "forcing" or membership conditions:

$u:T\left(  \pi\vee\sigma\right)  $ iff $u:T\pi$ or $u:T\sigma$;

$u:T\left(  \pi\wedge\sigma\right)  $ iff $u:T\pi$ and $u:T\sigma$;

$u:T\left(  \sigma\Rightarrow\pi\right)  $ iff $u:F\sigma$ or $u:T\pi$; and

$u:T\left(  \sigma\mid\pi\right)  $ iff $u:F\sigma$ or $u:F\pi$.

These conditions for a classical model of propositional logic are just
disguised versions of the usual truth tables but they make explicit the subset
interpretation of the logic. Each formula $\varphi$ would be interpreted in a
model by a subset $T_{\varphi}=\left\{  u|\text{ \ }u:T\varphi\right\}  $, and
the rules could be restated as membership conditions for a generic element,
e.g., $u\in T_{\pi\vee\sigma}$ iff $u\in T_{\pi}$ or $u\in T_{\sigma}$, and so forth.

In the usual treatment of Kripke structures, a formula $\varphi$ is
\textit{intuitionistically valid} if it is forced at every point in any Kripke
structure. But this is equivalent to saying that for any interpretation of the
atomic variables of $\varphi$ as open subsets of the model, the whole formula
evaluates to the universe set $U$. In the discrete or classical case, it means
that a formula is a classical or subset tautology if regardless of the subsets
of $U$ assigned to the atomic variables of the formula, the formula evaluates
to the universe set $U$ (for any non-empty $U$).

We now have sufficient motivation to define the analogous \textit{partition
forcing models}. We start with a universe set $U$ with two or more elements.
The points in the classical and Kripke structures are replaced by the pairs
$\left(  u,u^{\prime}\right)  $ of distinct points from $U$. Instead of using
an explicit forcing relation between pairs and formulas, we will again use
signed formulas so that:

\begin{center}
$\left(  u,u^{\prime}\right)  \models\varphi$ is written as $\left(
u,u^{\prime}\right)  :T\varphi$, and $\left(  u,u^{\prime}\right)
\not \models \varphi$ is written as $\left(  u,u^{\prime}\right)  :F\varphi$.
\end{center}

\noindent Unlike the points $u$ or $u^{\prime}$, the pairs $\left(
u,u^{\prime}\right)  $ have an internal structure; a pair $\left(
u,u^{\prime}\right)  $ can be reversed to $\left(  u^{\prime},u\right)  $ and
pairs can be connected in triangles as in $\left(  u,u^{\prime}\right)  $,
$\left(  u,a\right)  $, and $\left(  a,u^{\prime}\right)  $ or in longer
chains. Hence a partition forcing model has two structural conditions
reflecting the symmetry and anti-transitivity of partition relations:

if $\left(  u,u^{\prime}\right)  :T\varphi$, then $\left(  u^{\prime
},u\right)  :T\varphi$;

if $\left(  u,u^{\prime}\right)  :T\varphi$, then for any other $a$, $\left(
u,a\right)  :T\varphi$ or $\left(  a,u^{\prime}\right)  :T\varphi$.

\noindent No rule is needed to enforce the anti-reflexivity of partition
relations since the notation always assumes that $\left(  u,u^{\prime}\right)
$ is a pair of distinct elements.

The unstructured universe set $U$ still determines the complete undirected
graph $K\left(  U\right)  $ on $U$ which has a link $\left(  u,u^{\prime
}\right)  $ between any two distinct points. A\textit{\ }$u,u^{\prime}%
$\textit{-chain} is a finite sequence of links, $\left(  u_{1},u_{2}\right)
$,$\left(  u_{2},u_{3}\right)  $,...,$\left(  u_{n-1},u_{n}\right)  $, with
$u=u_{1} $ and $u^{\prime}=u_{n}$ as the endpoints. In particular, the
\textit{base pair} $\left(  u,u^{\prime}\right)  $ is a one-link $u,u^{\prime}
$-chain, and any third element $a$ gives the two-link $u,u^{\prime}$-chain
$\left(  u,a\right)  $ and $\left(  a,u^{\prime}\right)  $. Recall that the
\textit{Boolean condition} for any signed compound formula $\pi\ast\sigma$ is
the disjunction or conjunction of the pair of signed formulas that hold in a
classical model for the constituents $\pi$ and $\sigma$ where $\ast$ is any
binary operation.

Now the "forcing conditions" for a partition forcing model can be stated for
the $T$-signed formulas (with the $F$-rules obtained by contraposition).

$\left(  u,u^{\prime}\right)  :T\left(  \pi\vee\sigma\right)  $ iff $\left(
u,u^{\prime}\right)  :T\pi$ or $\left(  u,u^{\prime}\right)  :T\sigma$ (i.e.,
the Boolean condition holds at the base pair);

$\left(  u,u^{\prime}\right)  :T\left(  \sigma\Rightarrow\pi\right)  $ iff for
any $1$- or $2$-link $u,u^{\prime}$-chain, the Boolean condition (i.e.,
$F\sigma$ or $T\pi$) holds on some chain link;

$\left(  u,u^{\prime}\right)  :T\left(  \pi\wedge\sigma\right)  $ iff for any
$u,u^{\prime}$-chain, the Boolean condition (i.e., $T\pi$ and $T\sigma$) holds
on some chain link; and

$\left(  u,u^{\prime}\right)  :T\left(  \sigma\mid\pi\right)  $ iff for any
$u,u^{\prime}$-chain, the Boolean condition (i.e., $F\sigma$ or $F\pi$) holds
on some chain link.

The $T$-sets are $T_{\varphi}=\left\{  \left(  u,u^{\prime}\right)  |\text{
\ }\left(  u,u^{\prime}\right)  :T\varphi\right\}  $, and a partition validity
would be a formula whose $T$-set consisted of all pairs $\left(  u,u^{\prime
}\right)  $ of distinct elements in all partition forcing models.

These partition forcing models have been defined so that one can see the
analogies between Kripke structures (and classical structures as the discrete
special case). But we have met the partition forcing models before; they are
just a different presentation of the dit-set representation of the partition
algebras $\Pi(U)$:

\begin{center}
Partition forcing model = dit-set representation of $\Pi(U)$.
\end{center}

\noindent The $T$-sets are the dit sets since $\left(  u,u^{\prime}\right)
:T\varphi$ is the same as $\left(  u,u^{\prime}\right)  \in\operatorname*{dit}%
\left(  \varphi\right)  $ so that $T_{\varphi}=\operatorname*{dit}\left(
\varphi\right)  $.

The presentation of the dit-set representation as a "partition forcing model"
nevertheless brings out a number of analogies between the distinguishing-cut
and falsifying-chain results in partition and related results in
intuitionistic and classical logic. In a Kripke structure, the order
structural condition is that if $T\varphi$ holds at a point $u$, then it holds
at any higher point $u^{\prime}\geq u$. In a partition forcing model, the
anti-transitivity structure condition is that if $T\varphi$ holds at any pair
$\left(  u,u^{\prime}\right)  $, then it holds at some link on any
$u,u^{\prime}$-chain. Moreover, the conditions for the connectives provide a
stronger version of the analogy. Let $\ast$ be any operation such as $\vee$,
$\Rightarrow$, $\wedge$, or $\mid$.

\begin{center}%
\begin{tabular}
[c]{|c|}\hline
Partition forcing model\\\hline\hline
$\sigma\ast\pi$ distinguishes $\left(  u,u^{\prime}\right)  $, i.e., $\left(
u,u^{\prime}\right)  :T\left(  \sigma\ast\pi\right)  $\\
iff $\forall$ $u,u^{\prime}$-chains, the Boolean conditions for $T\left(
\sigma\ast\pi\right)  $\\
hold at some link on the chain.\\\hline
$\sigma\ast\pi$ identifies $\left(  u,u^{\prime}\right)  $, i.e., $\left(
u,u^{\prime}\right)  :F\left(  \sigma\ast\pi\right)  $\\
iff $\exists$ $u,u^{\prime}$-chain, with the Boolean conditions for $F\left(
\sigma\ast\pi\right)  $\\
holding at every link.\\\hline
\end{tabular}

--%

\begin{tabular}
[c]{|c|}\hline
Intuitionistic forcing model\\\hline\hline
$\sigma\ast\pi$ contains $u$, i.e., $u:T\left(  \sigma\ast\pi\right)  $\\
iff $\forall$ $u^{\prime}\geq u$, the Boolean conditions for $T\left(
\sigma\ast\pi\right)  $ hold at $u^{\prime}$.\\\hline
$\sigma\ast\pi$ does not contain $u$, i.e., $u:F\left(  \sigma\ast\pi\right)
$\\
iff $\exists$ $u^{\prime}\geq u$, such that the Boolean conditions for
$F\left(  \sigma\ast\pi\right)  $\\
hold at $u^{\prime}$.\footnote{The chain-cut results in partition logic have
an even closer analogy in the intuitionistic case if one uses the Beth
semantics of paths and bars in partially ordered sets, see \cite[p.
237]{van-d:il} or \cite[p. 276]{restall:sub}.}\\\hline
\end{tabular}

--%

\begin{tabular}
[c]{|c|}\hline
Classical forcing model\\\hline\hline
$\sigma\ast\pi$ contains $u$, i.e., $u:T\left(  \sigma\ast\pi\right)  $\\
iff the Boolean conditions for $T\left(  \sigma\ast\pi\right)  $ hold at
$u$.\\\hline
$\sigma\ast\pi$ does not contain $u$, i.e., $u:F\left(  \sigma\ast\pi\right)
$\\
iff the Boolean conditions for $F\left(  \sigma\ast\pi\right)  $ hold at
$u$.\\\hline
\end{tabular}

\end{center}

Some pains have been taken to emphasize the analogies between the Kripke
structure model and the classical and partition "forcing" models. But the
classical and partition models are just a fancy way to describe, respectively,
the membership conditions for subsets of a set $U$ and the distinction
conditions for partitions on a set $U$. Moreover, ordinary subset logic and
partition logic are at the same mathematical level in the sense that both
start with an unstructured set $U.$ The subsets of a set and the partitions on
a set can both be described without assuming any additional structure. In the
intuitionistic case, either a topology or a partial order (which induces the
topology of up-closed subsets as the open subsets) is assumed on the universe
set $U$.\footnote{Starting with Kripke structures as models for intuitionistic
and modal logics, there has recently been a vast proliferation of logics
modeled by sets with orderings or closure operations along with a variety of
compatibility and accessibility relations on the sets; see \cite{restall:sub}
for a survey. In contrast to this profusion of logics, partition logic, like
classical subset logic, is modeled using only unstructured sets $U$.}

\subsubsection{Tableau structural rules}

In general, the intuitionistic and partition $F$-rules will have a similar
form. For any connective $\ast$, the intuitionistic rule is that $u:F\left(
\pi\ast\sigma\right)  $ iff $\exists u^{\prime}\geq u$ such that the Boolean
condition for $F\left(  \pi\ast\sigma\right)  $ holds at $u^{\prime}$, while
the partition rule is that $\left(  u,u^{\prime}\right)  :F\left(  \pi
\ast\sigma\right)  $ iff $\exists u,u^{\prime}$-chain such that the Boolean
condition for $F\left(  \pi\ast\sigma\right)  $ holds at every link on the chain.

By the same token, we could formulate the intuitionistic and partition
$T$-rules as contrapositives. For the intuitionistic $T\ast$ rule, $u:T\left(
\pi\ast\sigma\right)  $ holds iff $\forall u^{\prime}\geq u$, the Boolean
condition for $T\left(  \pi\ast\sigma\right)  $ holds at $u^{\prime}$, and the
partition rule is that $\left(  u,u^{\prime}\right)  :T\left(  \pi\ast
\sigma\right)  $ iff $\forall u,u^{\prime}$-chains, there is a link on the
chain where the Boolean condition for $T\left(  \pi\ast\sigma\right)  $ holds
on that link.

But the $T$-rules are written in a simplified way where the Boolean condition
for $T\left(  \pi\ast\sigma\right)  $ holds at the base, and then is
transmitted to a new base with that Boolean condition also holding there. For
instance, the intuitionistic $T$-rule for $\pi\ast\sigma$ will be given in the
simplified form as $u:T\left(  \pi\ast\sigma\right)  $ implies the Boolean
condition for $T\left(  \pi\ast\sigma\right)  $ holds at $u$, together with a
$T$-transmitting rule so any $T$'s are transmitted to higher points in the
ordering. Similarly in the partition case, we have used the simplified rule
where $\left(  u,u^{\prime}\right)  :T\left(  \pi\ast\sigma\right)  $ implies
the Boolean condition for $T\left(  \pi\ast\sigma\right)  $ also holds at
$\left(  u,u^{\prime}\right)  $ and then the following $T$-anti-transitivity
rule transmits any $T$'s to some link in any $u,u^{\prime}$-chain.

The two $T$-transmitting structural rules for the intuitionistic and partition
cases are as follows. The $T$-anti-transitivity rule splits into two
alternatives given by the vertical line
$\vert$%
. Context should suffice to avoid confusion between the vertical line
$\vert$
separating branches in the tableau tree and the Sheffer stroke $\mid$ of the
nand operation.

\begin{center}%
\begin{tabular}
[c]{|c|c|}\hline
\multicolumn{1}{|c|}{%
\begin{tabular}
[c]{c}%
$u:T\varphi$\\\hline
$\forall a>u$, $a:T\varphi$%
\end{tabular}
} &
\begin{tabular}
[c]{c}%
$\left(  u,u^{\prime}\right)  :T\varphi$\\\hline
$\forall a$, $\left(  u,a\right)  :T\varphi$
$\vert$
$\left(  a,u^{\prime}\right)  :T\varphi$%
\end{tabular}
\\
\multicolumn{1}{|c|}{Intuitionistic $T$-transmitting rule} & Partition
$T$-anti-transitivity rule\\\hline
\end{tabular}

\end{center}

\noindent An easy corollary implies that a $T\varphi$ holding at $\left(
u,u^{\prime}\right)  $ is transmitted to some link in any $u,u^{\prime}$-chain.

The $T$-transmitting rules can also be contraposited to derive "$F$%
-reflecting" rules.

\begin{center}%
\begin{tabular}
[c]{|c|c|}\hline
\multicolumn{1}{|c|}{%
\begin{tabular}
[c]{c}%
$\exists a>u$, $a:F\varphi$\\\hline
$u:F\varphi$%
\end{tabular}
} &
\begin{tabular}
[c]{c}%
$\exists a$, $(u,a):F\varphi$ and $(a,u^{\prime}):F\varphi$\\\hline
$\left(  u,u^{\prime}\right)  :F\varphi$%
\end{tabular}
\\
\multicolumn{1}{|c|}{Intuitionistic $F$-reflecting rule} & Partition
$F$-transitivity rule\\\hline
\end{tabular}

\end{center}

\noindent Thus if $F\varphi$ holds at each link on any $u,u^{\prime}$-chain,
then $\left(  u,u^{\prime}\right)  :F\varphi$ follows.

Partition relations and their complementary equivalence relations are
symmetric. Since we are using the ordered pairs $\left(  u,u^{\prime}\right)
$ rather than the unordered pairs $\left\{  u,u^{\prime}\right\}  $, we need
rules to enforce that symmetry for the ordered pairs.

\begin{center}%
\begin{tabular}
[c]{|c|c|}\hline
\multicolumn{1}{|c|}{%
\begin{tabular}
[c]{c}%
$(u,u^{\prime}):T\varphi$\\\hline
$\left(  u^{\prime},u\right)  :T\varphi$%
\end{tabular}
} &
\begin{tabular}
[c]{c}%
$(u,u^{\prime}):F\varphi$\\\hline
$\left(  u^{\prime},u\right)  :F\varphi$%
\end{tabular}
\\
\multicolumn{1}{|c|}{Partition $T$ symmetric rule} & Partition $F$ symmetric
rule\\\hline
\end{tabular}

\end{center}

Equivalence relations are reflexive and partition relations are anti-reflexive
but we don't need rules to enforce that since we have stipulated that the
ordered pairs $\left(  u,u^{\prime}\right)  $ in the rules are always of
distinct elements.

\subsubsection{Tableaus for the partition join}

The tableau rules are given, for comparison purposes, for the three logics:
subset, intuitionistic, and partition. The terms $u,u^{\prime},a,b,...,c$ are
now elements in the syntactic machinery of the tableau rules with the intended
interpretations that have been already given; $u:T\varphi$ would be
interpreted as $u$ is a member of the set that interprets $\varphi$ in the
classical and intuitionistic rules while $\left(  u,u^{\prime}\right)
:T\varphi$ would be interpreted as $\left(  u,u^{\prime}\right)  $ is a
distinction of the partition that interprets $\varphi$, and similarly for the
$F$-formulas. The four operations of $\vee$, $\wedge$, $\Rightarrow$, and
$\mid$ will be taken as primitive in all the logics with the constant $1$
defined as $\sigma\Rightarrow\sigma$ for any $\sigma$ and $0$ defined as
$1\mid1$.

In general, the syntactic eliminative $T$ rules give the left-to-right
implication in the "forcing" models described above, and the $F$ rules are
obtained by contrapositing the implication in the other direction. To compare
the tableaus for these three logics, we start with the join where the
eliminative tableaus are the most alike.

\begin{center}%
\begin{tabular}
[c]{|c|c|c|}\hline%
\begin{tabular}
[c]{c}%
$u:F\left(  \pi\vee\sigma\right)  $\\\hline
$u:F\pi,F\sigma$%
\end{tabular}
&
\begin{tabular}
[c]{c}%
$u:F\left(  \pi\vee\sigma\right)  $\\\hline
$u:F\pi,F\sigma$%
\end{tabular}
&
\begin{tabular}
[c]{c}%
$\left(  u,u^{\prime}\right)  :F\left(  \pi\vee\sigma\right)  $\\\hline
$\left(  u,u^{\prime}\right)  :F\pi,F\sigma$%
\end{tabular}
\\
Classical $F\vee$ rule & Intuitionistic $F\vee$ rule & Partition $F\vee$
rule\\\hline
\end{tabular}

\end{center}

The $T\vee$ rules use the notion of a splitting of alternatives which is
indicated by a vertical line.

\begin{center}%
\begin{tabular}
[c]{|c|c|c|}\hline%
\begin{tabular}
[c]{c}%
$u:T\left(  \pi\vee\sigma\right)  $\\\hline
$u:T\pi$
$\vert$
$u:T\sigma$%
\end{tabular}
&
\begin{tabular}
[c]{c}%
$u:T\left(  \pi\vee\sigma\right)  $\\\hline
$u:T\pi$
$\vert$
$u:T\sigma$%
\end{tabular}
&
\begin{tabular}
[c]{c}%
$\left(  u,u^{\prime}\right)  :T\left(  \pi\vee\sigma\right)  $\\\hline
$\left(  u,u^{\prime}\right)  :T\pi$
$\vert$
$\left(  u,u^{\prime}\right)  :T\sigma$%
\end{tabular}
\\
Classical $T\vee$ rule & Intuitionistic $T\vee$ rule & Partition $T\vee$
rule\\\hline
\end{tabular}

\end{center}

The close analogies between the classical and intuitionistic rules on the one
hand and the partition rules on the other hand are all by virtue of turning
the lattice of partitions right side up.

\subsubsection{Tableaus for the partition implication}

The complications arise in the $F$ rules so we begin with the $T$ rules.

\begin{center}%
\begin{tabular}
[c]{|c|c|c|}\hline%
\begin{tabular}
[c]{c}%
$u:T\left(  \sigma\Rightarrow\pi\right)  $\\\hline
$u:F\sigma$
$\vert$
$u:T\pi$%
\end{tabular}
&
\begin{tabular}
[c]{c}%
$u:T\left(  \sigma\Rightarrow\pi\right)  $\\\hline
$u:F\sigma$
$\vert$
$u:T\pi$%
\end{tabular}
&
\begin{tabular}
[c]{c}%
$\left(  u,u^{\prime}\right)  :T\left(  \sigma\Rightarrow\pi\right)  $\\\hline
$\left(  u,u^{\prime}\right)  :F\sigma$
$\vert$
$\left(  u,u^{\prime}\right)  :T\pi$%
\end{tabular}
\\
Classical $T\Rightarrow$ rule & Intuitionistic $T\Rightarrow$ rule & Partition
$T\Rightarrow$ rule\\\hline
\end{tabular}

\end{center}

The classical rules never leave the base point $u$ so $u$ is usually left
implicit. In the intuitionistic $F\Rightarrow$ rule, a new element $a$ may be
introduced. Since the Beth-style tableau rules, in effect, try to construct a
model of a formula using syntactic machinery, the ordering between the points
in a Kripke structure must already be introduced as an ordering between
elements. In particular, in the intuitionistic $F\Rightarrow$ rule, the new
element $a$ introduced in the consequence of the rule is higher in the
ordering of elements than the base point used in the premise of the rule. In
other treatments of the intuitionistic rule $F\Rightarrow$ as in Fitting
\citeyear{fit:il}%
, the elements such as $u$ and $a$ are also left implicit but the rules that
require leaving the base point to move higher in the ordering (i.e., the
$F\Rightarrow$ and $F\lnot$ rules) are indicated by dropping any other
$F$-formulas in the premise and keeping only the $T$-formulas since only the
$T$-formulas are transmitted to points higher in the ordering. We will not
fully develop our version of the intuitionistic tableaus but we are presenting
them to bring out the analogies with the partition tableaus.

In the partition $F\Rightarrow$ we may introduce a new element $a$ but there
is no ordering on the elements. There is always the notion of a chain of pairs
of elements, and the partition $F\Rightarrow$ rule says that the Boolean
condition for $F\left(  \sigma\Rightarrow\pi\right)  $ holds on each link of
the chain $\left(  u,a\right)  ,\left(  a,u^{\prime}\right)  $.

\begin{center}%
\begin{tabular}
[c]{|c|c|}\hline%
\begin{tabular}
[c]{c}%
$u:F\left(  \sigma\Rightarrow\pi\right)  $\\\hline
$u:T\sigma,F\pi$%
\end{tabular}
&
\begin{tabular}
[c]{c}%
$u:F\left(  \sigma\Rightarrow\pi\right)  $\\\hline
$\exists a\geq u$, $a:T\sigma,F\pi$%
\end{tabular}
\\
Classical $F\Rightarrow$ rule & Intuitionistic $F\Rightarrow$ rule\\\hline
\end{tabular}

\begin{tabular}
[c]{|c|}\hline%
\begin{tabular}
[c]{c}%
$\left(  u,u^{\prime}\right)  :F\left(  \sigma\Rightarrow\pi\right)  $\\\hline
$\exists u,u^{\prime}$-chain (1 or 2 links) with$\ $ $T\sigma,F\pi$ on each
link
\end{tabular}
\\
Partition $F\Rightarrow$ rule\\\hline
\end{tabular}

\end{center}

Since this is the first partition tableau rule that might introduce a new
element, we have to be more explicit about how the tableau rules will be used
here. We are given some partition formula $\varphi$ and we begin a tableau for
$\varphi$ with the statement $\left(  u_{0},u_{1}\right)  :F\varphi$. Since a
tableau can branch like an upside-down tree, this initial statement $\left(
u_{0},u_{1}\right)  :F\varphi$ is the root of the tree. The application of the
tableau rules attempts to construct a partition on some model set $U$
containing $u_{0}$ and $u_{1}$ where $\left(  u_{0},u_{1}\right)  :F\varphi$
holds, i.e., to construct a countermodel for $\varphi$. The universe set
starts at $U_{0}=\left\{  u_{0},u_{1}\right\}  $, and each application of a
rule introducing one or more new elements will take the developing model from
some $U_{n}$ to $U_{n+1}$ which is $U_{n}$ plus the new elements.\footnote{It
may be useful to keep in mind the analogies with the development of models in
classical first-order logic using tableaus \cite{smul:fol}. We are from the
outset seeing the new "constants" being introduced as elements in a potential
model (in a manner reminiscent of the L\"{o}wenheim-Skolem theorem in
classical first-order logic).} Each $U_{n}$ might be called a \textit{stage}
of the developing model.

New elements should be introduced only as a last resort. Since new element
might be introduced only by $F$-rules, before introducing new elements to make
a falsifying chain, we need to first check that a falsifying chain could not
be formed using the existing elements. Given the premise $\left(  u,u^{\prime
}\right)  :F\left(  \sigma\Rightarrow\pi\right)  $, the "base pair" $\left(
u,u^{\prime}\right)  $ might be a one-link falsifying chain if $\left(
u,u^{\prime}\right)  :T\sigma,F\pi$, or the falsifying chain might be
constructed using an element $a\in U_{n}$ of the evolving universe set at that
stage. For instance, we might have already derived $\left(  u,u^{\prime
}\right)  :F\sigma,F\pi$ so the base pair was not a one-link falsifying chain,
but there might be an element $a\in U_{n}$ in the evolving universe set at
that stage and on that branch of the tableau where, say, $\left(  u,a\right)
:F\pi$ held. Then the $F$-transitivity rule given below would yield $\left(
a,u^{\prime}\right)  :F\pi$. If we then initiated a new branch with the
assumption $\left(  u,a\right)  :T\sigma$ then the $T$-anti-transitivity rule
given below would imply $\left(  a,u^{\prime}\right)  :T\sigma$ and we would
have a falsifying chain for $\left(  u,u^{\prime}\right)  :F\left(
\sigma\Rightarrow\pi\right)  $ without introducing any new elements. Such a
falsifying chain using existing elements might be called a \textit{back-chain}%
. Thus there are a finite number of options to establish a falsifying
back-chain before taking the "last option" of introducing a new element and
thus a new stage in the developing countermodel. Each of these options may
create a branch in the tree. Since each $U_{n}$ is a finite set, there are
only a finite number of possible back-chains (including the one-link
back-chain of the base pair) so only a finite number of branches might be
created by applying the rule.

We know from the previous falsifying-chain result that when the atomic
variables of some formula $\sigma\ast\pi$ are interpreted as partitions on
some universe set $U$, then $F\left(  \sigma\ast\pi\right)  $ will hold at
some pair $\left(  u,u^{\prime}\right)  $ iff there is a falsifying chain with
the Boolean conditions for $F\left(  \sigma\ast\pi\right)  $ holding at each
link on the chain. The $F\Rightarrow$ rule gives us the syntactic version of
that semantic theorem in the following sense. Each branch resulting from
applying the rule to $\left(  u,u^{\prime}\right)  :F\left(  \sigma
\Rightarrow\pi\right)  $ will have the statements for a falsifying chain
either at the base pair $\left(  u,u^{\prime}\right)  :T\sigma,F\pi$ or on a
two-link falsifying chain $\left(  u,a\right)  ,\left(  a,u^{\prime}\right)
:T\sigma,F\pi$ which might be a back-chain if $a\in U_{n}$ or a new chain if
$a$ is a new element that yields $U_{n+1}=U_{n}\cup\left\{  a\right\}  $. For
the operations of the meet $\wedge$ and nand $\mid$, the falsifying chains
could have more than two links and thus involve two or more elements other
than the base pair $\left(  u,u^{\prime}\right)  $. In that case, the
possibility arises of \textit{mixed chains} using some existing elements in
$U_{n}$ and some new elements.

We also know from the distinguishing cut result that when the atomic variables
of some formula $\sigma\ast\pi$ are interpreted as partitions on some universe
set $U$, then $T\left(  \sigma\ast\pi\right)  $ will hold at some pair
$\left(  u,u^{\prime}\right)  $ iff for every $u,u^{\prime}$-chain, there is a
link $\left(  a,b\right)  $ where the Boolean conditions for $T\left(
\sigma\ast\pi\right)  $ hold. The $T$-rules together with the $T$%
-anti-transitivity rule ensure that the corresponding formulas are derived in
the developing branch of a tableau. For instance, in the present case of the
implication $\left(  u,u^{\prime}\right)  :T\left(  \sigma\Rightarrow
\tau\right)  $, the $T$-anti-transitivity rule implies that for any
$u,u^{\prime}$-chain using the elements of $U_{n}$, there is a link $\left(
a,b\right)  $ where $T\left(  \sigma\Rightarrow\pi\right)  $ holds and then
the $T\Rightarrow$ rule implies that either $\left(  a,b\right)  :F\sigma$ or
$\left(  a,b\right)  :T\pi$ holds--which are the Boolean conditions for
$T\left(  \sigma\Rightarrow\pi\right)  $ holding at $\left(  a,b\right)  $.

Similar remarks apply to all the $T\ast$ and $F\ast$ rules where $\ast$ is
$\vee$, $\Rightarrow$, $\wedge$, or $\mid$.

\subsubsection{Tableaus for the partition meet}

All three of the $T\wedge$ rules are rather standard.

\begin{center}%
\begin{tabular}
[c]{|c|c|c|}\hline%
\begin{tabular}
[c]{c}%
$u:T\left(  \pi\wedge\sigma\right)  $\\\hline
$u:T\pi,T\sigma$%
\end{tabular}
&
\begin{tabular}
[c]{c}%
$u:T\left(  \pi\wedge\sigma\right)  $\\\hline
$u:T\pi,T\sigma$%
\end{tabular}
&
\begin{tabular}
[c]{c}%
$\left(  u,u^{\prime}\right)  :T\left(  \pi\wedge\sigma\right)  $\\\hline
$\left(  u,u^{\prime}\right)  :T\pi,T\sigma$%
\end{tabular}
\\
Classical $T\wedge$ rule & Intuitionistic $T\wedge$ rule & Partition $T\wedge$
rule\\\hline
\end{tabular}

\end{center}

The classical and intuitionistic rules for $F\wedge$ are standard while the
partition $F\wedge$ is complicated since it involves a chain of elements with
the Boolean condition, $F\pi$ or $F\sigma$, holding on each link. In the
eliminative rule for the universal quantifier in classical first-order logic,
we go from a premise $u:\left(  \forall x\right)  \varphi\left(  x\right)  $
to a conclusion of either $u:\varphi\left(  a\right)  $ where $a$ is a
constant in the developing model or $u:\varphi\left(  x^{\prime}\right)  $
where $x^{\prime}$ is a variable that can latter be replaced by a constant. In
the partition $F\wedge$ rule, we have a similar situation when there is no
falsifying back-chain so we need to introduce new elements to be strung
together to make a falsifying chain. How many new elements should be
introduced? In each branch of a tableau, we may eventually arrive at a
contradiction in the form $\left(  a,b\right)  :T\sigma,F\sigma$ at some pair
in which case the branch would \textit{close}. Along that branch, no
countermodel can be constructed so the branch is terminated. But a branch
might be "falsely" terminated if we don't introduce enough new links in the
falsifying chain of the $F\wedge$ rule. For instance, suppose we also had
$\left(  u,u^{\prime}\right)  :T\phi_{1},T\phi_{2},T\phi_{3}$ in the branch
and any two of these formulas holding at the same pair would give rise to a
contradiction. Then if we had only introduced one new element to give the
two-link falsifying chain $\left(  u,a\right)  :F\sigma$ and $\left(
a,u^{\prime}\right)  :F\pi$, then the $T$-anti-transitivity rule would have to
"transmit" two of the three formulas $T\phi_{1},T\phi_{2},T\phi_{3}$ to one of
the links in the chain and we would seem to have a closure of the branch. But
we could just as well have introduced two new elements $a$ and $b $ so we had
a falsifying chain of the three links $\left(  u,a\right)  ,\left(
a,b\right)  ,\left(  b,u^{\prime}\right)  $ and then each of the three
formulas could be transmitted to a different link avoiding the contradiction.
A crude upper bound on the number of necessary links is the number of
subformulas of the formula $\varphi$ in the root of the tree.

Hence when a branch closes, we must be sure that it would still close
regardless of the length of the falsifying chain introduced in the $F\wedge$
rule. This can be done by ensuring that any falsifying chain from the
$F\wedge$ rule in a closed branch could have been treated as a "variable" or
generic chain so that whenever some $T\phi_{i}$ holding at $\left(
u,u^{\prime}\right)  $ is transmitted to the chain, then it must have its
"own" link and must not be forced to unnecessarily share a link with some
other $T\phi_{j}$. If a branch does not close, then we need to construct a
countermodel from the elements introduced in that branch (see the Satisfaction
Theorem below) so we need to have introduced specific elements in an open branch.

In the $F\wedge$ rule, the elements $a,b,...,c$ form a $u,u^{\prime}$-chain,
$\left(  u,a\right)  ,\left(  a,b\right)  ,...,\left(  c,u^{\prime}\right)  $.

\begin{center}%
\begin{tabular}
[c]{|c|c|}\hline%
\begin{tabular}
[c]{c}%
$u:F\left(  \sigma\wedge\pi\right)  $\\\hline
$u:F\sigma$
$\vert$
$u:F\pi$%
\end{tabular}
&
\begin{tabular}
[c]{c}%
$u:F\left(  \sigma\wedge\pi\right)  $\\\hline
$u:F\sigma$
$\vert$
$u:F\pi$%
\end{tabular}
\\
Classical $F\wedge$ rule & Intuitionistic $F\wedge$ rule\\\hline
\end{tabular}

\begin{tabular}
[c]{|c|}\hline%
\begin{tabular}
[c]{c}%
$\left(  u,u^{\prime}\right)  :F\left(  \sigma\wedge\pi\right)  $\\\hline
$\exists a,b,...c$ so the $u,u^{\prime}$-chain has $F\sigma$ or $F\pi$ on each
link
\end{tabular}
\\
Partition $F\wedge$ rule\\\hline
\end{tabular}

\end{center}

\noindent By the $F$-transitivity rule, two consecutive $F\sigma$ links could
be shorted to one $F\sigma$ link so we may assume that the links of the
falsifying chain are alternating. As in the case of the $F\Rightarrow$ rule,
the falsifying chain might be a back-chain established using the elements of
the current stage $U_{n}$ without introducing new elements or a mixed chain
with some old and some new elements. When new assumptions are made to have a
falsifying back-chain, that creates a new branch. When new elements are
introduced and $T$-formulas are transmitted to the links, then each way this
could be done is a new branch. The possibilities quickly multiply but they are
always finite at each stage.

\subsubsection{Tableaus for the partition nand}

All three of the $T\mid$ rules are rather standard.

\begin{center}%
\begin{tabular}
[c]{|c|c|c|}\hline%
\begin{tabular}
[c]{c}%
$u:T\left(  \pi\mid\sigma\right)  $\\\hline
$u:F\pi$
$\vert$
$u:F\sigma$%
\end{tabular}
&
\begin{tabular}
[c]{c}%
$u:T\left(  \pi\mid\sigma\right)  $\\\hline
$u:F\pi$
$\vert$
$u:F\sigma$%
\end{tabular}
&
\begin{tabular}
[c]{c}%
$\left(  u,u^{\prime}\right)  :T\left(  \pi\mid\sigma\right)  $\\\hline
$\left(  u,u^{\prime}\right)  :F\pi$
$\vert$
$\left(  u,u^{\prime}\right)  :F\sigma$%
\end{tabular}
\\
Classical $T\mid$ rule & Intuitionistic $T\mid$ rule & Partition $T\mid$
rule\\\hline
\end{tabular}

\end{center}

In the "intuitionistic" $F\mid$ rule (which we have invented since the nand
operation is not ordinarily used in intuitionistic logic), a new element $a$
is introduced so that $a\geq u$ in the partial ordering of elements so that
the Boolean condition for $F\left(  \pi\mid\sigma\right)  $, i.e.,
$T\pi,T\sigma$, holds at that point. In the $F\mid$ rule for partitions we
already know that four links suffice in an falsifying chain so we only need to
introduce at most three new elements $a,b,c$ to form the falsifying
$u,u^{\prime}$-chain where the same Boolean conditions hold at each link.

\begin{center}%
\begin{tabular}
[c]{|c|c|}\hline%
\begin{tabular}
[c]{c}%
$u:F\left(  \pi\mid\sigma\right)  $\\\hline
$u:T\pi,T\sigma$%
\end{tabular}
&
\begin{tabular}
[c]{c}%
$u:F\left(  \pi\mid\sigma\right)  $\\\hline
$\exists a\geq u$, $a:T\pi,T\sigma$%
\end{tabular}
\\
Classical $F\mid$ rule & Intuitionistic $F\mid$ rule\\\hline
\end{tabular}

\begin{tabular}
[c]{|c|}\hline%
\begin{tabular}
[c]{c}%
$\left(  u,u^{\prime}\right)  :F\left(  \pi\mid\sigma\right)  $\\\hline
$\exists$ $u,u^{\prime}$-chain (at most four links) with $T\pi,T\sigma$ on
each link
\end{tabular}
\\
Partition $F\mid$ rule\\\hline
\end{tabular}

\end{center}

\noindent As before, the falsifying chain could be a back-chain. For the
option where new elements are introduced, at most three elements need to be
introduced since four links suffice in any falsifying chain for the nand
$\pi\mid\sigma$.

\subsubsection{Examples of proofs and countermodels using the $F\wedge$ rule}

Starting with the assumption that a "root" formula $\varphi$ does not
distinguish a generic pair $\left(  u_{0},u_{1}\right)  $, i.e., $\left(
u_{0},u_{1}\right)  :F\varphi$, the tableau rules for the connectives (as
opposed to the structural rules) eliminate the main connective of a formula at
each step. If all branches terminate with a contradiction such as
$T\sigma,F\sigma$ at some pair, then the tableau constitutes a proof of the
formula $\varphi$, i.e., $\varphi$ is a \textit{theorem} of the tableau
system. If a branch arrives at atomic signed formulas without any
contradiction but where all the possible rules have been applied, then the
open tableau branch will give a model of $\left(  u_{0},u_{1}\right)
:F\varphi$, i.e., a countermodel to $\varphi$ being a partition tautology.

The $F\wedge$ rule will be illustrated by developing tableaus for two related
formulas, $\varphi_{1}=\left(  \sigma\wedge\left(  \sigma\Rightarrow
\pi\right)  \right)  \Rightarrow\left(  \sigma\wedge\pi\right)  $ and
$\varphi_{2}=\sigma\Rightarrow\left(  \left(  \sigma\Rightarrow\pi\right)
\Rightarrow\left(  \sigma\wedge\pi\right)  \right)  $, where both formulas are
subset tautologies but only the first is a partition tautology. To save space,
we have ignored the base pair and back-chain branches for $F\Rightarrow$ and
$F\wedge$ since we show that the branches with new multiple-link falsifying
chains close. Hence the base pair and back-chain branches would, a fortiori,
close since they allow even fewer possibilities to avoid contradictions. When
a formula appears on a branch with both signs, e.g., $\left(  u_{0},b\right)
:F\sigma,T\sigma$, then the branch closes as indicated with an X.

\begin{center}%
\begin{tabular}
[c]{|c|c|c|}\hline
1 & $\left(  u_{0},u_{1}\right)  :F\left[  \left(  \sigma\wedge\left(
\sigma\Rightarrow\pi\right)  \right)  \Rightarrow\left(  \sigma\wedge
\pi\right)  \right]  $ & Rules used\\\hline
2 & $\exists a,\left(  u_{0},a\right)  ,\left(  a,u_{1}\right)  :T\left(
\sigma\wedge\left(  \sigma\Rightarrow\pi\right)  \right)  ,F\left(
\sigma\wedge\pi\right)  $ & $F\Rightarrow$\\\hline
& Continuing the analysis at $\left(  u_{0},a\right)  $ & \\\hline
3 & $\exists b,c$, $\left(  b,c\right)  :F\sigma,T\left(  \sigma\wedge\left(
\sigma\Rightarrow\pi\right)  \right)  $
$\vert$
cont. & $F\wedge$ and $T$-anti-trans.\\\hline
& $\left(  b,c\right)  :F\pi,T\left(  \sigma\wedge\left(  \sigma\Rightarrow
\pi\right)  \right)  $ & $F\wedge$ and $T$-anti-trans.\\\hline
4 & $\left(  b,c\right)  :T\sigma,T\left(  \sigma\Rightarrow\pi\right)  $ X
$\vert$
$\left(  b,c\right)  :T\sigma,T\left(  \sigma\Rightarrow\pi\right)  $ &
$T\wedge$ both branches\\\hline
5 & X
$\vert$
$\left(  b,c\right)  :F\sigma$X
$\vert$%
$\vert$
$\left(  b,c\right)  :T\pi$ X & $T\Rightarrow$\\\hline
\end{tabular}

Closed tableau for: $\left(  \sigma\wedge\left(  \sigma\Rightarrow\pi\right)
\right)  \Rightarrow\left(  \sigma\wedge\pi\right)  $
\end{center}

\noindent In the second line, there was only one $T$-formula $T\left(
\sigma\wedge\left(  \sigma\Rightarrow\pi\right)  \right)  $ to transmit to a
link in the falsifying chain for $F\left(  \sigma\wedge\pi\right)  $ so a
two-link chain would suffice to give $T\left(  \sigma\wedge\left(
\sigma\Rightarrow\pi\right)  \right)  $ the alternatives of going to a
$F\sigma$ link (the left-hand alternative) or to a $F\pi$ link (the right-hand
alternative). But we use the example to illustrate a generic $u_{0},a$-chain
with a link $\left(  b,c\right)  $ in the chain. No matter how long the chain
is, there are only two alternatives created since $T\left(  \sigma
\wedge\left(  \sigma\Rightarrow\pi\right)  \right)  $ is either transmitted to
an $F\sigma$ link (the left branch) or to a $F\pi$ link (the right branch). In
the last line, two vertical lines
$\vert$%
$\vert$
were used to indicate a second branching in the right-hand branch. The use of
multiple vertical lines helps one to keep track of the level of branching in
the tree.

In the second row of the above tableau, the same Boolean conditions would hold
on $\left(  a,u_{1}\right)  $ as hold on $\left(  u_{0},a\right)  $. They are
related by an "and" and are not alternatives. Hence if contradictions can be
obtained on all branches resulting from analyzing $\left(  u_{0},a\right)
$--as indeed happened--then one does not need any more analysis on $\left(
a,u_{1}\right)  $.

When a tableau has an open branch, a branch where the formulas have been
"atomized" with no contradictions appearing and all the rules have been
exhausted, then we will see that a countermodel can be constructed using the
branch. If the formula is not a subset tautology, then one can stick entirely
to the original base pair since there is a countermodel with $\left\vert
U\right\vert =2$. But if the formula is a subset tautology but not a partition
tautology, then a multiple-link falsifying chain is required at some point.
For the formula $\sigma\Rightarrow\left(  \left(  \sigma\Rightarrow\pi\right)
\Rightarrow\left(  \sigma\wedge\pi\right)  \right)  $, which is a subset but
not partition tautology, there is an open branch where a multiple-link
falsifying chain was only used for the $F\wedge$ rule.

\begin{center}%
\begin{tabular}
[c]{|c|c|c|}\hline
1 & $\left(  u_{0},u_{1}\right)  :F\left[  \sigma\Rightarrow\left(  \left(
\sigma\Rightarrow\pi\right)  \Rightarrow\left(  \sigma\wedge\pi\right)
\right)  \right]  $ & Rules used\\\hline
2 & $\left(  u_{0},u_{1}\right)  :T\sigma,F\left(  \left(  \sigma
\Rightarrow\pi\right)  \Rightarrow\left(  \sigma\wedge\pi\right)  \right)  $ &
$F\Rightarrow$ (base pair)\\\hline
3 & $\left(  u_{0},u_{1}\right)  :T\sigma,T\left(  \sigma\Rightarrow
\pi\right)  ,F\left(  \sigma\wedge\pi\right)  $ & $F\Rightarrow$ (base
pair)\\\hline
4 & $\left(  u_{0},u_{1}\right)  :T\pi$
$\vert$
$\left(  u_{0},u_{1}\right)  :F\sigma$ X & $T\Rightarrow$\\\hline
5 & $\exists a$, $\left(  u_{0},a\right)  :F\sigma,T\left(  \sigma
\Rightarrow\pi\right)  $ and $\left(  a,u_{1}\right)  :F\pi,T\sigma$
$\vert$%
$\vert$
...%
$\vert$
X & $F\wedge$ and $T$-anti-trans.\\\hline
6 & $\left(  u_{0},a\right)  :T\pi$
$\vert$%
$\vert$%
$\vert$
$\left(  u_{0},a\right)  :F\sigma$
$\vert$%
$\vert$
...%
$\vert$
X & $T\Rightarrow$\\\hline
\end{tabular}

Simple tableau for $\sigma\Rightarrow\left(  \left(  \sigma\Rightarrow
\pi\right)  \Rightarrow\left(  \sigma\wedge\pi\right)  \right)  $ with an open
branch.{}
\end{center}

Taking the left branches at the three splittings, which terminates with
$\left(  u_{0},a\right)  :T\pi$, we can use the atomic signed formulas on each
branch to construct a "countermodel", namely a model where $\left(
u_{0},u_{1}\right)  :F\left[  \sigma\Rightarrow\left(  \left(  \sigma
\Rightarrow\pi\right)  \Rightarrow\left(  \sigma\wedge\pi\right)  \right)
\right]  $ holds so that the formula cannot be a partition tautology.

But to construct the model, the branch needs to be "completed" by applying the
eliminative rules to any signed compound formulas in the branch until signed
atomic formulas are reached, and by assigning signed atomic variables to any
remaining branches in a manner consistent with $T$-anti-transitivity and
$F$-transitivity (symmetry is assumed as a matter of course).

How does one know if this is always possible? If an assignment of signed
atomic variables to the other pairs was not possible given the signed formulas
that already have to hold at the pairs, then either there is some
contradiction that could be derived using the rules so the branch was not
really open, or the rules are incomplete (so that one has a partition
tautology where the rules were unable to close all the branches)--the latter
possibility being ruled out by the satisfaction theorem below. In the case at
hand, there is already a consistent assignment of signed atomic variables to
all the links, i.e., $\left(  u_{0},u_{1}\right)  :T\sigma,T\pi$, $\left(
u_{0},a\right)  :F\sigma,T\pi$, and $\left(  a,u_{1}\right)  :T\sigma,F\pi$.
This immediately generates the partitions $\sigma=\left\{  \left\{
u_{0},a\right\}  ,\left\{  u_{1}\right\}  \right\}  $ and $\pi=\left\{
\left\{  u_{0}\right\}  ,\left\{  u_{1},a\right\}  \right\}  $. Then
$\sigma\wedge\pi=0$, $\sigma\Rightarrow\pi=\pi$, $\left(  \sigma\Rightarrow
\pi\right)  \Rightarrow\left(  \sigma\wedge\pi\right)  =0$, and the whole
formula $\sigma\Rightarrow\left(  \left(  \sigma\Rightarrow\pi\right)
\Rightarrow\left(  \sigma\wedge\pi\right)  \right)  $ then also evaluates to
$0$ which gives a model for $\left(  u_{0},u_{1}\right)  :F\left[
\sigma\Rightarrow\left(  \left(  \sigma\Rightarrow\pi\right)  \Rightarrow
\left(  \sigma\wedge\pi\right)  \right)  \right]  $ and thus a countermodel to
that formula being a partition tautology.

Thus the two similar subset tautologies, $\left(  \sigma\wedge\left(
\sigma\Rightarrow\pi\right)  \right)  \Rightarrow\left(  \sigma\wedge
\pi\right)  $ and $\sigma\Rightarrow\left(  \left(  \sigma\Rightarrow
\pi\right)  \Rightarrow\left(  \sigma\wedge\pi\right)  \right)  $, give rather
different results for partitions since only the first formula is a partition
tautology. The difference in the two cases was that for the first formula, we
had $T\left(  \sigma\wedge\left(  \sigma\Rightarrow\pi\right)  \right)  $
being transmitted to some link in the falsifying chain for $F\left(
\sigma\wedge\pi\right)  $, where a contradiction would then arise. But in the
second formula, it was the pair of $T$-formulas, $T\sigma,T\left(
\sigma\Rightarrow\pi\right)  $, which were being transmitted so there was no
necessity that they be transmitted to the same link in the falsifying chain.
By spreading them out with $T\sigma$ going to a $F\pi$ link and $T\left(
\sigma\Rightarrow\pi\right)  $ going to a $F\sigma$ link, no contradiction
arose and in fact a countermodel could be constructed.

\subsubsection{Tableaus for partition negation}

It may be useful to also have tableau rules for negation which can be derived
from the other rules. Since we are only taking the four operations $\vee$,
$\wedge$, $\Rightarrow$, and $\mid$ as primitive, we could define the constant
$1$ as $\sigma\Rightarrow\sigma$ for any atomic variable $\sigma$ and we could
define $0$ as $1\mid1$. Then we could define negation (as in intuitionistic
logic) as $\lnot\sigma=\sigma\Rightarrow0$. But since we have the nand
operation, it is far simpler to equivalently define negation as: $\lnot
\sigma=\sigma\mid\sigma$. Then the tableau rules for negation are just a
special case of the rules for the nand.

\begin{center}%
\begin{tabular}
[c]{|c|c|c|}\hline%
\begin{tabular}
[c]{c}%
$u:T\left(  \lnot\sigma\right)  $\\\hline
$u:F\sigma$%
\end{tabular}
&
\begin{tabular}
[c]{c}%
$u:T\left(  \lnot\sigma\right)  $\\\hline
$u:F\sigma$%
\end{tabular}
&
\begin{tabular}
[c]{c}%
$\left(  u,u^{\prime}\right)  :T\left(  \lnot\sigma\right)  $\\\hline
$\left(  u,u^{\prime}\right)  :F\sigma$%
\end{tabular}
\\
Classical $T\lnot$ rule & Intuitionistic $T\lnot$ rule & Partition $T\lnot$
rule\\\hline
\end{tabular}

\end{center}

We know for the nand that four links suffice in any falsifying chain for
$F\left(  \pi\mid\sigma\right)  $, and it can easily be shown the only two
links suffice if $\sigma=\pi$. The same holds if we had defined the negation
as the implication to $0$.

\begin{center}%
\begin{tabular}
[c]{|c|c|}\hline%
\begin{tabular}
[c]{c}%
$u:F\left(  \lnot\sigma\right)  $\\\hline
$u:T\sigma$%
\end{tabular}
&
\begin{tabular}
[c]{c}%
$u:F\left(  \lnot\sigma\right)  $\\\hline
$\exists a\geq u$, $a:T\sigma$%
\end{tabular}
\\
Classical $F\lnot$ rule & Intuitionistic $F\lnot$ rule\\\hline
\end{tabular}

\begin{tabular}
[c]{|c|}\hline%
\begin{tabular}
[c]{c}%
$\left(  u,u^{\prime}\right)  :F\left(  \lnot\sigma\right)  $\\\hline
$\exists$ $u,u^{\prime}$-chain (one or two links) with $T\sigma$ on each link.
\end{tabular}
\\
Partition $F\lnot$ rule\\\hline
\end{tabular}

\end{center}

The $T\left(  \lnot\sigma\right)  $ rule is an example of a $T$-formula
implying an $F$-formula. In any such case, the $F$-formula has to hold
everywhere. If we consider any other $a\in U_{n}$, then by the $T$%
-anti-transitivity rule, $\left(  u,u^{\prime}\right)  :T\left(  \lnot
\sigma\right)  $ implies either $\left(  u,a\right)  :T\left(  \lnot
\sigma\right)  $ or $\left(  a,u^{\prime}\right)  :T\left(  \lnot
\sigma\right)  $. Whichever one holds, it implies that $F\sigma$ holds on the
link which together with $\left(  u,u^{\prime}\right)  :F\sigma$ implies that
$F\sigma$ holds on the other link by $F$-transitivity. Similarly for any other
$b\in U_{n}$, and then $\left(  a,b\right)  :F\sigma$ follows from $\left(
u,a\right)  ,\left(  u,b\right)  :F\sigma$ by $F$-transitivity where $\left(
a,b\right)  $ is any link in the complete graph $K\left(  U_{n}\right)  $.

\subsubsection{Possibility of infinite open branches: the Devil's tableau}

In the usual treatment of intuitionistic tableaus \cite{fit:il}, the elements
of the developing potential countermodel are left implicit and another device
is used to construct a countermodel when a tableau does not close. However, in
the partition tableaus we have treated the pairs $\left(  u,u^{\prime}\right)
$ quite explicitly. But then we need the $T$-anti-transitivity and
$F$-transitivity rules which do not reduce the complexity of formulas. The
cost is that we do not have the usual proof of the finiteness of tableaus
based on the fact that each of the non-structural rules for the connectives
reduces the complexity of formulas so each branch must terminate after a
finite number of steps in either a contradiction or in an open branch. That
argument is unavailable due to the two complexity-preserving rules.

Moreover, the $T$-anti-transitivity rule leads to the possibility of cycles
that can introduce an infinite sequence of stages: $U_{0}\subseteq...\subseteq
U_{n}\subseteq U_{n+1}\subseteq...$. The $F$-transitivity rule never forces
the introduction of new elements. If we had a chain $\left(  u,a\right)
,\left(  a,b\right)  ,...,\left(  c,u^{\prime}\right)  $ of elements in
$U_{n}$ with $F\sigma$ holding at each link with falsifying chains in $U_{n}$,
then we can simply hook the chains together to give a falsifying chain for
$\left(  u,u^{\prime}\right)  :F\sigma$, the conclusion of the $F$%
-transitivity rule. Hence the $F$-transitivity rule would never force new
elements to be added to $U_{n}$.

But we have seen that $T\left(  \lnot\varphi\right)  $ implies $F\varphi$, and
for an appropriate $\varphi$, the $F\varphi$ might imply new elements yield a
falsifying chain. And then the cycle repeats itself. The formula
$\sigma\Rightarrow\left(  \pi\Rightarrow\left(  \sigma\wedge\pi\right)
\right)  $ is a subset tautology that is not a partition tautology. But is it
a weak partition tautology so that its double negation would be a partition
tautology? That tableau would have the following infinite branch.

\begin{center}%
\begin{tabular}
[c]{|c|c|}\hline
$\left(  u_{0},u_{1}\right)  :F\left(  \lnot\lnot\left(  \sigma\Rightarrow
\left(  \pi\Rightarrow\left(  \sigma\wedge\pi\right)  \right)  \right)
\right)  $ & Rules used\\\hline
$\left(  u_{0},u_{1}\right)  :T\left(  \lnot\left(  \sigma\Rightarrow\left(
\pi\Rightarrow\left(  \sigma\wedge\pi\right)  \right)  \right)  \right)  $ &
$F\lnot$ (base pair)\\\hline
$\left(  u_{0},u_{1}\right)  :F\left(  \sigma\Rightarrow\left(  \pi
\Rightarrow\left(  \sigma\wedge\pi\right)  \right)  \right)  $ & $T\lnot
$\\\hline
$\exists u_{2},\left(  u_{0},u_{2}\right)  ,\left(  u_{2},u_{1}\right)
:T\sigma,F\left(  \pi\Rightarrow\left(  \sigma\wedge\pi\right)  \right)  $ &
$F\Rightarrow$\\\hline
$\left(  u_{0},u_{2}\right)  :T\pi,F\left(  \sigma\wedge\pi\right)  $ &
$F\Rightarrow$ (base pair)\\\hline
$\left(  u_{0},u_{1}\right)  :T\pi$ & $T$-anti-trans.\\\hline
$\left(  u_{0},u_{1}\right)  :F\sigma$ and $\left(  u_{2},u_{1}\right)  :F\pi$
& $F\wedge$ back-chain $u_{0},u_{1},u_{2}$\\\hline
$\left(  u_{0},u_{2}\right)  :T\left(  \lnot\left(  \sigma\Rightarrow\left(
\pi\Rightarrow\left(  \sigma\wedge\pi\right)  \right)  \right)  \right)  $ &
$T$-anti-trans.\\\hline
$\left(  u_{0},u_{2}\right)  :F\left(  \sigma\Rightarrow\left(  \pi
\Rightarrow\left(  \sigma\wedge\pi\right)  \right)  \right)  $ & $T\lnot
$\\\hline\hline
$\exists u_{3},\left(  u_{0},u_{3}\right)  ,\left(  u_{3},u_{2}\right)
:T\sigma,F\left(  \pi\Rightarrow\left(  \sigma\wedge\pi\right)  \right)  $ &
$F\Rightarrow$\\\hline
$\left(  u_{0},u_{3}\right)  :T\left(  \lnot\left(  \sigma\Rightarrow\left(
\pi\Rightarrow\left(  \sigma\wedge\pi\right)  \right)  \right)  \right)  $ &
$T$-anti-trans.\\\hline
Cycle repeats... & \\\hline
\end{tabular}

Infinite open branch in tableau
\end{center}

This tableau adds a single new element at each stage: $U_{0}=\left\{
u_{0},u_{1}\right\}  \subseteq U_{1}=\left\{  u_{0},u_{1},u_{2}\right\}
\subseteq...$ so the universe set associated with the infinite branch is the
union $U=%
{\textstyle\bigcup}
U_{n}$. For the branch to be \textit{finished}, then at each stage, each rule
needs to be applied wherever possible. For instance, at the end of stage 1
(the double line in the table), the $F$-transitivity rule could be applied to
derive $\left(  u_{2},u_{1}\right)  :F\left(  \sigma\Rightarrow\left(
\pi\Rightarrow\left(  \sigma\wedge\pi\right)  \right)  \right)  $ and $\left(
u_{0},u_{1}\right)  :F\left(  \pi\Rightarrow\left(  \sigma\wedge\pi\right)
\right)  $ but the status of $\sigma$ on $\left(  u_{0},u_{1}\right)  $ is undetermined.

To construct the countermodel, the partitions are defined on $U$ by using all
the atomic $F$-statements so that $a$ and $b$ are in the same block of the
$\alpha$ if $\left(  a,b\right)  :F\alpha$ occurred at some finite stage. If
for two elements $u,u^{\prime}\in U$, the formula $\left(  u,u^{\prime
}\right)  :F\alpha$ never occurs at any stage, then those two elements would
be in separate blocks of $\alpha$. Otherwise, there would have been a finite
$u,u^{\prime}$-chain where $\left(  a,b\right)  :F\alpha$ holds at each link
$\left(  a,b\right)  $ in the chain. But then at some finite stage, all the
links and the statements $\left(  a,b\right)  :F\alpha$ would be present so
$\left(  u,u^{\prime}\right)  :F\alpha$ would be implied by the $F$%
-transitivity rule at that stage.

By the satisfaction theorem proven below, this will provide a countermodel for
$\lnot\lnot\left(  \sigma\Rightarrow\left(  \pi\Rightarrow\left(  \sigma
\wedge\pi\right)  \right)  \right)  $. But the tableau construction of an
infinite model does not show the absence of any finite models. Indeed, the
above tableau could have been stopped at the double line, the end of the first
stage. For instance, the formula $\left(  u_{0},u_{2}\right)  :F\left(
\sigma\Rightarrow\left(  \pi\Rightarrow\left(  \sigma\wedge\pi\right)
\right)  \right)  $ was used to introduce a new element $u_{3}$ and to move to
another stage. But that formula is already satisfied at its base pair $\left(
u_{0},u_{2}\right)  $ so the introduction of a new element was unnecessary. If
we stop at the double line (after applying some more rules to "finish" that
stage), the model on $U_{1}=\left\{  u_{0},u_{1},u_{2}\right\}  $ given by the
atomic $F$-statements is: $\sigma=\left\{  \left\{  u_{0},u_{1}\right\}
,\left\{  u_{2}\right\}  \right\}  $ and $\pi=\left\{  \left\{  u_{0}\right\}
,\left\{  u_{1},u_{2}\right\}  \right\}  $ which, in this case, provides a
countermodel for $\lnot\lnot\left(  \sigma\Rightarrow\left(  \pi
\Rightarrow\left(  \sigma\wedge\pi\right)  \right)  \right)  .$ This shows,
incidentally, that $\sigma\Rightarrow\left(  \pi\Rightarrow\left(
\sigma\wedge\pi\right)  \right)  $ is not even a weak partition tautology.

It is easy to see why this sort of an infinite branch generated by a simple
cycle was unnecessary. On the links of the chain introduced by the new element
$a$, the Boolean conditions for $F\left(  \sigma\Rightarrow\left(
\pi\Rightarrow\left(  \sigma\wedge\pi\right)  \right)  \right)  $ had to hold.
But when $T\left(  \lnot\left(  \sigma\Rightarrow\left(  \pi\Rightarrow\left(
\sigma\wedge\pi\right)  \right)  \right)  \right)  $ was sent to one of the
links and $F\left(  \sigma\Rightarrow\left(  \pi\Rightarrow\left(
\sigma\wedge\pi\right)  \right)  \right)  $ again derived, then its Boolean
conditions would hold at that link so it was unnecessary to introduce a new element.

There is a much more devilish pattern that can generate an infinite branch, a
pattern we might call the "Devil's tableau." The idea is to take two formulas
with complementary Boolean conditions, such as $F\left(  \sigma\wedge
\pi\right)  $ and $F\left(  \sigma\mid\pi\right)  $, where one or both might
introduce new elements. Thus one or the other of the formulas would not have
their Boolean conditions satisfied at the base pair. In the following Devil's
tableau, we develop an infinite open branch taking the set of elements being
introduced as the natural numbers $%
\mathbb{N}
$.

\begin{center}%
\begin{tabular}
[c]{|c|c|}\hline
$\left(  0,1\right)  :F\left(  \lnot\lnot\left[  \left(  \sigma\wedge
\tau\right)  \vee\left(  \sigma\mid\tau\right)  \right]  \right)  $ & Rules
used\\\hline
$\left(  0,1\right)  :T\left(  \lnot\left[  \left(  \sigma\wedge\tau\right)
\vee\left(  \sigma\mid\tau\right)  \right]  \right)  $ & $F\lnot$\\\hline
$\left(  0,1\right)  :F\left[  \left(  \sigma\wedge\tau\right)  \vee\left(
\sigma\mid\tau\right)  \right]  $ & $T\lnot$\\\hline
$\left(  0,1\right)  :F\left(  \sigma\wedge\tau\right)  ,F\left(  \sigma
\mid\tau\right)  $ & $F\vee$\\\hline
$\left(  0,1\right)  :T\sigma,T\tau$ & $F\mid$ (base pair)\\\hline
$\exists2,\left(  0,2\right)  :F\sigma$ and $\left(  1,2\right)  :F\tau$ &
$F\wedge$\\\hline
$\left(  0,2\right)  :F\left(  \sigma\wedge\tau\right)  ,F\left(  \sigma
\mid\tau\right)  $ & $T$-anti-trans. etc.\\\hline
$\exists3,\left(  0,3\right)  ,\left(  2,3\right)  :T\sigma,T\tau$ & $F\mid
$\\\hline
$\left(  0,3\right)  :F\left(  \sigma\wedge\tau\right)  ,F\left(  \sigma
\mid\tau\right)  $ & $T$-anti-trans. etc.\\\hline
$\exists4,\left(  0,4\right)  :F\sigma$ and $\left(  3,4\right)  :F\tau$ &
$F\wedge$\\\hline
$\left(  0,4\right)  :F\left(  \sigma\wedge\tau\right)  ,F\left(  \sigma
\mid\tau\right)  $ & $T$-anti-trans. etc.\\\hline\hline
$\exists5,\left(  0,5\right)  ,\left(  4,5\right)  :T\sigma,T\tau$ & $F\mid
$\\\hline
... & \\\hline
\end{tabular}

Infinite open branch of a Devil's tableau
\end{center}

The even stages $U_{0}=\left\{  0,1\right\}  $, $U_{2}=\left\{
0,1,2,3\right\}  $, ... use $F\wedge$ to introduce a new element and the odd
stages use $F\mid$ to introduce a new element. This generates the pattern

\begin{center}
$\left(  0,even\right)  :F\sigma$ and $\left(  even-1,even\right)  :F\tau$ and
$\left(  0,even-1\right)  :T\sigma,T\tau$

$\left(  0,odd\right)  ,\left(  odd-1,odd\right)  :T\sigma,T\tau$ and $\left(
0,odd-1\right)  =\left(  0,even\right)  :F\sigma$.
\end{center}

The union of the stages $U_{n}=\left\{  0,1,...,n+1\right\}  $ is the natural
numbers $%
\mathbb{N}
$ and the partitions defined by the atomic $F$-statements are:

\begin{center}
$\sigma=\left\{  \left\{  0,2,4,6,...\right\}  ,\left\{  1\right\}  ,\left\{
3\right\}  ,\left\{  5\right\}  ,...\right\}  $

$\tau=\left\{  \left\{  0\right\}  ,\left\{  1,2\right\}  ,\left\{
3,4\right\}  ,\left\{  5,6\right\}  ,...\right\}  $.
\end{center}

This is indeed a model since $\sigma\wedge\tau=0=\sigma\mid\tau$. The fact
that $\sigma\wedge\tau=0$ is easily seen since $\sigma$ identifies all the
even numbers and $\tau$ identifies each odd number with its successor even
number. To see that $\sigma\mid\tau=0$, consider its graph which will have
links $n\sim m$ whenever $n$ and $m$ are distinguished by both partitions.
Thus in that graph $even\sim even+1$ ($=odd$) and $odd\sim odd+2$ so there is
a finite chain connecting any $n,m\in%
\mathbb{N}
$.

By alternating between the two potentially element-introducing $F$-formulas,
$F\left(  \sigma\wedge\tau\right)  $ and $F\left(  \sigma\mid\tau\right)  $,
the Devil's tableau avoids having both formulas satisfied at the base pair at
the same time. But there is still the possibility that both formulas could be
satisfied by back-chains at the same time--so that there would be no need to
introduce any new constants and the branch could be terminated there. Indeed,
that is the case with this Devil's tableau. If we stop the tableau at the
double line where the stage is $U_{3}=\left\{  0,1,2,3,4\right\}  $, then the
partitions are: $\sigma=\left\{  \left\{  0,2,4\right\}  ,\left\{  1\right\}
,\left\{  3\right\}  \right\}  $ and $\tau=\left\{  \left\{  0\right\}
,\left\{  1,2\right\}  ,\left\{  3,4\right\}  \right\}  $ and $\sigma
\wedge\tau=0=\sigma\mid\tau$ as well. Hence in this case, the consideration of
back-chains gives a finite tableau that provides a finite countermodel, but
the question of whether there is always a finite countermodel is left open
along with the related question of the decidability of the set of partition
tautologies. The necessity of considering back-chains in order to have a
finite open branch of the Devil's tableau shows why back-chains are included
in the "$\exists u,u^{\prime}$-chain" clause in the conclusions of the
element-introducing $F$ rules.

\subsubsection{More proofs and countermodels using tableaus}

A few more examples may be helpful. The partition tautology $\pi
\Rightarrow\left(  \sigma\Rightarrow\pi\right)  $ provides a simple example.
But even for this example, tableau trees expand rapidly without shortcuts and
symmetry arguments.

\begin{center}%
\begin{tabular}
[c]{|c|c|}\hline
$\left(  u_{0},u_{1}\right)  :F\left[  \pi\Rightarrow\left(  \sigma
\Rightarrow\pi\right)  \right]  $ & Rules used:\\\hline
$\left(  u_{0},u_{1}\right)  :T\pi,F\left(  \sigma\Rightarrow\pi\right)  $
$\vert$%
& $F\Rightarrow$\\\hline
$\exists a$ $\left(  u_{0},a\right)  ,\left(  a,u_{1}\right)  :T\pi,F\left(
\sigma\Rightarrow\pi\right)  $ and $\left(  u_{0},u_{1}\right)  :F\pi$ &
$F\Rightarrow$, $F$-trans.\\\hline
(continuing with the right branch) & \\\hline
...%
$\vert$
$\left(  u_{0},a\right)  :T\sigma,F\pi$X
$\vert$%
$\vert$
$\exists b,\left(  u_{0},b\right)  ,\left(  b,a\right)  :T\sigma,F\pi$ and
$\left(  u_{0},a\right)  :F\pi$X & $F\Rightarrow$, $F$-trans.\\\hline
(right branch closed so picking up the left branch) & \\\hline
$\left(  u_{0},u_{1}\right)  :T\pi,F\left(  \sigma\Rightarrow\pi\right)  $
$\vert$
X
$\vert$%
$\vert$
X & \\\hline
$\left(  u_{0},u_{1}\right)  :T\sigma,F\pi$X
$\vert$%
$\vert$
cont. & $F\Rightarrow$\\\hline
$\exists a,\left(  u_{0},a\right)  ,\left(  a,u_{1}\right)  :T\sigma,F\pi$ and
$\left(  u_{0},u_{1}\right)  :F\pi$ X
$\vert$
X
$\vert$%
$\vert$
X & $F\Rightarrow$, $F$-trans.\\\hline
\end{tabular}

Closed tableau for:\ $\pi\Rightarrow\left(  \sigma\Rightarrow\pi\right)  $
\end{center}

\noindent Taking the leftmost branches, we stay at the base pair $\left(
u_{0},u_{1}\right)  $ and have essentially the classical closing tableau since
this formula is a subset tautology. Since the other element-introducing
branches also close, the formula is a partition tautology (assuming the
correctness theorem proved below).

Peirce's law, $((\sigma\Rightarrow\pi)\Rightarrow\sigma)\Rightarrow\sigma$, is
a good example of non-closing tableau which must generate a model where the
formula does not distinguish some pair.

\begin{center}%
\begin{tabular}
[c]{|c|c|}\hline
$\left(  u_{0},u_{1}\right)  :F\left[  ((\sigma\Rightarrow\pi)\Rightarrow
\sigma)\Rightarrow\sigma\right]  $ & Rules used\\\hline
$\left(  u_{0},u_{1}\right)  :T\left[  ((\sigma\Rightarrow\pi)\Rightarrow
\sigma)\right]  ,F\sigma$ & $F\Rightarrow$ (base)\\\hline
$\left(  u_{0},u_{1}\right)  :F\left(  \sigma\Rightarrow\pi\right)  $
$\vert$
$\left(  u_{0},u_{1}\right)  :T\sigma$ X & $T\Rightarrow$\\\hline
$\exists a,\left(  u_{0},a\right)  ,\left(  a,u_{1}\right)  :T\sigma,F\pi$ and
$\left(  u_{0},u_{1}\right)  :F\pi$%
$\vert$
X & $F\Rightarrow$, $F$-trans.\\\hline
$\left(  u_{0},a\right)  :T\left[  ((\sigma\Rightarrow\pi)\Rightarrow
\sigma)\right]  $%
$\vert$
X & $T$-anti-trans.\\\hline
$\left(  u_{0},a\right)  :F\left(  \sigma\Rightarrow\pi\right)  $
$\vert$%
$\vert$
$\left(  u_{0},a\right)  :T\sigma$
$\vert$
X & $T\Rightarrow$\\\hline
\end{tabular}

Non-closed tableau for: $((\sigma\Rightarrow\pi)\Rightarrow\sigma
)\Rightarrow\sigma$
\end{center}

\noindent The branch terminating with $\left(  u_{0},a\right)  :T\sigma$ in
the last row is an open branch (atomic formulas with no contradiction) so it
may be used to generate of model of $F\left[  ((\sigma\Rightarrow
\pi)\Rightarrow\sigma)\Rightarrow\sigma\right]  $, i.e., a countermodel to
Peirce's law being a partition tautology. To generate the model, we need to
fill out the atomic signed formulas on all the links but that is already done
on the indicated branch. The universe set is the three elements used in the
tableau: $U=\left\{  u_{0},u_{1},a\right\}  $. The partition $\sigma$ has
$\left(  u_{0},u_{1}\right)  :F\sigma$ while $T\sigma$ holds at $\left(
u_{0},a\right)  $ and $\left(  a,u_{1}\right)  $. Thus $\sigma=\left\{
\left\{  u_{0},u_{1}\right\}  ,\left\{  a\right\}  \right\}  $. The partition
$\pi$ has $F\pi$ on all links so $\pi$ is the blob: $\pi=\left\{  \left\{
u_{0},u_{1},a\right\}  \right\}  =0$. The compound partitions are then:
$\sigma\Rightarrow\pi=\pi=0$ (since no non-singleton block of $\pi$ is
contained in a block of $\sigma$), $\left(  \sigma\Rightarrow\pi\right)
\Rightarrow\sigma=1$ (since all blocks of $\sigma$ are contained in the blob),
and finally $\left(  \left(  \sigma\Rightarrow\pi\right)  \Rightarrow
\sigma\right)  \Rightarrow\sigma=\sigma$ (since $1\Rightarrow\sigma=\sigma$)
so that $\left(  u_{0},u_{1}\right)  :F\left[  ((\sigma\Rightarrow
\pi)\Rightarrow\sigma)\Rightarrow\sigma\right]  $ holds and Peirce's law is
not a partition tautology.

Essentially the same argument as in the common-dits theorem yields a powerful
result that can be used to close branches of a tableau. It gives conditions
under which a contradiction has to exist on some link without forcing one to
work through all the possibilities on sub-branches to show they close.

\begin{lemma}
[Branch-closing lemma]Suppose $\left(  a,b\right)  :T\tau,F\pi$ and $\left(
c,d\right)  :T\varphi,F\pi$ where there is a chain connecting the two links
that has $F\pi$ holding at each link on the chain. Then there exists a link
where $T\tau,T\varphi,F\pi$ all hold on the link.
\end{lemma}

\noindent Proof: The $F\pi$-chain needs to connect $a$ or $b$ with $c$ or $d
$. If it connects, say, $a$ and $c$, then by $F$-transitivity, $\left(
a,c\right)  :F\pi$. Then we have the following situation regarding those four points.%

\begin{center}
\includegraphics[
natheight=76.027397bp,
natwidth=195.018707bp,
height=78.3125pt,
width=197.6875pt
]%
{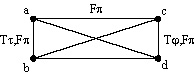}%
\\
Figure 7: Branch-closing lemma
\end{center}

\noindent By $F$-transitivity again, $F\pi$ has to hold at all the links
connecting the four points. Consider the triangle formed by $a$, $b$, and $c$.
By $T$-anti-transitivity, $T\tau$ has to hold on either $\left(  a,c\right)  $
or $\left(  b,c\right)  $. Case 1. If $T\tau$ holds on $\left(  a,c\right)  $,
then by considering the triangle formed by $a$, $c$, and $d$, then $T\tau$ has
to hold at $\left(  a,d\right)  $ or $\left(  c,d\right)  $. If it holds at
$\left(  c,d\right)  $, then we are finished so suppose it holds on $\left(
a,d\right)  $. But since $T\varphi$ holds on $\left(  c,d\right)  $, by
$T$-anti-transitivity again, $T\varphi$ has to hold at either $\left(
a,c\right)  $ or $\left(  a,d\right)  $ and we are finished in either case.
Case 2. If $T\tau$ holds on $\left(  b,c\right)  $, then we repeat the same
argument but for the triangle formed by $b$, $c$, and $d$. $\blacksquare$

The principal use of the branch-closing lemma is when, on a branch of a
tableau, we have signed formulas $T\sigma,F\pi$ on one link, $T\left(
\sigma\Rightarrow\pi\right)  ,F\pi$ on another link, with a chain connecting
the two links with $F\pi$ holding on each link of the chain. Then by the
branch-closing lemma, there exists a link where $T\sigma,T\left(
\sigma\Rightarrow\pi\right)  ,F\pi$ all hold and thus the branch closes since
there is a contradiction on that link regardless of whether $T\left(
\sigma\Rightarrow\pi\right)  $ is developed as $F\sigma$ or $T\pi$ by the
$T\Rightarrow$ rule.

The single $\pi$-negation transform, $\overset{\pi}{\lnot}\sigma
\vee\overset{\pi}{\lnot}\overset{\pi}{\lnot}\sigma$, of the law of excluded
middle, $\sigma\vee\lnot\sigma$, is an example of a partition tautology that
is not an intuitionistic validity. It is the $\pi$-negation version of
$\lnot\sigma\vee\lnot\lnot\sigma$, the weak law of excluded middle. The
tableau proof of the formula is also an example of using the branch-closing lemma.

\begin{center}%
\begin{tabular}
[c]{|c|c|}\hline
$\left(  u_{0},u_{1}\right)  :F\left[  \left(  \sigma\Rightarrow\pi\right)
\vee\left(  \left(  \sigma\Rightarrow\pi\right)  \Rightarrow\pi\right)
\right]  $ & Rules used\\\hline
$\left(  u_{0},u_{1}\right)  :F\left(  \sigma\Rightarrow\pi\right)  ,F\left[
\left(  \sigma\Rightarrow\pi\right)  \Rightarrow\pi\right]  $ & $F\vee
$\\\hline
$\left(  u_{0},u_{1}\right)  :T\sigma,F\pi$
$\vert$
$\exists a,$ $\left(  u_{0},a\right)  ,(a,u_{1}):T\sigma,F\pi$ and $\left(
u_{0},u_{1}\right)  :F\pi$ & $F\Rightarrow$, $F$-trans.\\\hline
...
$\vert$
$\exists b,$ $\left(  u_{0},b\right)  ,\left(  b,u_{1}\right)  :T\left(
\sigma\Rightarrow\pi\right)  ,F\pi$ and $\left(  u_{0},u_{1}\right)  :F\pi$ &
$F\Rightarrow$, $F$-trans.\\\hline
...%
$\vert$
$F\pi$ holds on chain $u_{0},u_{1},a,b$ so $\left(  u_{0},a\right)  :T\sigma$
(cont.) & \\\hline
and $\left(  u_{0},b\right)  :T\left(  \sigma\Rightarrow\pi\right)  $ collide
with $F\pi$. X & B-C lemma\\\hline
\end{tabular}

Closed tableau for: $\overset{\pi}{\lnot}\sigma\vee\overset{\pi}{\lnot
}\overset{\pi}{\lnot}\sigma$
\end{center}

\noindent The branch with both falsifying chains allows maximal freedom from
contradiction but it still closes (by the branch-closing lemma) so the left
branch starting with a base-pair application of $F\Rightarrow$ would, a
fortiori, close.

To see why this formula is not intuitionistically valid, we could develop its
intuitionistic tableau. In the partition case, we have used repeatedly the
fact that when $F\left(  \sigma\Rightarrow\pi\right)  $ is satisfied by a
falsifying chain, then $F$-transitivity implies that $F\pi$ has to hold at the
base pair. A similar result holds in the intuitionistic case. By the
$F\Rightarrow$ rule, $u:F\left(  \sigma\Rightarrow\pi\right)  $ implies that
the Boolean condition $T\sigma,F\pi$ has to hold at some higher point $a\geq
u$. But if $u:T\pi$ held, then the structural rule would imply that $T\pi$ had
to hold at all higher points (contradicting $a:F\pi$), so $u:F\pi$ must hold.

\begin{center}%
\begin{tabular}
[c]{|c|c|}\hline
$u:F\left[  \left(  \sigma\Rightarrow\pi\right)  \vee\left(  \left(
\sigma\Rightarrow\pi\right)  \Rightarrow\pi\right)  \right]  $ & Rules
used\\\hline
$u:F\left(  \sigma\Rightarrow\pi\right)  ,F\left[  \left(  \sigma
\Rightarrow\pi\right)  \Rightarrow\pi\right]  $ & $F\vee$\\\hline
$\exists a\geq u,$ $a:T\sigma,F\pi$ and $\exists b\geq u,$ $b:T\left(
\sigma\Rightarrow\pi\right)  ,F\pi$ & $F\Rightarrow$ twice\\\hline
$b:F\sigma$
$\vert$
$b:T\pi$ X & $T\Rightarrow$\\\hline
\end{tabular}

Open intuitionistic tableau for: $\overset{\pi}{\lnot}\sigma\vee
\overset{\pi}{\lnot}\overset{\pi}{\lnot}\sigma$
\end{center}

\noindent As with partition tableaus, a model can be constructed from an open
branch of an intuitionistic tableau. There are three points in $U=\left\{
u,a,b\right\}  $ and the partial ordering is given by $u\leq a$ and $u\leq b$.
Then $F\pi$ holds at all points so $\pi$ is modeled by the empty set
$\emptyset$. $T\sigma$ holds at $a$ but cannot hold at $b$ and thus cannot
hold at $u$. Hence $\sigma$ is modeled by the up-closed set $\left\{
a\right\}  $. The sets formed by the connectives are then: $\sigma
\Rightarrow\pi=\operatorname*{int}\left(  \sigma^{c}\cup\pi\right)  =\left\{
b\right\}  $ and $\left(  \left(  \sigma\Rightarrow\pi\right)  \Rightarrow
\pi\right)  =\operatorname*{int}\left(  \left\{  b\right\}  ^{c}\cup
\pi\right)  =\left\{  a\right\}  $ so that: $\left(  \sigma\Rightarrow
\pi\right)  \vee\left(  \left(  \sigma\Rightarrow\pi\right)  \Rightarrow
\pi\right)  =\left\{  b\right\}  \cup\left\{  a\right\}  =\left\{
a,b\right\}  \not =U$ and we have a model for $u:F\left[  \left(
\sigma\Rightarrow\pi\right)  \vee\left(  \left(  \sigma\Rightarrow\pi\right)
\Rightarrow\pi\right)  \right]  $.

The reason why the intuitionistic tableau does not close is that once $u$
branches to the two points $a$ and $b$, those branches in the ordered set $U$
do not need to interact so the "conflict" between the two branches in the
ordering never gives a contradiction to close the tableau. However in
partition logic, for any two links, there is always a direct connection so the
conflict becomes a contradiction. For instance, in the partition tableau for
this formula, the potential conflict at the two separate links $\left(
u_{0},a\right)  :T\sigma,F\pi$ and $\left(  b,u_{1}\right)  :T\left(
\sigma\Rightarrow\pi\right)  ,F\pi$ is connected by the link $\left(
u_{0},u_{1}\right)  :F\pi$ so the branch-closing lemma brings out the contradiction.

\subsection{Correctness theorem for partition tableaus}

A tableau for $\varphi$, i.e., a tableau with the root $\left(  u_{0}%
,u_{1}\right)  :F\varphi$, \textit{closes} if all the possible branches
terminate with a contradiction $\left(  a,b\right)  :T\pi,F\pi$ at some pair
$\left(  a,b\right)  $ for some subformula $\pi$. But this definition requires
special attention to the $F\wedge$ rule. If a branch does not close with a
contradiction, then the branch should generate a countermodel which requires
any element-introducing use of the $F\wedge$ rule to introduce specific
elements in the falsifying chain. But if a branch is to close with a
contradiction for each alternative, then it is not enough to have it close
from some finite set of specific falsifying chains since there is an infinite
set of possible finite falsifying chains (and the $F\wedge$ rule would not
have the finite-branching property). This is why some special attention is
required in a tableau that uses the $F\wedge$ and that closes. The $F\wedge$
rule is interpreted as only introducing a generic finite chain of finite
length, and the links in the chain only become specific when the
$T$-anti-transitivity rule transmits a $T$-formula to some link. By taking it
to be the shortest falsifying chain we could ensure that the links are
alternating. Thus if $\left(  u,u^{\prime}\right)  :F\left(  \sigma\wedge
\pi\right)  $, then the links would alternate between $F\sigma,T\pi$ and
$T\sigma,F\pi$. If there were, say, three other $T$-formulas, $T\phi_{1}$,
$T\phi_{2}$, and $T\phi_{3}$, holding at $\left(  u,u^{\prime}\right)  $, then
each $T\phi_{i}$ could be transmitted to either a $F\sigma$ link or a $F\pi$
link (and always to different $F\sigma$ or $F\pi$ links from the other
$T\phi_{j}$ formulas).\footnote{If a branch would close when the $T\phi_{i}$
formulas were spread out on different links, then it would, a fortiori, close
when some of the formulas were bunched together on the same type of link of
the falsifying chain, so those alternatives may be ignored.} Hence there are
only $2^{3}=8$ branches generated by the $F\wedge$ that would ultimately need
to close for the tableau to close.

The correctness theorem for tableaus asserts that if the tableau for
$F\varphi$ closes, then $\varphi$ is a partition tautology, and the
completeness theorem proves the reverse. The strategy of the proof of the
correctness theorem is to show that if there is an interpretation in $\Pi(U)$
of the premise of a tableau rule, then there is an interpretation of the
conclusion. Hence if the tableau closes, then since there can be no
interpretation of the conclusions that close a tableau, there can be no
interpretation of the beginning of the tableau, $\left(  u_{0},u_{1}\right)
:F\varphi$ and thus $\varphi$ is a partition tautology.

An \textit{interpretation} or \textit{model} of the formulas has a universe
set $U$ with two or more elements, interprets the atomic variables as
partitions on $U$, and interprets the operation symbols $\vee$, $\wedge$,
$\Rightarrow$, and $\mid$ as those operations in the partition algebra $\Pi(U)
$. When convenient, we use the dit-set representation of $\Pi(U)$ so the
variables and formula would refer to dit sets or partition relations rather
than set-of-blocks partitions. Statements like $B\in\pi$ are interpreted in
the obvious manner without pedantically saying that $B$ is a block in the
partition interpreting the symbol "$\pi$" and so forth. We are also already
accustomed to using statements like "$T\varphi$ holds at $\left(  u,u^{\prime
}\right)  $" as saying that $\left(  u,u^{\prime}\right)  $ is a distinction
of the partition interpreting $\varphi$, and similarly for $F$ statements.

\begin{theorem}
[Correctness of partition tableaus]If the tableau for $F\varphi$ closes, then
$\varphi$ is a partition tautology.
\end{theorem}

\noindent Proof: We assume we have an interpretation of the formulas in a
universe set $U$ where the premises of the rules hold, and then we show that
one of the possible conclusions holds.

All the $T$ rules can be handled in a uniform way. Where $\ast$ is $\vee$,
$\wedge$, $\Rightarrow$, or $\mid$, if $T\left(  \sigma\ast\pi\right)  $ holds
at $\left(  u,u^{\prime}\right)  $, then the Boolean conditions for $T\left(
\sigma\ast\pi\right)  $ must hold at some link on any $u,u^{\prime}$-chain
which means they must hold at the one-link chain $\left(  u,u^{\prime}\right)
$ which are the conclusions in the four $T$ rules.

All the $F$ rules have the general form that the premise $\left(  u,u^{\prime
}\right)  :F\left(  \sigma\ast\pi\right)  $ implies the existence of a
$u,u^{\prime}$-chain where the Boolean conditions for $F\left(  \sigma\ast
\pi\right)  $ hold at every link of the chain. The assumption is that at a
certain stage where the set of elements or "constants" is $U_{n}$, then the
elements $u$ and $u^{\prime}$ of $U_{n}$ are interpreted in $U$ and there are
partitions on $U$ interpreting the atomic variables so that $F\left(
\sigma\ast\pi\right)  $ holds at $u,u^{\prime}\in U_{n}$. Then by the
falsifying-chain theorem, there is a finite $u,u^{\prime}$-chain of elements
of $U$ where the Boolean conditions for $F\left(  \sigma\ast\pi\right)  $ hold
at each link. In terms of $U_{n}$, that falsifying chain could be a
back-chain, a mixed chain, or a chain of new elements linking $u$ and
$u^{\prime}$. Thus by adding a finite number of new elements of $U$ to $U_{n}$
if necessary to have $U_{n+1}$, one of the alternatives of the $F\ast$ rule is
the set of assignments to the links of that chain that hold in the model on
$U$.

The structural rules are also correct by similar reasoning. \noindent In any
interpretation, $\left(  u,u^{\prime}\right)  :T\varphi$ means that $\left(
u,u^{\prime}\right)  \in\operatorname*{dit}\left(  \varphi\right)  $ which is
a partition relation and thus anti-transitive so the conclusion of the $T$
anti-transitivity rule holds. If $\left(  u,a\right)  ,\left(  a,u^{\prime
}\right)  :F\varphi$, then $\left(  u,a\right)  ,\left(  a,u^{\prime}\right)
\in\operatorname*{indit}\left(  \varphi\right)  $ which is an equivalence
relation so its transitivity gives the conclusion of the $F$ transitivity
rule. \noindent In any interpretation, both dit sets and indit sets are
symmetric so if the premise holds, then the conclusion holds in each of the
symmetry rules.

Hence if the premise in any of the rules has an interpretation, then so does
one of the alternatives in the conclusion. Since the conclusions of the closed
branches have no interpretation, a closed tableau for $F\varphi$ implies there
is no interpretation for the premise of $\left(  u_{0},u_{1}\right)
:F\varphi$ so that $\left(  u_{0},u_{1}\right)  :T\varphi$ holds for any pair
in any interpretation and thus $\varphi$ is a partition tautology.
$\blacksquare$

\subsection{Completeness theorem for partition tableaus}

\subsubsection{Completing a tableau}

The correctness theorem shows that if all branches in the tableau with the
root $\left(  u_{0},u_{1}\right)  :F\varphi$ close, then $\varphi$ is a
partition tautology. The goal now is to prove the converse: if $\varphi$ is a
partition tautology, then there is a tableau for $F\varphi$ where all branches
close, i.e., $\varphi$ is a theorem of the tableau system. It would be
equivalent to prove the contrapositive that if there was an open branch (i.e.,
a branch that could not be closed), then the branch would provide a
countermodel to $\varphi$, i.e., a model for $\left(  u_{0},u_{1}\right)
:F\varphi$.

A branch of a tableau is \textit{closed} if for some pair $\left(  a,b\right)
$ and some formula $\pi$, both $\left(  a,b\right)  :T\pi$ and $\left(
a,b\right)  :F\pi$ occur on the branch. In terms of stages, a closed branch
must close at some finite stage and the branch terminates at that stage. If a
tableau with the root $\left(  u_{0},u_{1}\right)  :F\varphi$ is
\textit{closed} in the sense that all branches are closed, then, since all
rules are finitely-branching (using the generic falsifying chain in the
$F\wedge$ rule which only branches for the finite number of possible ways that
$T$-formulas $T\phi_{i}$ could be transmitted to the chain), a closed tableau
is finite and thus constitutes a tableau proof of $\varphi$.

A branch of a tableau is \textit{complete at stage }$n$ with the universe set
$U_{n}$ if all applications of the rules that can be made have been made.
There are two types of rules, the connective rules for the four connectives
and the structural rules ($T$-anti-transitivity, $F$-transitivity, and the
symmetry rules). When a connective rule with a premise $\left(  u,u^{\prime
}\right)  :F\phi$ or $\left(  u,u^{\prime}\right)  :T\phi$ has been used then
it can be checked ($\checkmark$) once. But the same premise could also be used
in the premise for many structural rules so a premise would get a second check
mark when all the structural rules have been applied at that stage. In the
order of applying rules systematically, the structural rules and the
non-element-introducing rules should be used first and should involve only
elements from the universe set $U_{n}$ at that stage. Then the potentially
element-introducing rules ($F\Rightarrow$, $F\wedge$, and $F\mid$) are used.
Then any applicable rules may need to be applied again if any new formulas
$\left(  a,b\right)  :F\phi$ or $\left(  a,b\right)  :T\phi$ were introduced
for old elements $a,b\in U_{n}$. This cycling over the rules at each stage
terminates after a finite number of steps since there are a finite number of
elements in each $U_{n}$ and we are not yet considering any new elements
introduced into the next stage. No infinite regress (like in the Devil's
tableau) is possible since we are only considering formulas at pairs of
elements of the given finite stage.

Being "complete" is defined stage by stage since when new elements are
introduced, there is a new stage and new applications of the structural rules
of $T$-anti-transitivity and $F$-transitivity to premises of former stages may
occur. Thus at each new stage, the second check mark on the old formulas is
erased until all the new applications using pairs involving new elements have
been used. For instance, the premise $\left(  u,u^{\prime}\right)  :T\phi$
might be checked a second time when applied to all $a\in U_{n}$ to yield
$\left(  u,a\right)  :T\phi$ or $\left(  a,u^{\prime}\right)  :T\phi$ but then
could be applied again using new $b\in U_{n+1}$.

When a stage is complete but new elements were introduced, then the same
process continues at the next stage. If a stage is completed with no new
elements introduced, then the branch is complete (with no further stages). A
branch of the tableau is \textit{complete} when it is complete at each of its
stages. A tableau is \textit{completed} if every branch is either complete or closed.

\subsubsection{Satisfaction and completeness theorems}

A completed tableau that is not closed must have at least one open complete branch.

\begin{theorem}
[Satisfaction theorem]An open complete branch of a partition tableau with the
root formula $\left(  u_{0},u_{1}\right)  :F\varphi$ gives a model where the
root formula is satisfied, i.e., $\left(  u_{0},u_{1}\right)  \in
\operatorname*{indit}\left(  \varphi\right)  $ in the model.
\end{theorem}

\noindent Proof: An open complete branch of a tableau will be used to define a
model on a set $U$. If the complete open branch terminated at the stage
$U_{n}$, then take $U=U_{n}$. Otherwise, there is an infinite sequence of
stages $U_{0}\subseteq U_{1}\subseteq...$ and $U=%
{\textstyle\bigcup\limits_{n}}
U_{n}$. The partitions on $U$ are defined by the formulas $\left(  a,b\right)
:F\alpha$ occurring in the branch for the atomic variables $\alpha$ occurring
in the root formula $\varphi$. Using the graph machinery, these atomic
$F$-formulas occurring in the branch define the links of a graph on the node
set $U$, and the blocks of the partition $\alpha$ are the sets of nodes in the
connected components of the graph. This defines the partitions interpreting
the atomic variables of $\varphi$ and then the partition operations of
$\Pi\left(  U\right)  $ will give an interpretation of $\varphi$ using
partitions on $U$. We need to show that $\left(  u_{0},u_{1}\right)
\in\operatorname*{indit}\left(  \varphi\right)  $ under that interpretation.

The proof is by induction over the complexity of the subformulas of $\varphi$.
The basis step is that every signed atomic formula which occurs in the branch
is true in the model. If $\left(  u,u^{\prime}\right)  :F\alpha$ occurs in the
branch then it is true by definition in the model, i.e., $\left(  u,u^{\prime
}\right)  \in\operatorname*{indit}\left(  \alpha\right)  $. If $\left(
u,u^{\prime}\right)  :T\alpha$ occurs in the branch but does not hold in the
model, i.e., $\left(  u,u^{\prime}\right)  \in\operatorname*{indit}\left(
\alpha\right)  $, then using the graph constructed for $\alpha$ and using the
falsifying-chain theorem, there is a finite $u,u^{\prime}$-chain with $\left(
u_{i},u_{i+1}\right)  :F\alpha$ holding at each link. Moreover, there is a
finite stage $U_{n}$ where all these formulas would have occurred. But
completeness at that stage would then imply, by using the $F$ -transitivity
rule, that the formula $\left(  u,u^{\prime}\right)  :F\alpha$ held at that
stage which would contradict $\left(  u,u^{\prime}\right)  :T\alpha$ holding
at some stage on the complete open branch. Hence if $\left(  u,u^{\prime
}\right)  :T\alpha$ did occur in the open branch, then $\left(  u,u^{\prime
}\right)  \in\operatorname*{dit}\left(  \alpha\right)  $ in the model.

The induction steps can be efficiently treated using the graph machinery.
Suppose $\left(  u,u^{\prime}\right)  :T\left(  \sigma\ast\pi\right)  $ occurs
in the complete open branch. In order for $\left(  u,u^{\prime}\right)
\in\operatorname*{dit}\left(  \sigma\ast\pi\right)  $ in the model on $U$,
then for every finite $u,u^{\prime}$-chain in $U$, the Boolean conditions for
$T\left(  \sigma\ast\pi\right)  $ must hold at some link in the chain. Suppose
not, so there is a $u,u^{\prime}$-chain where the complementary Boolean
conditions for $F\left(  \sigma\ast\pi\right)  $ hold at each link. There is a
finite stage $U_{n}$ of the branch in which all the elements of that chain
have appeared and where $\left(  u,u^{\prime}\right)  :T\left(  \sigma\ast
\pi\right)  $ also occurs. But then by the completeness of applying the
$T$-anti-transitivity at that stage, there is a link $\left(  u_{i}%
,u_{i+1}\right)  $ in the chain where $\left(  u_{i},u_{i+1}\right)  :T\left(
\sigma\ast\pi\right)  $ holds. Then by completeness and the connective rule
for $T\left(  \sigma\ast\pi\right)  $, the formulas for the Boolean conditions
for $T\left(  \sigma\ast\pi\right)  $ holding at $\left(  u_{i},u_{i+1}%
\right)  $ would be in the branch at that stage as well. But they are formulas
of lower complexity than $\sigma\ast\pi$, so by the induction hypothesis,
those formulas must hold in the model which contradicts the complementary
Boolean conditions holding at all links of that chain in the model. Hence
$\left(  u,u^{\prime}\right)  \in\operatorname*{dit}\left(  \sigma\ast
\pi\right)  $ holds in the model.

Suppose $\left(  u,u^{\prime}\right)  :F\left(  \sigma\ast\pi\right)  $ occurs
at some stage $U_{n}$ in the open branch. Then by completeness and the
connective rule for $F\left(  \sigma\ast\pi\right)  $, there is a finite
$u,u^{\prime}$-chain in $U_{n+1}$ (or in $U_{n}$ if no new elements were
introduced) where the formulas for the Boolean conditions for $F\left(
\sigma\ast\pi\right)  $ occur at each link in the chain. But all those
formulas are of lower complexity than $\sigma\ast\pi$ so by the induction
hypothesis, they are true in the model on $U$, which in turn implies that
$\left(  u,u^{\prime}\right)  \in\operatorname*{indit}\left(  \sigma\ast
\pi\right)  $ in that model.

Since the formula $\left(  u_{0},u_{1}\right)  :F\varphi$ occurs in every
branch, the open complete branch supplies a model where $\left(  u_{0}%
,u_{1}\right)  \in\operatorname*{indit}\left(  \varphi\right)  $.
$\blacksquare$

\begin{theorem}
[Completeness theorem for partition tableaus]If $\varphi$ is a partition
tautology, then any completed tableau beginning with $\left(  u_{0}%
,u_{1}\right)  :F\varphi$ must close, and thus every partition tautology is
provable by the tableau method.
\end{theorem}

\noindent Proof: If a completed tableau beginning with $\left(  u_{0}%
,u_{1}\right)  :F\varphi$ had a complete open branch, then by the satisfaction
theorem there would be an interpretation where $\left(  u_{0},u_{1}\right)
\in\operatorname*{indit}\left(  \varphi\right)  $ and thus $\varphi$ is not a
partition tautology. Hence if $\varphi$ is a partition tautology, then any
completed tableau beginning with $\left(  u_{0},u_{1}\right)  :F\varphi$ must
close so that $\varphi$ is a theorem by the tableau method. $\blacksquare$

\section{Concluding remarks}

The view of "propositional" logic as being about subsets goes back to the
beginning of modern logic in Boole and DeMorgan.

\begin{quote}
The algebra of logic has its beginning in 1847, in the publications of Boole
and De Morgan. This concerned itself at first with an algebra or calculus of
classes, to which a similar algebra of relations was later added. Though it
was foreshadowed in Boole's treatment of "Secondary Propositions," a true
propositional calculus perhaps first appeared from this point of view in the
work of Hugh MacColl, beginning in 1877. \cite[pp. 155-156.]{church:ML}
\end{quote}

\noindent Indeed, in Boole's treatment of "Secondary Propositions," he noted
that propositions could be substituted for subsets and the same laws would hold.

\begin{quote}
But while the laws and processes of the method remain unchanged, the rule of
interpretation must be adapted to new conditions. Instead of classes of
things, we shall have to substitute propositions, and for the relations of
classes and individuals, we shall have to consider the connexions of
propositions or of events. \cite[p. 162]{boole:lot}
\end{quote}

The key mathematical fact that allowed the specific propositional
interpretation of Boolean logic to eventually overshadow the general subset
interpretation is that the logical operations on subsets can be modeled using
just the subsets $0$ and $1$ of the one element set. As Boole himself emphasized:

\begin{quote}
\textit{We may in fact lay aside the logical interpretation of the symbols in
the given equation; convert them into quantitative symbols, susceptible only
of the values }$0$\textit{\ and }$1$\textit{; perform upon them as such all
the requisite processes of solution; and finally restore them to their logical
interpretation.} \cite[p. 70]{boole:lot}
\end{quote}

In more recent times, the subset interpretation of Boolean logic has been
reemphasized by categorical logic in the topos of sets and by the
Kripke-structure and topological interpretations of intuitionistic
"propositional" logic. Boolean logic works with the subsets of an unstructured
universe set $U$ and those interpretations of intuitionistic logic add
structure to the universe set $U$ to define a topology (e.g., the up-closed
subsets from a partial ordering on the universe set) so that the relevant
subsets for the interpretation are the open subsets. Boolean logic can then be
seen as the special case with the discrete topology on $U$ so that all subsets
are open and the intuitionistic operations reduce to the Boolean ones.

Partition logic, like Boolean logic, starts with an unstructured universe set
$U$ (two or more elements). The subsets of the powerset Boolean algebra
$\mathcal{P}(U)$ and the partitions of the partition algebra $\Pi(U)$ are both
defined simply on the basis of the set $U$. Thus subset logic and partition
logic are at the same mathematical level, and are based on the dual concepts
of subsets and partitions. Partition logic provides a dual semantics for
Boolean logic formulas, a semantics based on the distinctions of partitions
rather than the elements of subsets.

One can go further with the subsets-partitions and elements-distinctions
duality. Probability theory can be seen as a conceptual continuation of subset
logic. Probability theory conceptually starts with the finite case where the
probability is the ratio of the number of elements in a subset ("event") to
the size of the finite universe $U$ ("sample space") of equiprobable outcomes.
This conceptual continuation from subset logic to finite probability theory
was there from the beginning in Boole. Quoting Poisson, Boole defined "the
measure of the probability of an event [as] the ratio of the number of cases
favourable to that event, to the total number of cases favourable and
unfavourable, and all equally possible." \cite[p. 253]{boole:lot}

What arises from the analogous continuation of partition logic? Replacing
elements and subsets with distinctions and partitions yields a \textit{logical
information theory} where the \textit{logical entropy} of a partition is
defined as the ratio of the number of distinctions of the partition to the
size of the finite closure space $U\times U$. The resulting logical
information theory provides a logical-conceptual foundation for Shannon's
information theory \cite{ell:cd}.

Finally, we might speculate about why it has taken so long for partition logic
to be developed. The subset interpretation dates back to the beginning of
modern logic in the mid-nineteenth century. The subset-partition duality is at
least as old as category theory (mid-twentieth century) and in any case has
long been evident in the interplay of subobjects and quotient objects
throughout algebra and in the subset-partition analogies of combinatorial
theory. There seems to be a cluster of reasons for the delay.

From the side of logic, we have already noted that the propositional
interpretation of Boolean logic has all but eclipsed the general subset
interpretation so that most non-category-theoretic treatments of logic give
only the propositional interpretation.

Moreover, the progression from "propositional" logic to "quantification
theory" is based on the specific propositional interpretation of Boolean logic
where the propositions are quantified formulas. Tarski's semantics developed
as model theory has been very successful in applications. \ Model theory
interprets open formulas and atomic relations as subsets of an $n$-fold
product $U^{n}$ of some underlying universe set, and then closed formulas are
propositions which are either true or false. But Lawvere's development of
categorical logic brings out the general setting in the category of sets.
Given a set map $f:V\rightarrow U$ between two universe sets, the two
quantifiers are given as left and right adjoints mapping subsets of $V$ to
subsets of $U$.\footnote{The technical details are not relevant to our point
here since this paper does not deal with "quantifier" morphisms for partition
logic. The categorical logic treatment of the subset quantifiers is covered in
Mac Lane
\citeyear{mac:cwm}%
, Lawvere and Rosebrugh
\citeyear{law:sfm}%
, Awodey
\citeyear{awod:ct}%
, and Mac Lane and Moerdijk
\citeyear{macm:sh}%
.} In the special case of classical quantification theory, quantifying over a
variable in effect takes the set map as the projection $U^{n}\rightarrow
U^{n-1}$ that leaves out the variable so that the subset quantifiers carry
subsets of $U^{n}$ to subsets of $U^{n-1}$. When $n=1$, quantifying over the
single variable is usually interpreted as turning an open single-variable
formula into a closed formula or proposition which is true or false, but the
interpretation in categorical logic is mapping subsets of $U^{1}$ to subsets
of $U^{0}=1$ where the subsets of $1$ behave like the usual propositional
truth values.

The propositional special case has been so important that the general case of
subset logic (not to mention subset quantifier-morphisms) has been rather
neglected. The part has been taken as the whole. The point here is that since
propositions do not have a dual notion of partitions, the idea of a dual logic
of partitions does not arise in the conventional treatment of Boolean logic as
"propositional" logic. Indeed, one might wager that the dual interplay of
subsets and quotient sets is sufficiently well-known so that \textit{if}
Boolean logic had been commonly understood as subset logic, then partition
logic would not be far behind. In that sense, the hegemony of the
propositional interpretation of subset logic seems to be the principal reason
for the late development of partition logic.

From the partition side, it has long been known that partitions form a lattice
just like subsets. But the "lattice of partitions" was traditionally defined
"upside down" as (isomorphic to) the lattice of equivalence relations rather
than its opposite. However, the element-distinction duality makes it clear
that the lattice of partitions should use the partial ordering given by the
set of distinctions (dit set) of a partition rather than its set of
indistinctions (just as the lattice of subsets uses the partial ordering given
by the set of elements of a subset rather than its set of non-elements). This
is what allowed the direct comparison of formulas in subset, intuitionistic,
and partition logic as well as the proof-theoretic parallels between the
tableaus for the three logics.

Another reason is that (at least to our knowledge) the implication, nand and
other new binary operations on partitions (aside from the join and meet) have
not been previously studied. In a recent paper in a commemorative volume for
Gian-Carlo Rota, the three authors remark that in spite of the importance of
equivalence relations, only the operations of join and meet have been studied.

\begin{quotation}
Equivalence relations are so ubiquitous in everyday life that we often forget
about their proactive existence. Much is still unknown about equivalence
relations. Were this situation remedied, the theory of equivalence relations
could initiate a chain reaction generating new insights and discoveries in
many fields dependent upon it.

This paper springs from a simple acknowledgement: the only operations on the
family of equivalence relations fully studied, understood and deployed are the
binary join $\vee$ and meet $\wedge$ operations. \cite[p. 445]{bmp:eqrel}
\end{quotation}

Yet the new operations, particularly the implication, are crucial to the whole
development. The only partition tautologies with only lattice operations are
trivialities such as $1$ and $1\vee\pi$. Without the non-lattice operations,
one can always study identities in the partition lattice such as $\pi
\preceq\pi\vee\sigma$ (which corresponds to the tautology $\pi\Rightarrow
\pi\vee\sigma$). But it has been shown \cite{whitman:lat} that partition
lattices are so versatile that any formula in the language of lattices (i.e.,
without the implication or other non-lattice operations) that is an identity
in all partition lattices (or lattices of equivalence relations) is actually a
general lattice-theoretic identity. Hence the logic taking models in all
partition algebras $\Pi\left(  U\right)  $ only became interesting by moving
beyond the lattice operations on partitions.

Throughout his career, Gian-Carlo Rota emphasized the analogies between the
Boolean lattice of subsets of a set and the lattice of equivalence relations
on a set. Partition logic, with the heavy emphasis on the analogies with
subset logic, should be seen as a continuation of that Rota program. The
closest earlier work in the vein of the partition logic tableaus was indeed by
Rota and colleagues [\cite{rota:logiccers} and \cite{Haim:LinLog}], but it
used the lattice of equivalence relations and did not define the partition
implication (which would be the difference operation on equivalence relations)
or other non-lattice operations. It was restricted to the important class of
commuting equivalence relations \cite{dub:cer} where identities hold which are
not general lattice-theoretic identities.

In sum, the subset interpretation of Boolean logic (so the subset-partition
duality would come into play), the turning of the lattice of partitions right
side up, and the introduction of the non-lattice operations (particularly the
implication) were all important in the development of partition logic.

\bibliographystyle{chicago}
\bibliography{Partition-Logic}

\end{document}